\documentclass[11pt]{article}
\usepackage{graphicx,epstopdf,subfigure,mathtools,mathrsfs, arydshln, amsmath, amssymb} 
\usepackage[font=small,labelfont=bf]{caption}
\usepackage{float}
\usepackage{authblk}
\usepackage[title]{appendix}
\PassOptionsToPackage{usenames,dvipsnames}{xcolor}
\usepackage[usenames,dvipsnames]{xcolor}
\usepackage[margin=1in]{geometry}
\usepackage[normalem]{ulem}

\usepackage{amsfonts}

\newcommand{\red}[1]{\color{red}#1\normalcolor}
\newcommand{\delete}[1]{}
\newcommand{\change}[1]{\color{black}#1\normalcolor}
\newcommand{\rev}[1]{\color{black}#1\normalcolor}

\newcommand{\V}[1]{\boldsymbol{#1}}                 
\newcommand{\M}[1]{\boldsymbol{#1}}
\newcommand{\Lop}[1]{\boldsymbol {\mathcal{#1}}}
\newcommand{\widebar}[1]{%
   \hbox{%
     \vbox{%
       \hrule height 0.5pt 
       \kern0.5ex
       \hbox{%
         \kern-0.1em
         \ensuremath{#1}%
         \kern-0.1em
       }%
     }%
   }%
}

\global\long\def\DO{\V{D}^{(1)}}
\global\long\def\DT{\V{D}^{(2)}}

\global\long\def\norm#1{\left\Vert #1\right\Vert }
\global\long\def\tt#1{#1_\text{f}}
\global\long\def\tr#1{#1_\text{n}}
\global\long\def\rt#1{#1_\text{f}}
\global\long\def\rr#1{#1_\text{n}}
\global\long\def\Xs{\V{\tau}}
\global\long\def\DR{\V{D}^{(3)}}
\global\long\def\Slet#1{\mathbb{S}\left(#1\right)}
\global\long\def\Dlet#1{\mathbb{D}\left(#1\right)}

\global\long\def\EPMI{\frac{1}{8\pi\mu}}
\global\long\def\tdisc#1#2{#1^{#2}}
\global\long\def\grad{\nabla}
\global\long\def\Rhat{\widehat{\V{R}}}
\global\long\def\twmod{\gamma}
\global\long\def\follows{\rightarrow}
\global\long\def\dt#1{\partial_t #1}
\global\long\def\ds#1{\partial_s #1}
\global\long\def\ind#1#2{#1^{(#2)}}
\global\long\def\Eul#1{#1_E}
\global\long\def\eps{{\hat{a}}}
\global\long\def\epsRS{\hat{\epsilon}}
\global\long\def\rc{a}
\global\long\def\epsc{\epsilon}

\global\long\def\dOm{\Delta \V{\Psi}}
\global\long\def\dOmpar{\Delta \Omega^\parallel}
\global\long\def\dF{\Delta \V{\Lambda}}

\global\long\def\dU{\Delta \V{U}}
\global\long\def\Mtt{\M{M}_\text{tt}}
\global\long\def\Mtr{\M{M}_\text{tr}}
\global\long\def\Mrt{\M{M}_\text{rt}}
\global\long\def\Mrr{\M{M}_\text{rr}}
\global\long\def\Mbtt{\mathbb{M}_\text{tt}}
\global\long\def\Mbtr{\mathbb{M}_\text{tr}}
\global\long\def\Mbrt{\mathbb{M}_\text{rt}}
\global\long\def\Mbrr{\mathbb{M}_\text{rr}}
\global\long\def\Mctt{\Lop{M}_\text{tt}}
\global\long\def\Mctr{\Lop{M}_\text{tr}}
\global\long\def\Mcrt{\Lop{M}_\text{rt}}
\global\long\def\Mcrr{\Lop{M}_\text{rr}}

\usepackage{hyperref}
\hypersetup{
    colorlinks=false,
    pdfborder={0 0 0},
}

\title{The hydrodynamics of a twisting, bending, inextensible fiber in Stokes flow}
\author[1]{Ondrej Maxian}
\author[1]{Brennan Sprinkle}
\author[1]{Charles S. Peskin}
\author[1]{Aleksandar Donev}
\affil[1]{Courant Institute of Mathematical Sciences, New York University, New York, NY 10012}

\begin{document}
\maketitle

\begin{abstract}
In swimming microorganisms and the cell cytoskeleton, inextensible fibers resist bending and twisting, and interact with the surrounding fluid to cause or resist large-scale fluid motion. In this paper, we develop a novel numerical method for the simulation of cylindrical fibers by extending our previous work on inextensible bending fibers [Maxian et al., \textit{Phys. Rev. Fluids} \textbf{6} (1), 014102] to fibers with twist elasticity. In our ``Euler'' model, twist is a scalar function that measures the deviation of the fiber cross section relative to a twist-free frame, the fiber exerts only torque parallel to the centerline on the fluid, and the perpendicular components of the rotational fluid velocity are discarded in favor of the translational velocity. In the first part of this paper, we justify this model by comparing it to another commonly-used ``Kirchhoff'' formulation where the fiber exerts both perpendicular and parallel torque on the fluid, and the perpendicular angular fluid velocity is required to be consistent with the translational fluid velocity. Through asymptotics and numerics, we show that the perpendicular torques in the Kirchhoff model are small, thereby establishing the asymptotic equivalence of the two models for slender fibers. We then develop a spectral numerical method for the hydrodynamics of the Euler model. We define hydrodynamic mobility operators using integrals of the Rotne-Prager-Yamakawa tensor, and evaluate these integrals through a novel slender-body quadrature, which requires on the order of 10 points along a smooth fiber to obtain several digits of accuracy. We demonstrate that this choice of mobility removes the unphysical negative eigenvalues in the translation-translation mobility associated with asymptotic slender body theories, and ensures strong convergence of the fiber velocity and weak convergence of the fiber constraint forces. We pair the spatial discretization with a semi-implicit temporal integrator to confirm the negligible contribution of twist elasticity to the relaxation dynamics of a bent fiber. We also study the instability of a twirling fiber, and demonstrate that rotation-translation and translation-rotation coupling change the critical twirling frequency at which a whirling instability onsets by $\sim$20\%, which explains some (but not all) of the deviation from theory observed experimentally. \change{When twirling is unstable, we show that the dynamics transition to a steady overwhirling state, and quantify the amplitude and frequency of the resulting steady crankshafting motion.}
\end{abstract}

\section{Introduction}
The focus of this paper is on the simulation of cylindrical, slender, inextensible fibers immersed in a Stokes (zero Reynolds number) fluid, as found in cytoskeletal fibers \cite{alberts}, bacterial flagella and cilia \cite{brennen1977fluid, lauga2009hydrodynamics}, and fiber composites \cite{agarwal2017analysis}. Our particular goal is to study numerically how \emph{twist elasticity} contributes to the dynamics of these filaments. Previous work has shown that twist elasticity controls the steady state structure and energy storage in actin filament bundles \cite{claessens2008helical, ma2018structural}, and that a spinning filament that resists twist is unstable to perturbations, which can lead to whirling and swimming motions \cite{wolgemuth2000twirling, powers2010dynamics}.

In modeling twisted fibers immersed in a fluid, there are two largely independent choices to be made: the elastic rod model and the hydrodynamic mobility (force-velocity relationship). With regard to the first of these, most studies that explicitly (or fully) couple the rod to the surrounding fluid use some form of the Kirchhoff rod model \cite{dill1992kirchhoff}, in which the inextensible and unshearable rod resists both bending (along two principal axes) and twisting (along the third axis parallel to the fiber centerline). As we discuss, there is an ambiguity in the equilibrium equations of the Kirchhoff rod in which the torque density perpendicular to the fiber centerline can be recast into a force density; many studies refer to this parallel-torque-only version as ``the Kirchhoff rod model'' \cite[Sec.~2.3]{vogel2012bacterial}. In either case, the Kirchhoff model is a special case of the Cosserat model \cite{rubin2013cosserat, gargslender}, which also allows for extensibility and shear of the rod, and has recently been used to model cytoskeletal fibers without dynamics \cite{floyd2021stretching}.

For hydrodynamic coupling, the numerical methods used to simulate slender fibers fall into two broad categories: methods in which the fibers are represented by a discrete set of regularized singularities  or ``blobs,'' and asymptotic methods, in which the true three-dimensional dynamics are approximated in the limit of very slender fibers. The blob-based class of methods includes the immersed boundary method \cite{peskin1972flow, peskin2002acta, kallemov2016immersed}, the method of regularized Stokeslets \cite{cortez2001method, cortez2005method}, and the force coupling method \cite{fcm03}. \rev{In some of these, at least for an unbounded domain, the Stokes equations can be solved analytically for the chosen regularization function, thereby analytically eliminating the fluid equations}. Asymptotic methods are also referred to as slender body theories (SBTs), and were originally developed as analytical tools \cite{krub, johnson, gotz2001interactions}, but later used in numerical simulations \cite{shelley2000stokesian, ts04, ehssan17, maxian2021integral}. Numerical experiments have shown that regularized singularities can accurately reproduce the results of slender body theory \cite{ttbring08} or analytical solutions to the Stokes equations in the infinite cylinder limit \cite[Sec.~5.2]{balboa2017hydrodynamics}, thus justifying that the two hydrodynamic approaches are equivalent. We recently showed \cite[Appendix~A]{maxian2021integral} that the continuum limit of infinitely many blobs distributed on the fiber centerline corresponds to Keller-Rubinow slender body theory \cite{krub}, if the Rotne-Prager-Yamakawa (RPY) \cite{RPYOG} regularization radius $\eps$ is an $\mathcal{O}(1)$ multiple of $\rc$, the true cylinder radius, $\eps = \rc e^{3/2}/4\approx 1.12 \rc$; Cortez and Nicholas obtained a related result for regularized Stokeslets \cite{cortez2012slender}. Throughout this paper, we will use $\epsRS=\eps/L$ for the aspect ratio of the regularized cylinder of length $L$, which is $\approx 12\%$ larger than the true cylinder aspect ratio $\epsc = \rc/L$.

The numerical use of slender body theory came after regularized singularity methods to solve a problem in the latter: when the fiber is discretized as a chain of blobs, in order to accurately \rev{hydrodynamically} model the fiber as a continuous curve, and not as a series of disconnected beads, the spacing between the blobs must be on the order of the fiber radius. This means that the number of degrees of freedom scales as $L/\eps=\epsRS^{-1}$. Slender body theory circumvents this problem, because the fiber is modeled as a continuous and smooth curve, with arclength parameterization $\V{X}(s,t)$. The \rev{hydrodynamics} can be reformulated in terms of integrals along the curve, and spectral methods can be used to achieve high accuracy \cite{ehssan17, maxian2021integral}. Unfortunately, while SBT methods are the most efficient way to simulate the dynamics of smooth slender fibers, they have never been fully extended to include twist elasticity, \rev{as discussed in \cite{TwistSBT}}.\delete{Although the original SBT of Keller and Rubinow \cite{krub} did include fiber twist and dilation, their matched asymptotic procedure yields an expression for the fiber velocity that depends on the intermediate point in the fluid where the matching is done. } As a result, a disconnect has emerged in the literature in which immersed boundary methods \cite{lim2008dynamics}, regularized Stokeslets \cite{inexRS, olson2013modeling}, and the force coupling method \cite{keavRPY}, have accounted for the full hydrodynamics (including rotation-translation (rot-trans) and translation-rotation (trans-rot) coupling) of twisting fibers, while methods that use slender body theory \cite{powers2010dynamics, wada2006non, nguyen2018impacts} have (at best) included nonlocal hydrodynamics only in the translation-translation (trans-trans) mobility, with only local hydrodynamics in the rotation-rotation (rot-rot) mobility, and no trans-rot/rot-trans coupling. 

A related, and perhaps more important, ambiguity is how the rod elasticity couples to the fluid mechanics, for which there are again two unreconciled choices in the literature. The main question is: If the rod moves with the local fluid velocity, does the cross section also rotate with the local angular fluid velocity? If it does, then for the rod to be unshearable, the fluid velocity has to be constrained so that the translational velocity $\V{U}(s)$ and angular velocity $\V{\Psi}(s)$ are consistent with each other, where $s$ denotes arclength. Specifically, there are two possible velocities of the tangent vector $\Xs=\ds{\V{X}}$:  $\dt{\Xs} = \ds{\V{U}}$, and $\dt{ \Xs} = \V{\Psi} \times \Xs$. If we require that the two match, i.e., if we impose the constraint $\ds{\V{U}}=\V{\Psi} \times \Xs$ , we get what we refer to as the \emph{Kirchhoff model} for historical reasons \cite{keavRPY, inexRS, westwood2021coordinated}. To avoid constraints, Lim et al.\ \cite{lim2008dynamics} introduced the \emph{generalized} Kirchhoff rod model, in which the cross section is shearable and the fiber is extensible, both with an associated penalty. This model has been used frequently with the original IB method and with regularized Stokeslets \cite{olson2013modeling, ishimoto2018elasto}, but its consistency with the constrained Kirchhoff formulation as all penalty parameters go to infinity has not yet been established. Of course, imposing the constraint $\ds{\V{U}}=\V{\Psi} \times \Xs$ requires knowledge of the full rotational and translational fluid velocity, which so far has only been possible in regularized singularity methods, since SBT methods have not yet properly accounted for rot-trans coupling.

In the second class of fluid coupling methods, the redundancy in the possible motions is treated not by enforcing a constraint, but by ignoring the rotational fluid velocity in perpendicular directions. The rod then has only four possible motions: translational velocity $\V{U}(s)$ and rotational velocity $\Omega^\parallel(s)$ about its axis. In this model, the tangent vector $\Xs$ is updated by the translational velocity, and then a local triad or some other measure of twist is rotated around $\Xs$ by the velocity $\Omega^\parallel$. This decouples solving for the translational and rotational velocity, simplifying the problem and reducing the size of the linear system to be solved. This model, which we shall refer to as the \emph{Euler model}, is presented in the review article by Powers \cite{powers2010dynamics}, and has been used in the simulation of both discrete \cite{nguyen2018impacts, reichert2006hydrodynamic, wada2006non, wada2007model,  huang2021numerical} and continuous \cite{man2017bundling} elastic rods. All of these methods use an SBT-style mobility, in which nonlocal hydrodynamics is included in the translation-translation (trans-trans) coupling (to give the translational velocity $\V{U}$), but only local drag is included for rotation-rotation (rot-rot, to give the parallel rotational velocity $\Omega^\parallel$). Indeed, a perpendicular angular velocity is not available because no slender body theory exists for translation-rotation (trans-rot) or rotation-translation (rot-trans) coupling, so only parallel rotational velocity can be used. 

This paper has two parts, each of which addresses one of the two main ambiguities in the literature. In part one, we compare the two approaches for coupling fluid mechanics and rod elasticity, which we refer to as the Euler and Kirchhoff models. For the fluid mechanics, in Section\ \ref{sec:two} we define the four components of the mobility operator (trans-trans, trans-rot, rot-trans, rot-rot) in terms of integrals of the Rotne-Prager-Yamakawa RPY kernel \cite{wajnryb2013generalization}, which is a popular, physics-inspired form of a regularized singularity. To compare the two models of fluid coupling, we use both asymptotic arguments and numerical tests to show that the Kirchhoff and Euler models give the same translational velocity $\V{U}$ a distance $\mathcal{O}(1)$ away from the fiber endpoints as the slenderness $\epsc\rightarrow 0$. The key idea is that the perpendicular torque in the Kirchhoff model, which enforces the constraint $\ds{\V{U}}=\V{\Psi} \times \Xs$, can be as small as $\mathcal{O}\left(\epsc^2\right)$ and still generate a total angular velocity $\V{\Psi}$ that is consistent with $\V{U}$. As a result, these Lagrange multipliers have strength of order $\epsc^2$, and therefore have a minimal impact on $\ds{\V{U}}$, so they can be dropped in the Euler model without changing the resulting centerline velocity $\V{U}$ to $\mathcal{O}\left(\epsc^2\right)$. Thus the first part of this paper establishes the asymptotic equivalence of the Euler and Kirchhoff models as $\epsc \rightarrow 0$. We will use the Euler model, since that choice makes it simpler to add twist to our previous model of inextensible bending fibers \cite{maxian2021integral} as an extra degree of freedom. \rev{The key equations are summarized in Section\ \ref{sec:EulEqns}.}

\rev{In the second part of this paper and a companion paper \cite{TwistSBT}, we address the gap in the SBT literature with regard to the hydrodynamic mobility and rot-trans coupling. In \cite{TwistSBT}, we develop a slender body theory for twisting which is based on evaluating the RPY integrals asymptotically in the (regularized) slenderness $\epsRS$.}\delete{While this theory is not an asymptotic treatment of the fully three-dimensional problem (like those of Keller and Rubinow \cite{krub}, Koens and Lauga \cite{koens2018boundary, gargslender}, or Mori et. al \cite{mori2018theoretical}), it does provide a treatment of hydrodynamic mobility that is in the spirit of SBT and captures the essential physics and leading-order behavior (up to constants).} Unfortunately, our asymptotic evaluation of the RPY integrals is still plagued by a longstanding problem in SBT-style methods: the ill-posedness of the trans-trans mobility, which has negative eigenvalues corresponding to force densities with lengthscales less than $\epsc$ \cite{gotz2001interactions, ts04, mori2020accuracy}. This is problematic in constrained problems, in which the mobility appears in a linear system to be solved for the constraint force. To avoid this problem with the asymptotic mobilities, in this paper we will instead use integrals of the RPY kernels along the centerline to define the action of the mobility operators (with the exception of rot-rot coupling). To evaluate these near-singular integrals to high accuracy, in Section\ \ref{sec:specEul} and Appendix\ \ref{sec:quads} we develop ``slender body'' quadrature schemes which give three digits of accuracy for smooth fibers with on the order of 10 points discretizing the fiber centerline. This represents a significant improvement over existing regularized singularity approaches, which require $1/\epsRS$ points to resolve the \rev{hydrodynamics}, and, unlike asymptotic methods, our approach reduces the required degrees of freedom \emph{without} introducing negative eigenvalues in the mobility operators. Coupling our quadratures with the Euler model gives a spectral numerical method for the twisting and bending dynamics of an inextensible, slender, smooth filament with as few degrees of freedom as possible. In particular, in Section\ \ref{sec:staticconv} \rev{we consider the relaxation of bent and twisted filaments with varying curvature, and show that our spectral discretization gives a more accurate velocity with $\rev{32-40}$ points on the fiber centerline than a standard second-order, blob-based method with $1/\epsRS = 1000$ points}. We show empirically that the constraint force enforcing inextensibility converges weakly, even though it is not smooth near the fiber endpoints.

In Section\ \ref{sec:dyn} of this paper, we extend our numerical method to dynamic problems by introducing a (constrained) backward Euler temporal integration scheme. We apply this numerical method to two problems involving twisting fibers: a relaxing bent fiber, where we explore the feedback of twist on the fiber dynamics, and a twirling fiber, where we look at the critical frequency at which the fiber starts to whirl and \change{then transitions to a stable ``overwhirling'' state} \cite{wolgemuth2000twirling, lim2004simulations, wada2006non, lee2014nonlinear}. For fiber relaxation, we show that the assumption of twist being in quasi-equilibrium \cite{bergou2008discrete, man2017bundling} is satisfied for materials in which the twist and bending moduli are comparable, as is the case in most materials (for inextensible isotropic elastic cylinders, the ratio is 2/3 \cite{nguyen2018impacts}). For a twirling fiber, we show how the systematic neglect of rot-trans coupling in prior studies based on the Euler model \cite{wolgemuth2000twirling, wada2006non} leads to an overestimation in the critical frequency of 20\%. Our results explain half of the deviation between prior theoretical/computational studies and what is observed experimentally \cite{bruss2019twirling}.

\section{Model formulation and governing equations \label{sec:two}}
In this section, we formulate the governing equations for both the Kirchhoff and Euler models in a viscous fluid. For brevity of notation, we omit explicit time dependence; it is implied that all functions of space also depend on time $t$. We also omit dependence on arclength $s$ whenever there is no ambiguity. We start from the classical immersed boundary equations \cite{peskin2002acta}, and take a continuum limit to model the rod as a line of regularized singularities. After introducing the hydrodynamic model in Section\ \ref{sec:gmob}, we discuss the two models of inextensible, bending, twisting fibers. This is necessarily preceded by a brief introduction to inextensibility and twist elasticity in Section\ \ref{sec:fibmech}, where we introduce the Bishop frame \cite{bergou2008discrete}. This reference frame allows for a simple representation of twist as a scalar angle $\theta(s)$, instead of the full orthogonal triad used in previous formulations \cite{lim2008dynamics, olson2013modeling}. We use this scalar angle as our measure of twist in Section\ \ref{sec:krod}, where we introduce the Kirchhoff rod model. A key feature of this model is the consistency condition between the angular velocity $\V{\Psi}(s)$ and translational velocity $\V{U}(s)$ of the cross section. This requires a Lagrange multiplier force density with \emph{three} independent components for the three directions of the velocity. 

We compare the Kirchhoff formulation to the Euler model, which we present in Section\ \ref{sec:euler}. In the Euler model, we solve for the translational velocity, which we use to evolve the tangent vectors \cite{maxian2021integral}, then compute parallel rotational velocity to evolve the twist in a second, independent step. The Euler formulation completely disregards the perpendicular components $\V{\Psi}^\perp$ of the fluid angular velocity that could in principle be obtained from our hydrodynamic model. Since the Euler model does not enforce consistency of $\V{\Psi}^\perp$ with the translational velocity $\V{U}$, it requires only a scalar tension as a Lagrange multiplier. \rev{In Section\ \ref{sec:compareMeth}, we show that the two discarded perpendicular components of the Lagrange multiplier are small, thereby establishing the equivalence of the Kirchhoff and Euler models as $\epsc \rightarrow 0$.}

\subsection{The hydrodynamic mobility \label{sec:gmob}}
In immersed boundary (IB) methods \cite{peskin2002acta, cortez2001method}, the fiber is discretized by a series of markers or ``blobs'' of radius $\eps$, each of which exerts a regularized point force on the fluid. In grid-free approaches, the regularization is done over a spherically-symmetric function $\delta_\eps(\V{r}) \equiv \delta_\eps(r)/(4\pi\eps^2)$ of width $\sim \eps$. Following \cite[Appendix~A]{maxian2021integral}, we will pass to a continuum limit, so that instead of a sum of blobs, we obtain an integral over a continuous curve. For translational velocity from force, we showed in \cite[Appendix~A]{maxian2021integral} that the continuum limit of regularized singularties can be made equivalent to slender body theory for a cylinder of radius $\rc = \epsc L$ \cite{krub,johnson} via a judicious choice of the $\mathcal{O}(1)$ constant ratio $\eps/\rc$, specifically
\begin{equation}
\label{eq:ahat}
\eps = \rc \frac{e^{3/2}}{4} \approx 1.1204 \rc.
\end{equation}

Introducing the fiber as a curve $\V{X}(s)$ with tangent vector $\Xs(s)=\ds{\V{X}}(s)$, force density (units force per length) $\V{f}(s)$, and torque density $\V{n}(s)$, the fluid velocity $\V{u}(\V{x})$ and pressure $\pi(\V{x})$ are given by the solution of the Stokes equations, $\nabla \cdot \V{u}=0$ and 
\begin{equation}
\label{eq:Stokes}
\grad\pi(\V{x})=\mu \grad^{2}\V u(\V{x})+ \int_0^L \left(\left(\V{f}(s) + \frac{\nabla}{2} \times \V{n}(s)\right)\delta_{\eps}\left(\V{x}-\V X(s)\right)\right) \, ds,
\end{equation}
where $\mu$ is the fluid viscosity. In this paper we consider an unbounded domain with the fluid at rest at infinity, but other boundary conditions can be used. Because the Stokes equations\ \eqref{eq:Stokes} are linear, the solution of\ \eqref{eq:Stokes} can be written analytically as a convolution with the Green's function, also called the Stokeslet or Oseen tensor in free space,
\begin{equation}
\label{eq:Slet}
\Slet{\V{x},\V{y}}=\EPMI \dfrac{\M{I}+\Rhat \Rhat}{R},
\end{equation}
where $\V{R}=\V{x}-\V{y}$, $R=\norm{\V{R}}$, $\Rhat = \V{R}/R$, and $\Rhat \Rhat$ denotes the outer product of $\Rhat$ with itself. In this context, it is also useful to define the doublet (Laplacian of the Stokeslet), which is given by 
\begin{equation}
\label{eq:Dlet}
\Dlet{\V{x},\V{y}}=\EPMI \dfrac{\M{I}-3\Rhat \Rhat}{R^3}.
\end{equation}
The translational and rotational velocities of the fiber centerline can then be obtained from the fluid velocity $\V{u}(\V{x})$, as we detail next. 

\subsubsection{Translational velocity}
In most blob-based methods,\footnote{An exception is the method of regularized Stokeslets, which takes $\V{U}(s)=\V{u}\left(\V{X}(s)\right)$.} the translational velocity of the fiber centerline, which we denote by $\V{U}(s)=\dt{\V{X}}(s)$, is a convolution of $\V{u}(\V{x})$ with the blob function $\delta_\eps$, 
\begin{equation}
\label{eq:Uavg}
\V U \left(s\right)= \int \V{u}(\V{x}) \, \delta_{\eps}\left(\V X(s)-\V x\right) \, d\V{x}.
\end{equation}
To write this velocity as a convolution with the Stokeslet $\mathbb{S}$, we first define the kernels
\begin{align}
\label{eq:Ukernels}
\Mbtt\left(\V{x},\V{y}\right) & =  \int \delta_{\eps}\left(\V x-\V z\right)  \int \Slet{\V{z},\V{w}}\delta_{\eps}\left(\V w-\V y\right) \, d\V{w} \, d\V{z}\\ \nonumber
\Mbtr\left(\V{x},\V{y}\right) & =  \int \delta_{\eps}\left(\V x-\V z\right)  \int \Slet{\V{z},\V{w}} \frac{\nabla_{\V{w}}}{2} \times\delta_{\eps}\left(\V w-\V y\right) \, d\V{w} \, d\V{z},
\end{align}
which express the locally-averaged velocity\ \eqref{eq:Uavg} at point $\V{x}$ due to a regularized force/torque at point $\V{y}$. Using superposition, the translational fiber velocity\ \eqref{eq:Uavg} can be written as the integral
\begin{align}
\label{eq:UIBdef}
\V{U}(s)& = \int_0^L \left(\Mbtt\left(\V{X}(s),\V{X}(s^\prime)\right)\V{f}(s^\prime)+\Mbtr\left(\V{X}(s),\V{X}(s^\prime)\right)\V{n}(s^\prime)\right) \, ds^\prime  \\ \nonumber
& :=\tt{\V{U}}(s)+\tr{\V{U}}(s) \\ \nonumber
&:=\left(\Mctt \V{f} \right)(s)+\left(\Mctr \V{n} \right)(s),
\end{align}
where in the process we split the velocity $\V{U}$ into a component due to force, $\tt{\V{U}}$, and a component due to torque, $\tr{\V{U}}$, and define corresponding mobility (linear) operators $\Mctt\left[\V{X}\right]$ and $\Mctr \left[\V{X}\right]$.

In the classical IB method \cite{peskin2002acta}, the fluid equations\ \eqref{eq:Stokes} are solved numerically on a Cartesian grid, and $\delta_\eps$ is a regularized Delta function that communicates between the fibers and the fluid grid. The double convolution\ \eqref{eq:UIBdef} is therefore done numerically (as is the evaluation of the Stokeslet), and any analytical formulas for them are only approximate \cite[Sec.~IV]{kallemov2016immersed}. When the domain is unbounded, however, for some forms of $\delta_\eps$ we can bypass the fluid solver entirely and write exact formulas for the $3 \times 3$ mobility tensors $\Mbtt$ and $\Mbtr$. This is the case for the Rotne-Prager-Yamakawa (RPY) tensor \cite{RPYOG, wajnryb2013generalization}, for which the regularized delta function is a surface delta function on an $\eps$-radius sphere, $\delta_\eps(r)=\delta(r-\eps)$. In Appendix\ \ref{sec:rpysing}, we give the analytical formulas for the RPY mobility tensors, which we will use throughout this paper.

\subsubsection{Rotational velocity}
In exactly the same way, we define the angular velocity of the fiber centerline by
\begin{equation}
\label{eq:OmIBdef}
\V{\Psi}\left(s\right)=\int \delta_{\eps}\left(\V X(s)-\V x\right)\frac{\grad}{2}\times\V u\left(\V x\right) \, d\V{x},
\end{equation}
which can be written as an integral involving the two kernels
\begin{flalign}
\label{eq:Psikernels}
&&\Mbrt\left(\V{x},\V{y}\right) & =  \int \delta_{\eps}\left(\V x-\V z\right)  \frac{\nabla_{\V{z}}}{2} \times \int \Slet{\V{z},\V{w}}\delta_{\eps}\left(\V w-\V y\right) \, d\V{w} \, d\V{z}=\Mbtr^T\left(\V{y},\V{x}\right)\\ \nonumber
&&\Mbrr\left(\V{x},\V{y}\right) & =  \int \delta_{\eps}\left(\V x-\V z\right)  \frac{\nabla_{\V{z}}}{2} \times \int \Slet{\V{z},\V{w}} \frac{\nabla_{\V{w}}}{2} \times\delta_{\eps}\left(\V w-\V y\right) \, d\V{w} \, d\V{z},\\ 
\label{eq:PsiMobs}
\text{Specifically,} && \V{\Psi}(s)& = \int_0^L \left(\Mbrt\left(\V{X}(s),\V{X}(s^\prime)\right)\V{f}(s^\prime)+\Mbrr\left(\V{X}(s),\V{X}(s^\prime)\right)\V{n}(s^\prime)\right) \, ds^\prime  \\ \nonumber
&& & :=\tt{\V{\Psi}}(s)+\tr{\V{\Psi}}(s) \\ \nonumber
&& &:=\left(\Mcrt \V{f} \right)(s)+\left(\Mcrr \V{n} \right)(s).
\end{flalign}
As for translational velocity, we split $\V{\Psi}$ into the contributions from force $\rt{\V{\Psi}}$ and from torque $\rr{\V{\Psi}}$, and define continuum mobility operators $\Mcrt\left[\V{X}\right]$ and $\Mcrr\left[\V{X}\right]$ associated with each.  In Appendix\ \ref{sec:rpysing}, we give the analytical formulas for the RPY mobility tensors $\Mbrt$ and $\Mbrr$.

\subsubsection{Slender-body asymptotics for rot-rot mobility\label{sec:RPYasymp}}
Because the mobility integrals in\ \eqref{eq:UIBdef} and\ \eqref{eq:PsiMobs} are nearly singular, they are difficult to evaluate numerically, and so it is appealing to look for an approach in which they can be evaluated asymptotically in $\epsRS=\eps/L \ll 1$, \rev{as we do in \cite{TwistSBT}}. \delete{The asymptotic evaluation splits the mobility into a local drag part, which gives the velocity at $s$ from force/torque $\mathcal{O}(\eps)$ away from $s$, and a nonlocal finite part integral which gives the velocity due to force/torque on the rest of the fiber.}However, as we discuss in the Introduction, asymptotic evaluation brings its own set of problems, since it introduces negative eigenvalues in the trans-trans mobility operator and makes its inversion ill-posed. Because of this, our numerical method in Section\ \ref{sec:specEul} uses slender-body \emph{quadratures} to evaluate the action of the three mobilities $\Mctt,\Mctr$, and $\Mcrt$ numerically, rather than asymptotically. 

The rot-rot mobility is different from the other three because the local drag part of the mobility is more dominant. In the case of the trans-trans, rot-trans, and trans-rot mobilities, the local pieces contribute $\log{\epsRS}$ to the mobility, and the rest of the fiber contributes $\mathcal{O}(1)$. In the rot-rot mobility, the local pieces contribute $\epsRS^{-2}$, and the rest of the fiber contributes $\mathcal{O}(1)$. Because the local drag formula for the rot-rot mobility is sufficiently accurate for our purposes (see Fig.\ \ref{fig:QuadConv}), we will use it in place of the full integral\ \eqref{eq:PsiMobs} for the rot-rot mobility $\Mcrr$. The local drag formula we use is 
\begin{gather}
\label{eq:OmMid}
\rr{\V{\Psi}}(s) =\left(p_I(s) \M{I}+p_\tau(s) \M{\Xs}(s)\Xs(s)\right)\V{n}(s),
\end{gather}
where $p_I$ and $p_\tau$ are derived using matched asymptotics up to and including the fiber endpoints \rev{in \cite[Appendix~C]{TwistSBT}}.

In the Euler model (see Section\ \ref{sec:euler}), we are interested in the parallel rotational velocity $\rr{\Psi}^\parallel$ due to parallel torque $\V{n}(s)=n^\parallel(s) \Xs(s)$, for which we use the simplified equation
\begin{equation}
\label{eq:UinRRpar}
8\pi\mu  \rr{\Psi}^\parallel(s) =(p_I+p_\tau)n^\parallel(s)\approx \frac{9n^\parallel(s)}{4 \eps^2},
\end{equation}
in the fiber interior, \rev{$2\eps \leq s \leq L-2\eps$ (expressions for the endpoints are in \cite[Appendix~C]{TwistSBT})}. \delete{since it is obtained by dropping the $\mathcal{O}(\eps^2)$ terms from the $2\eps < s < L-2\eps$ terms in the expression\ \eqref{eq:EPRR}. In our numerical method, we will use the full expressions for $p_I$ and $p_\tau$ given in\ \eqref{eq:EPRR}.} \rev{In this paper, we use\ \eqref{eq:ahat} to set the RPY regularization radius as $\eps = \epsc e^{3/2}/4$. Substituting into\ \eqref{eq:UinRRpar} gives $8\pi\mu  \rr{\Psi}^\parallel \approx 1.79n^\parallel/\rc^2$, which is approximately 10\% different from the formula for the rotational drag on a radius $\rc$ cylinder, $8\pi\mu  \rr{\Psi}^\parallel = 2n^\parallel/\rc^2$ \cite[Eq.~(62)]{powers2010dynamics}.}

\subsection{Fiber mechanics \label{sec:fibmech}}
We simulate inextensible fibers that resist bending and twisting. In this section, we discuss pieces of the fiber mechanics that are common to both the Kirchhoff and Euler models, in particular twist elasticity and inextensibility. The way bending elasticity appears in the equations depends on the underlying model of the fiber, and so we defer our discussion of bending to Sections\ \ref{sec:krod} and\ \ref{sec:euler}.

\subsubsection{Twist elasticity \label{sec:twistint}}
In previous immersed-boundary formulations of twist \cite{lim2008dynamics, olson2013modeling}, the rod is defined by an orthonormal material frame $\left(\DO(s), \DT(s), \DR(s)\right)$. In the standard model for slender fibers, it is assumed that the fibers are unshearable, so that cross sections of the fiber maintain their shape (e.g., circular) and remain perpendicular to the tangent vector. This implies the constraint
\begin{equation}
\label{eq:DrXs}
\DR=\Xs.
\end{equation}
The fiber cross sections are allowed to twist relative to the tangent vector. The twist density $\psi(s)$ is defined as
\begin{equation}
\label{eq:sctwist}
\psi:=\DT \cdot \ds{\DO}
\end{equation}
with boundary conditions
\begin{equation}
\label{eq:twBC}
\psi(0)=\psi(L)=0
\end{equation}
for a free filament \cite{powers2010dynamics}. The use of the material frame is somewhat cumbersome, however, since it requires us to keep track of three separate axes throughout the course of a simulation. 

We will instead use a formulation based on the twist-free Bishop frame \cite{bergou2008discrete}, in which the state variables are the tangent vector $\Xs(s)$ and a degree of twist $\theta(s)$ relative to the Bishop frame. The Bishop frame is an orthonormal frame $\left(\V{b}^{(1)}(s),\V{b}^{(2)}(s), \Xs(s)\right)$ which satisfies
\begin{equation}
\label{eq:twFrBish}
\V{b}^{(2)} \cdot \ds{\V{b}^{(1)}}= -\V{b}^{(1)} \cdot \ds{\V{b}^{(2)}}= 0,
\end{equation}
making it twist free. The Bishop frame along the entire fiber can be constructed by choosing  $\V{b}^{(1)}(0)$ and $\V{b}^{(2)}(0)$ such that $\left(\V{b}^{(1)}(0),\V{b}^{(2)}(0), \Xs(0)\right)$ is an orthonormal frame, and then solving the ordinary differential equation (ODE) \cite{bergou2008discrete}
\begin{equation}
\label{eq:BishODE}
\ds{\V{b}}^{(1)}= \left(\Xs \times \ds{\Xs}\right) \times \V{b}^{(1)},
\end{equation}
for $\V{b}^{(1)}(s)$, with $\V{b}^{(2)}=\Xs(s) \times \V{b}^{(1)}$. The material frame and Bishop frame are related by \cite{bergou2008discrete}
\begin{gather}
\label{eq:D1}
\DO = \V{b}^{(1)} \cos{\theta}  + \V{b}^{(2)}\sin{\theta}, \qquad \DT = - \V{b}^{(1)}\sin{\theta}  + \V{b}^{(2)} \cos{\theta}, \qquad \DR = \Xs,
\end{gather}
so that $\theta(s)$ is the angle of rotation of the material frame relative to the Bishop frame. \rev{Solving\ \eqref{eq:BishODE} is also called ``parallel transporting'' the vector $\V{b}^{(1)}$, since its tangential component is kept at zero by performing rotations about the binormal \cite{bergou2008discrete}.} Using\ \eqref{eq:D1}, it can be shown that the twist\ \eqref{eq:sctwist} simplifies to \cite{bergou2008discrete}
\begin{equation}
\label{eq:twist}
\psi=\ds{\theta}.
\end{equation}
If $\theta$ is a constant, the material frame is just a constant rotation of the twist-free Bishop frame, so the twist density $\psi=\ds{\theta}=0$ as well.

\subsubsection{Inextensibility}
The fibers that we consider are inextensible, meaning that the constraint
\begin{equation}
\label{eq:inex1}
\norm{\Xs}^2=\Xs \cdot \Xs=1
\end{equation}
has to preserved by the fiber evolution. This implies that the variable $s$ is arclength. The form of the constraint\ \eqref{eq:inex1} that we will use is the time derivative 
\begin{equation}
\label{eq:inex2}
\dt{\Xs}\cdot \Xs= \ds{\V{U}} \cdot \Xs=0, 
\end{equation}
i.e., that the velocity of the tangent vector is orthogonal to itself; recall that we omit the explicit dependence on $t$ for brevity. This can be enforced using either a penalty \cite{lim2008dynamics, olson2013modeling} or a constrained \cite{ts04, keavRPY, inexRS, maxian2021integral} approach, as we discuss further in what follows.

\subsubsection{Evolution of the material frame and twist angle}
Finally, we consider the time evolution of the material frame $\left(\DO(s), \DT(s), \DR(s)\right)$. Because the frame is orthonormal and stays orthonormal for all time, its velocity is constrained to satisfy
\begin{equation}
\label{eq:Devolve}
\dt{\V{D}^{(i)}} = \V{\Omega}\times \V{D}^{(i)},
\end{equation}
for some $\V{\Omega}(s)$. In the particular case of $i=3$, we have the evolution equation for $\Xs$, 
\begin{equation}
\label{eq:Xsevolve}
\dt{\Xs} = \V{\Omega} \times \Xs,
\end{equation}
which automatically satisfies the inextensibility constraint\ \eqref{eq:inex2}. 

The evolution of the twist \rev{density} can be written in terms of the angular velocity $\V{\Omega}$ as
\begin{align}
\label{eq:twODE}
\dt{\psi}= \ds{\Omega^\parallel}-\left(\V{\Omega} \cdot \ds{\Xs}\right),
\end{align}
where $\Omega^\parallel=\V{\Omega} \cdot \Xs$. Equation\ \eqref{eq:twODE} is derived using variational arguments in \cite{powers2010dynamics}; in Appendix\ \ref{sec:twderiv} we rederive it using vector calculus. In the next two sections, we discuss how to use the Kirchhoff and Euler models to obtain $\V{\Omega}$, which closes the evolution equation for $\psi$. 

\subsection{The Kirchhoff model \label{sec:krod}}
The distinguishing feature of what we refer to as the \emph{Kirchhoff model} is that the rotational velocity of the tangent vector must be exactly equal to the local angular fluid velocity\ \eqref{eq:OmIBdef},
\begin{equation}
\label{eq:fluvelK}
\V{\Omega} = \V{\Psi} \follows
\dt{\Xs} = \V{\Psi} \times \Xs.
\end{equation}
However, because $\dt{\Xs}= \ds{\V{U}}$, the assumption\ \eqref{eq:fluvelK} implies a local constraint that relates the translational and rotational velocity of the rod
\begin{equation}
\label{eq:simplest_cons}
\ds{\V{U}} =\V{\Psi} \times \Xs.
\end{equation}
Thus the fluid velocity is (locally) constrained so that the integrity of the rod cross section is preserved, i.e., that the rod is unshearable \cite{powers2010dynamics}, and inextensible, since\ \eqref{eq:simplest_cons} implies the inextensibility constraint\ \eqref{eq:inex2}. 

In the Kirchhoff model as presented in \cite{lim2008dynamics, olson2013modeling, keavRPY}, the total force density $\V{f}(s)$ and torque density $\V{n}(s)$ exerted by the fiber on the fluid are given by the force and torque balances \cite{keavRPY}
\begin{gather}
\label{eq:fforce}
\V{f} = \ds{\V{F}}, \qquad \V{n}=\ds{\V{N}} + \Xs \times \V{F}.
\end{gather}
Here the torque density $\ds \V{N}(s)$ comes from the internal moments on the rod $\V{N}(s)$, which are given in terms of the bending modulus $\kappa$ and twist modulus $\twmod$, both of which we assume to be constant along the cylindrical rod. The internal moment $\V{N}(s)$ can be written in terms of the material frame as \cite{lim2008dynamics, olson2013modeling, keavRPY}
\begin{align}
\label{eq:KirchN}
\V{N}&= \kappa \left(-\DT \cdot \ds{\Xs} \right) \DO+ \kappa\left( \ds{\Xs} \cdot \DO \right) \DT + \twmod \psi \Xs \\ \nonumber
&= \kappa \left(\left(\DO \times \DT\right) \times \ds{\Xs} \right) + \twmod \psi\Xs \\ \nonumber 
& = \kappa \left(\Xs \times \ds{\Xs} \right) + \twmod \psi \Xs,
\end{align}
subject to the free fiber boundary conditions \cite{powers2010dynamics}
\begin{equation}
\label{eq:XssBC}
\rev{\V{N}(s=0,L)=\V{0} \follows \ds{\Xs}(s=0,L)=\V{0}, }
\end{equation}
and twist boundary conditions\ \eqref{eq:twBC}. The torque density $\ds{\V{N}}(s)$ in\ \eqref{eq:fforce} is
\begin{equation}
\label{eq:Ns}
\ds{\V{N}}=  \kappa \left(\Xs \times \ds^2{\Xs} \right) + \twmod \left(\Xs \ds{\psi}+ \psi \ds{\Xs}\right),
\end{equation}
with the remaining torque density coming from the force $\V{F}(s)$, which is subject to the free fiber boundary conditions
\begin{equation}
\label{eq:FBC}
\V{F}(s=0,L)=\V{0}.
\end{equation}

\rev{As discussed in \cite[p.~5]{keavRPY}, the force $\V{F}$ is a Lagrange multiplier to enforce the position-based constraint $\ds{\V{X}}=\Xs$. It turns out, however, that the final solution for $\V{F}$ depends on the internal torque\ \eqref{eq:Ns}, since the perpendicular torque from $\V{F}$ \emph{must} eliminate the perpendicular torque coming from $\ds{\V{N}}$ for the constraint\ \eqref{eq:simplest_cons} to hold. This is because perpendicular angular velocity $\rr{\V{\Psi}}^\perp \sim \epsRS^{-2} \V{n}^\perp$, while translational velocity $\tr{\V{U}}\sim \V{n}^\perp \log{\epsRS }$. Thus, it is impossible for\ \eqref{eq:simplest_cons} to hold if there is a $\mathcal{O}(1)$ perpendicular torque on the fluid, so $\V{F}$ must cancel the $\mathcal{O}(1)$ perpendicular torque in $\ds{\V{N}}$.}


We now write an explicit formula for $\V{F}$ that yields a torque density $\V{n}$ which is only in the parallel direction, plus a Lagrange multiplier torque in the perpendicular direction. Because this formula eliminates perpendicular torque from the constitutive model, it also eliminates the dependency of torque on the bending modulus $\kappa$, and leaves only parallel torque depending on the twist modulus $\twmod$. The main step is to introduce the forces from bending $\V{F}^{(\kappa)}(s)$ and twisting $\V{F}^{(\twmod)}(s)$, such that when we take their cross product with $\Xs$ and add it to $\ds{\V{N}}$ in\ \eqref{eq:ndecomp}, the perpendicular components of torque will be eliminated
\begin{gather}
\label{eq:Ftwkap}
\V{F}^{(\kappa)} = -\kappa \ds^2{\Xs}, \qquad \V{F}^{(\twmod)}=\gamma \psi \left(\Xs \times \ds{\Xs}\right).
\end{gather}
We then redefine the force $\V{F}$ in the Kirchhoff model as
\begin{gather}
\label{eq:Fdecomp}
\V{F} = \V{F}^{(\kappa)}+\V{F}^{(\twmod)}+\V{\Lambda}, 
\end{gather}
where now $\V{\Lambda}(s)$ is the Lagrange multiplier that enforces \rev{the velocity-based constraint}\ \eqref{eq:simplest_cons}. For free fibers, the force $\V{F}$ must vanish at the endpoints (see\ \eqref{eq:FBC}), which, combined with\ \eqref{eq:XssBC}, leads to the free fiber boundary conditions
\begin{equation}
\label{eq:EulBC}
\ds{\Xs}(s=0,L)=\V{0}, \quad \ds^2{\Xs}(s=0,L)=\V{0}, \quad \V{\Lambda}(s=0,L)=\V{0}, \quad \text{and} \quad \psi(s=0,L)=0.
\end{equation}
Substituting the new form of $\V{F}$ in\ \eqref{eq:Fdecomp} into\ \eqref{eq:fforce}, we have the torque density applied to the fluid as
\begin{align}
\label{eq:ndecomp}
\V{n} &= \ds{\V{N}} + \Xs \times \left(\V{F}^{(\kappa)}+\V{F}^{(\twmod)}+\V{\Lambda}\right)\\ \label{eq:torqK}
& =\left(\twmod \ds{\psi}\right) \Xs + \Xs \times \V{\Lambda}=n^\parallel \Xs + \Xs \times \V{\Lambda}^\perp,
\end{align}
where we defined the parallel torque coming from twist as \cite{powers2010dynamics}
\begin{equation}
\label{eq:nparF}
n^\parallel = \twmod \ds{\psi}.
\end{equation}
Throughout this section, we have assumed that the filament is intrinsically straight and untwisted; expressions for $n^\parallel$ and $\V{F}$ in the case of an intrinsically bent and twisted rod are given in Appendix\ \ref{sec:Intrin}.

Note that only the (two) perpendicular components of $\V{\Lambda}$ enter in the torque\ \eqref{eq:torqK}. The third Lagrange multiplier required to enforce\ \eqref{eq:simplest_cons} is the parallel component of $\V{\Lambda}$, which enters in the force density applied to the fluid as $\V{\lambda}(s)$, 
\begin{gather}
\label{eq:forceK}
\V{f} = \ds{\left(\V{F}^{(\kappa)}+\V{F}^{(\twmod)}+\V{\Lambda}\right)}=\V{f}^{(\kappa)}+\V{f}^{(\twmod)}+\V{\lambda}, \qquad \text{where}\\ \nonumber
\V{f}^{(\kappa)}= -\kappa \ds^3{\Xs}, \qquad
\V{f}^{(\twmod)} = \twmod \ds{\left(\psi \left(\Xs \times \ds{\Xs} \right)\right)}, \quad \text{and} \quad \V{\lambda}=\ds{\V \Lambda}.
\end{gather}
In our reformulation of the Kirchhoff model, the force density on the fluid is given by\ \eqref{eq:forceK} and the torque is given by\ \eqref{eq:torqK}. \rev{This makes it obvious that the only perpendicular torque acting on the fluid comes from the Lagrange multipliers $\V{\Lambda}^\perp$, while the parallel torque comes from twist elasticity.} 

\subsubsection{Linear system for $\V{\Lambda}$}
To derive a linear system for $\V{\Lambda}$, we first write the constraint\ \eqref{eq:simplest_cons} in terms of the total force density $\V{f}(s)$ and torque density $\V{n}(s)$ applied by the rod to the fluid using the mobility operators defined in\ \eqref{eq:UIBdef} and\ \eqref{eq:PsiMobs},
\begin{equation}
\label{eq:kcons}
\ds{\left( \Lop{M}_\text{tt} \V{f} +  \Lop{M}_\text{tr} \V{n}\right)}= \left( \Lop{M}_\text{rt} \V{f} +  \Lop{M}_\text{rr} \V{n}\right)  \times \Xs. 
\end{equation}
We then substitute the representations for $\V{f}$ and $\V{n}$ from\ \eqref{eq:forceK} and\ \eqref{eq:torqK}, and define a linear operator $\Lop{C}$ such that $\left(\Lop{C}\V{x}\right)(s)= \Xs(s) \times \V{x}(s)$. This yields the second-order nonlocal boundary value problem (BVP)
\begin{gather}
\label{eq:kirch}
\left[\partial_s\left(\Lop{M}_\text{tt}\partial_s+\Lop{M}_\text{tr} \Lop{C}\right) + \Lop{C}\left(\Lop{M}_\text{rt}\partial_s+\Lop{M}_\text{rr} \Lop{C}\right)\right] \V{\Lambda} =\\ \nonumber -\ds{\left(\Lop{M}_\text{tt}  \left(\V{f}^{(\kappa)}+\V{f}^{(\twmod)}\right)+\Lop{M}_\text{tr} \left(n^\parallel \Xs\right)\right)} - \Lop{C}\left(\Lop{M}_\text{rt} \left(\V{f}^{(\kappa)}+\V{f}^{(\twmod)}\right)+\Lop{M}_\text{rr}\left(n^\parallel \Xs\right)\right)
\end{gather}
which can be solved for the Lagrange multipliers $\V{\Lambda}(s)$. Using the force $\V{f}$ in\ \eqref{eq:forceK} and torque in\ \eqref{eq:torqK} then gives the resulting velocities in the Kirchhoff model, 
\begin{gather}
\label{eq:UKdef}
\V{U}= \Lop{M}_\text{tt}\left(\V{f}^{(\kappa)}+\V{f}^{(\twmod)}+\V{\lambda}\right)+  \Lop{M}_\text{tr}\left(n^\parallel \Xs+\Xs \times \V{\Lambda}\right), \\ 
\label{eq:OmKdef}
\V{\Omega}=\V{\Psi}= \Lop{M}_\text{rt}\left(\V{f}^{(\kappa)}+\V{f}^{(\twmod)}+\V{\lambda}\right)+  \Lop{M}_\text{rr}\left(n^\parallel \Xs+\Xs \times \V{\Lambda}\right) .
\end{gather}
\rev{In Appendix\ \ref{sec:ftilde}, we show that these dynamics dissipate elastic energy by doing work on the fluid, with the Lagrange multiplier $\V{\Lambda}$ not doing any additional work as required by the principle of virtual work\footnote{\rev{For a
clamped rotating end, it can be shown that the motor rotating the fiber with angular frequency $\omega$ does additional work with power $N_0 \omega$, where $N_0$ is the applied torque at the clamped end, see\ \eqref{eq:torqBC}.}} \cite{varibook}.  }

\subsection{The Euler model \label{sec:euler}}
The main difference between the Kirchhoff model and what we refer to as the \emph{Euler model} is how the rotation rate of the cross section (equivalently, of the tangent vector) is obtained. In the Kirchhoff model, it does not matter whether we use $\ds{\V{U}}$ or $\V{\Psi} \times \Xs$ to update the tangent vector, since by\ \eqref{eq:simplest_cons} they must be the same. In the Euler model, by contrast, we do not set $\V{\Omega}$ to be equal to $\V{\Psi}$ as in\ \eqref{eq:fluvelK}, and consequently do not enforce the constraint\ \eqref{eq:simplest_cons} that follows from\ \eqref{eq:fluvelK}. Instead, we relate $\V{\Omega}^\perp(s)$ to the \emph{translational} velocity via
\begin{equation}
\label{eq:OmEul}
\V{\Omega}^\perp = \Xs \times \ds{\V{U}}, \qquad \text{and} \qquad \dt{\Xs}=\ds{\V{U}} = \V{\Omega}^\perp \times \Xs,
\end{equation}
where $\V{\Omega}_E^\perp(s) \cdot \Xs=0$, and, unlike the Kirchhoff model, $\V{\Omega}^\perp$ is not necessarily equal to $\V{\Psi}^\perp$. The evolution of the fiber can now be obtained by solving for $\V{U}(s)$ and specifying $\Omega^\parallel(s)=\Psi^\parallel(s)$.  

In the Euler model, there is only one constraint, the inextensibility constraint\ \eqref{eq:inex2}. This means that the Lagrange multiplier to enforce this constraint should have only one independent component, as opposed to the three independent components for the Lagrange multiplier in the Kirchhoff model. Specifically, the Lagrange multiplier $\V{\Lambda}(s)$ reduces to a scalar line tension $T(s)$, 
\begin{equation}
\label{eq:LamEul}
\V{\Lambda} = T \Xs, 
\end{equation}
with the free-fiber boundary condition\ \eqref{eq:EulBC} implying that 
\begin{equation}
\label{eq:Tfree}
T(s=0,L)=0.
\end{equation}
Since $\V{\Lambda}$ is tangential, the torque on the fluid in\ \eqref{eq:torqK} reduces to $\V{n}=n^\parallel \Xs$, where $n^\parallel$ is defined in\ \eqref{eq:nparF}. The force is still given by\ \eqref{eq:forceK}, but now with $\V{\Lambda}$ having only one independent component.

In the Euler model, the translational velocity $\V{U}(s)$ is obtained by solving a constrained system which includes only the parallel torque density. This constrained system is derived by substituting\ \eqref{eq:LamEul} into the Kirchhoff velocity\ \eqref{eq:UKdef}, 
\begin{gather}
\label{eq:eulerbeam}
\V{U}= \Lop{M}_\text{tt}\left(\V{f}^{(\kappa)}+\V{f}^{(\twmod)}+\V{\lambda}\right)+  \Lop{M}_\text{tr}\left(n^\parallel \Xs\right),
\end{gather}
where in the Euler model $\V{\lambda}=\ds{\left(T\Xs\right)}$. By substituting\ \eqref{eq:eulerbeam} \rev{into\ \eqref{eq:inex2}}, we can obtain a (second order, nonlocal) BVP for $T(s)$ \cite{ts04} similar to BVP\ \eqref{eq:kirch} for the Kirchhoff model. As discussed in Section\ \ref{sec:inext}, in spectral numerical methods it is better to solve a saddle-point linear system for $\V{U}$ and $\V{\lambda}$ instead of this tension equation for $T(s)$. The saddle point system is obtained from the principle of virtual work \cite[Sec.~3]{maxian2021integral}; see also Appendix\ \ref{sec:ftilde}.

The evolution of the twist angle in\ \eqref{eq:twODE} requires the rotation rate $\Omega^\parallel(s)$, which is computed in a post-processing step as the parallel component of the angular velocity in the Kirchhoff model\ \eqref{eq:OmKdef}, 
\begin{gather}
\label{eq:OmE}
\Omega^\parallel= \left( \Lop{M}_\text{rt} \left(\V{f}^{(\kappa)}+\V{f}^{(\twmod)}+\V{\lambda}\right)+  \Lop{M}_\text{rr}  \left(n^\parallel \Xs \right)\right)\cdot \Xs.
\end{gather}
The Euler model \emph{does not} enforce the constraint\ \eqref{eq:simplest_cons}, since only the translational velocity $\V{U}$ and parallel rotational velocity $\Omega^\parallel$ are tracked. In Section\ \ref{sec:compareMeth}, we study the effect of this on the velocity of the fiber centerline.

We emphasize that our choice to name the model of this section the Euler model and that of the previous section the Kirchhoff model is historical and arbitrary. The reason we choose the naming we do is because the Kirchhoff model as formulated in Section\ \ref{sec:krod} is what has been used before \cite{lim2008dynamics, olson2013modeling, keavRPY}, and the Euler model is the one that is closest to what we (and many others) have used before without twist \cite{ts04, ehssan17, maxian2021integral}. 

\subsection{\rev{Comparing the Kirchhoff and Euler formulations} \label{sec:compareMeth}}
In this section, we analyze the difference between the Kirchhoff model of Section\ \ref{sec:krod} and the Euler model of Section\ \ref{sec:euler}. Because the Kirchhoff model accounts for the angular fluid velocity and has more degrees of freedom, we refer to the translational velocity obtained from\ \eqref{eq:UKdef} as the ``exact'' velocity $\V{U}$, with corresponding tangent vector evolution $\dt{\Xs} = \ds{\V{U}}=\V{\Psi} \times \Xs=\V{\Omega} \times \Xs$. We refer to the translational velocity from the Euler model\ \eqref{eq:eulerbeam} as $\Eul{\V{U}}$ and the parallel rotational velocity from the Euler model\ \eqref{eq:OmE} as $\Eul{\Omega^{\parallel}}$. Since the rotation of the tangent vector in the Euler model is obtained from $\ds{\Eul{\V{U}}}$, the parallel rotational velocity $\Eul{\Omega^{\parallel}}$ represents the only remaining degree of freedom. To examine the difference between the models, we substitute the definitions\ \eqref{eq:UKdef},\ \eqref{eq:eulerbeam},\ \eqref{eq:OmKdef}, and\ \eqref{eq:OmE} to define
\begin{align}
\label{eq:dU}
\dU = \V{U}-\Eul{\V{U}}&= \Mctt \left(\V{\lambda}-\V{\lambda}_E\right)+ \Lop{M}_\text{tr} \left(\Xs \times \dF\right),\\ 
\label{eq:dOm}
\dOmpar = \Omega^\parallel-\Eul{\Omega^\parallel} &= \left(\Lop{M}_\text{rt}\left(\V{\lambda}-\V{\lambda}_E\right)+  \Lop{M}_\text{rr}\left(\Xs \times \dF \right)\right) \cdot \Xs:=\Delta \V \Psi \cdot \Xs
\end{align}
The differences in velocity\ \eqref{eq:dU} and\ \eqref{eq:dOm} depend on the difference in the constraint force
\begin{equation}
\label{eq:dF}
\dF := \V{\Lambda}-\V{\Lambda}_E,
\end{equation}
and its derivative, the force density $\V{\lambda}-\V{\lambda}_E$, where $\V{\Lambda}_E=T\Xs$ and $\V{\lambda}_E=\ds{\left(T\Xs\right)}$ are the constraint force and force density from the Euler model. 

To obtain $\dU$ and $\dOmpar$, we solve the ``mismatch'' problem
\begin{gather}
\label{eq:mprobT}
\partial_s \Delta \V{U} + \Xs \times \Delta \V{\Psi} = -\V{m} \\
\label{eq:mmT}
\V{m}= \ds{\left(\Lop{M}_\text{tt} \left(\ds{\Eul{\V{F}}}\right)+ \Lop{M}_\text{tr} (n^\parallel \Xs ) \right)}+ \Xs \times \left(\Lop{M}_\text{rt} \left(\ds{\Eul{\V{F}}}\right)+ \Lop{M}_\text{rr} \left(n^\parallel \Xs \right)\right), \\
\label{eq:F0}
\qquad \text{where} \qquad \V{F}_E = \V{F}^{(\kappa)}+\V{F}^{(\twmod)} + \V{\Lambda}_E,
\end{gather}
is the total force in the Euler model. This problem, which we obtain by substituting\ \eqref{eq:dF} into the Kirchhoff model\ \eqref{eq:kirch}, gives a second order nonlocal BVP for $\dF$ as a function of the inconsistency in the evolution of the tangent vector in the Euler formulation $\V{m}(s)$. Because the Euler formulation does not enforce\ \eqref{eq:simplest_cons}, the mismatch $\V{m}$ is nonzero in general. We therefore need the Lagrange multiplier $\dF$ to correct for it and make the total Kirchhoff velocities $\V{U}=\Eul{\V{U}}+\dU$ and $\V{\Omega}=\V{\Psi}=\Eul{\V{\Psi}}+\dOm$ consistent with\ \eqref{eq:simplest_cons}. This Lagrange multiplier can be split into two components: a perpendicular part $\dF^\perp$ that contributes to the torque $\Xs \times \dF$, and a parallel part $\dF^\parallel$ which preserves the inextensibility of the velocity $\dU$. Because the constraint force $\V{\Lambda}_E=T\Xs$ in the Euler model is tangential, we have $\dF^\perp=\V{\Lambda}^\perp$ and $\dF^\parallel = \V{\Lambda}^\parallel-T\Xs$, where $T(s)$ is the tension in the Euler model. 

In  Appendix\ \ref{sec:KENumer}, we use the asymptotic slender-body scalings \cite{TwistSBT}
\begin{gather}
\label{eq:LDmats}
\rev{\Lop{M}_\text{tt}\sim \log{\epsc}, \qquad \Lop{M}_\text{tr} \sim \log{\epsc}, \qquad
\Lop{M}_\text{rt}\sim \log{\epsc}, \qquad \Lop{M}_\text{rr} \sim \epsc^{-2},}
\end{gather}
to show that the velocity differences $\dU$ and $\dOmpar$ are small in the fiber interior. Specifically, we show that
\begin{gather}
\label{eq:dUEst}
\norm{\ds{\dU}}/\norm{\ds{\V{U}_E}} \sim  \epsc^2 \log{\epsc} \ll 1,\\ 
\label{eq:dOmEst}
\norm{\Delta \Omega^\parallel} / \norm{\Omega_E^\parallel} \sim \begin{cases} \epsc^4 \left(\log{\epsc}\right)^2 & \twmod > 0 \\ \epsc^2 \log{\epsc} & \twmod=0
\end{cases}
\end{gather}
away from the endpoints, so that the difference between the Kirchhoff and Euler translational and rotational velocities is very small for sufficiently slender fibers. \rev{In Appendix\ \ref{sec:KENumer}, we also verify these asymptotic estimates numerically and discuss how they change at the endpoints.}

\section{Euler model with dynamics\label{sec:EulEqns}}
\rev{In this section, we formulate the full set of dynamical equations for the Euler model of an inextensible filament in Stokes flow. These are the equations of Section\ \ref{sec:euler}, augmented with a parameterization of the set of inextensible fiber velocities, and an adjoint equation for the constraint forces~$\V{\lambda}$.}

\subsection{Continuum formulation for inextensible fibers \label{sec:inext}}
As discussed in Section\ \ref{sec:euler}, most approaches to inextensibility for the Euler model use a scalar Lagrange multiplier for tension $T(s)$ \cite{du2019dynamics}. This is problematic because it creates a complicated auxiliary integro-differential equation for the tension that is difficult to solve with spectral methods. The approach we follow here \rev{\cite{blum1979biophysics}} is detailed in \cite[Sec.~3]{maxian2021integral}. We enforce the velocity constraint\ \eqref{eq:eulerbeam} by restricting the fiber velocity to the space of inextensible motions using\ \eqref{eq:inex1}. Specifically, for free fibers we let\change{
\begin{equation}
\label{eq:Kcdef}
\V{U}(s)=\left(\Lop{K}\V{\alpha}\right)(s) = \V{U}_0+\int_0^{s} \V{\Omega}(s') \times \Xs(s') \, ds',}
\end{equation}
\change{for an arbitrary angular velocity $\V \Omega(s)$, where $\V{\alpha}=\left\{ \V{U}_0, \V{\Omega}\right\}$. The kinematic operator $\Lop{K}\left[\V{X}\right]$ gives a (complete) parameterization of the space of (free fiber) inextensible motions, so that any velocity $\Lop{K}\V{\alpha}$ automatically satisfies the constraint\ \eqref{eq:inex2}. In \cite{maxian2021integral}, we set $\V{\alpha}=\left\{\V{U}_0, \V{\Omega}^\perp\right\}$, but here we allow $\V{\Omega}$ to have a parallel component. This makes no difference in the continuum formula\ \eqref{eq:Kcdef}, but improves numerical robustness, as discussed in Appendix\ \ref{sec:DiscK}.

Instead of setting $\V{\lambda}=\ds{\left(T \Xs \right)}$ and solving for the tension $T(s)$, we solve for a (vector) force density $\V{\lambda}(s)$ that is constrained to satisfy the principle of virtual work (see Appendix\ \ref{sec:ftilde}), 
\begin{align}
\label{eq:kIP1}
\left\langle \Lop{K}\V{\alpha}, \V{\lambda} \right\rangle & = \int_0^L \left( \V{U}_0 + \int_0^s \left(\V{\Omega}(s') \times \Xs(s')\right) \, ds'\right) \cdot \V{\lambda}(s) \, ds \\
& = \V{U}_0 \cdot \int_0^L \V{\lambda}(s') \, ds' + \int_0^L \V{\Omega}(s) \cdot \left(\Xs(s) \times \cdot \int_{s}^L  \V{\lambda}(s') \, ds'\right)  \, ds\\ \nonumber
:&=\left\langle \V{\alpha},\Lop{K}^* \V{\lambda} \right\rangle =0
\end{align}
for all $\V \alpha$. The second equality, which comes from changing the order of integration in the first, immediately gives the constraints
\begin{gather}
\label{eq:KstPt}
\Lop{K}^*\V{\lambda}:=\begin{pmatrix}
\int_0^L \V{\lambda}(s) \, ds \\
\Xs(s) \times \int_s^L \V{\lambda}(s') \, ds'\end{pmatrix} = 
\begin{pmatrix} \V{0} \\ \V{0} \end{pmatrix},
\end{gather}
where the second line holds for all $s$. The second constraint in\ \eqref{eq:KstPt} implies $\V{\lambda}=\ds{\left(T\Xs\right)}$, with the first constraint then giving the appropriate boundary conditions\ \eqref{eq:Tfree} on tension.}

\subsubsection{Modifications for clamped ends}
In this work, we will consider two kinds of filaments: those with two free ends, and those with a clamped end at $s=0$ and a free end at $s=L$. For the latter, the $\Lop{K}$ (and therefore $\Lop{K}^*$) operators have to be modified to account for motion at $s=0$ being prohibited. For clamped fibers, the kinematic operator $\Lop{K} \equiv \Lop{K}_c$ is given by
\begin{equation}
\label{eq:Kclamp}
\change{\Lop{K}_c\V{\alpha}:=\Lop{K}_c \V{\Omega} = \int_0^{s} \V{\Omega}(s') \times \Xs(s') \, ds'.}
\end{equation}
These motions are a subset of those in\ \eqref{eq:Kcdef}, but with $\V{U}_0=\V{0}$. The corresponding virtual work constraint $\Lop{K}_c^*\V{\lambda}=\V 0$ is obtained \change{from\ \eqref{eq:KstPt} by dropping the first row, since the fiber is no longer force-free. While it is also possible to build the clamped condition $\Xs(0) \equiv \Xs_0$ into $\Lop{K}_c$ by restricting $\V{\Omega}^\perp(s=0)= \V 0$, thereby ensuring the clamped boundary condition is numerically satisfied to higher accuracy, doing this causes tension to be ill-defined at $s=0$. In the numerical method in this paper, we do not enforce $\V{\Omega}^\perp(0) = \V 0$ kinematically, but rather enforce it mechanically through the bending force as described in Section\ \ref{sec:XBCF}. We have numerically confirmed that, for sufficiently resolved simulations, the same solution is obtained regardless of whether we enforce $\V{\Omega}^\perp(0)=\V 0$ in $\Lop{K}_c$.}\footnote{\change{The ambiguity in enforcing $\V{\Omega}^\perp(0)=\V 0$ kinematically and/or mechanically has a physical origin in the nonzero length of a clamp. If we consider $s=0$ to be ``inside'' the clamp, then $\V{\Omega}^\perp(0)=\V 0$, but if $s=0$ is ``outside'' the clamp, then $\V{\Omega}^\perp(0)$ is free.}}

\subsection{\rev{Summary of equations}}
To summarize, the saddle-point system that we solve to obtain the motion of the fiber centerline is 
\begin{gather}
\label{eq:euler}
\begin{pmatrix} - \Lop{M}_\text{tt} & \Lop{K} \\ 
\Lop{K}^*& \M{0} \end{pmatrix}
\begin{pmatrix} \V{\lambda} \\ \V{\alpha} \end{pmatrix} 
= \begin{pmatrix}  \Lop{M}_\text{tt}  \V{f}+  \Lop{M}_\text{tr}  \left(n^\parallel \Xs\right)  \\ \V{0} \end{pmatrix}, \qquad \text{where} \\ \nonumber
\V{f}=\V{f}^{(\kappa)} +\V{f}^{(\twmod)}, \qquad
\V{f}^{(\kappa)} = -\kappa \ds^3{\Xs}, \qquad
\V{f}^{(\twmod)}= \twmod \ds{\left(\psi \left(\Xs \times \ds{\Xs} \right)\right)}, \qquad 
n^\parallel =\twmod \ds{\psi},
\end{gather}
subject to the free fiber BCs\ \eqref{eq:EulBC}. Solving the saddle-point system\ \eqref{eq:euler} yields the constraint forces $\V{\lambda}=\ds{\left(T\Xs\right)}$ and translational velocity $\V{U} = \Lop{K}\V{\alpha}$. If desired, the scalar tension $T(s)$ can be extracted from the constraint forces via integration of $\V{\lambda}(s)$. The perpendicular component of the rotation rate $\V{\Omega}(s)$ can be obtained from\ \eqref{eq:OmEul} by
\begin{equation}
\label{eq:OmperpFromK}
\V{\Omega}^\perp = \Xs \times \ds{\V{U}}=\Xs \times \ds{\left(\Lop{K}\V{\alpha}\right)}.
\end{equation}

The parallel component $\V{\Omega}^\parallel(s)$ is not determined by\ \eqref{eq:euler}, and is instead obtained from the post-processing step\ \eqref{eq:OmE}, 
\begin{equation}
\label{eq:OmParAgain}
\Omega^\parallel= \left( \Lop{M}_\text{rt} \V{f}+  \Lop{M}_\text{rr}  \left(n^\parallel \Xs \right)\right) \cdot \Xs,
\end{equation}
The evolution of the twist $\psi=\ds{\theta}$ is given by\ \eqref{eq:twODE}
\begin{align}
\label{eq:PsiEul}
\dt{\psi}&= \ds{\Omega^\parallel}-\left(\V{\Omega} \cdot \ds{\Xs}\right) = \ds{\Omega^\parallel}-\left(\V{\Omega}^\perp \cdot \ds{\Xs}\right)
\end{align}
where in the second equality we have used the formula $\Xs \cdot \ds{\Xs}=0$, which means that the parallel part of $\V{\Omega}$ makes no contribution to the dot product $\V{\Omega}\cdot \ds{\Xs}$.

If needed (e.g., for visualization or intrinsic curvature), the material frame vectors $\DO(s)$ and $\DT(s)$ can be obtained by computing the Bishop frame from\ \eqref{eq:BishODE}, then twisting it by the angle $\theta$. The Bishop frame has to be chosen at one point on the fiber, and the angle $\theta$ is only determined up to a constant. Our choice is to keep track of the material vector $\DO$ at the fiber midpoint\footnote{We use the midpoint because it is best conditioned in a spectral method; ill-conditioning of the spectral derivatives typically happens at the endpoints of the Chebyshev grid first.} by explicitly evolving it via rotation by $\V{\Omega}(L/2)$, 
\begin{equation}
\dt{\DO}(L/2) = \left(\V{\Omega}^\perp(L/2)+\Omega^\parallel(L/2)\Xs(L/2)\right) \times \DO(L/2). 
\end{equation}
We then assign $\theta(L/2)=0$, i.e., require that the Bishop frame be the same as the material frame at $L/2$
\begin{equation}
\label{eq:BishBC}
\V{b}^{(1)}(L/2) = \DO(L/2),
\end{equation}
which provides a boundary condition for the Bishop ODE\ \eqref{eq:BishODE}. See Appendix\ \ref{sec:BishNumer} for how we solve the Bishop ODE to spectral accuracy. 

\section{Spectral spatial discretization \label{sec:specEul}}
\delete{Broadly speaking, our goals for the discretization are to improve on existing numerical methods by reducing their cost. This was our guiding principle in \cite{maxian2021integral}, . As we review in Section\ \ref{sec:inext}, the main contributions of that work were to parameterize exactly the space of inextensible motions, thereby eliminating penalty parameters for inextensibility \cite{ts04},   The latter allows for spectral methods, and represents a major improvement over low-order regularized singularity discretizations in which the number of discrete points on the fiber is approximately $1/\epsRS$. }

In this section, we formulate a spectral spatial discretization of the Euler model. This is an extension of our work \cite{maxian2021integral}, where we designed a numerical method for the case of fibers \emph{without} twist elasticity. In that work, we \rev{followed the idea of Shelley et al.\ \cite{shelley1996nonlocal,shelley2000stokesian,ts04} in using slender body theory for the mobility, which decouples the number of degrees of freedom from the slenderness $\epsc$.} However, we left open a few issues which have plagued numerical methods for SBT for quite some time \cite{ts04, ehssan17}, including the ill-posedness of the continuum equations on lengthscales less than $\epsc$ \cite{gotz2001interactions}, and singularities at the endpoints for non-cylindrical fibers \cite[Sec.~2.1]{maxian2021integral}. As has been mentioned in Section\ \ref{sec:gmob}, our solution to both of these problems is to take the action of the mobility operator to be an \emph{integral} of the RPY kernel along the fiber centerline. \delete{Since this approach automatically gives a mobility operator with positive eigenvalues whose action is well-defined up to and including the endpoints, our task is to discretize it to accuracy sufficient to obtain a mobility \emph{matrix} $\Mtt$ that represents the trans-trans mobility \emph{operator} $\Lop{M}_\text{tt}$.}

In Section\ \ref{sec:exactmob} and Appendix\ \ref{sec:quads} we develop efficient quadrature schemes to evaluate the trans-trans, rot-trans, and trans-rot RPY integrals with on the order of ten points along the fiber centerline. For the short-ranged rot-rot mobility, we use the local drag formulas of Section\ \ref{sec:RPYasymp}, which are sufficiently accurate for our purposes. \delete{It is possible to use the asymptotic formulas of Appendix\ \ref{sec:EPs} for the other mobilities as well;  Section\ \ref{sec:asympmob} presents the minimal modifications required in this case. No matter the choice of mobility,}Since we expect the fiber shapes (and velocity on the fiber centerline) to be smooth, we represent the fiber centerline as a Chebyshev interpolant, and evaluate the mobility and forces at Chebyshev collocation points. Forces and torques at these points must be computed with proper treatment of the boundary conditions. Our choice is to use (a modified form of) rectangular spectral collocation \cite{dhale15, tref17}, as we discuss in Section\ \ref{sec:BCs}.

In Section\ \ref{sec:staticconv}, we show that the solution of the Euler problem \emph{with twist} can be well-approximated with a spectral method, which drastically lowers the number of degrees of freedom relative to previous methods \cite{keavRPY, lim2008dynamics}. We demonstrate this in detail in a sequence of steps, first computing the eigenvalues of the discrete mobility matrix, then the accuracy of the forward quadratures, and finally the convergence of the solution to the static Euler problem. 

\subsection{Spectral discretization \label{sec:specDisc}}
Our spectral spatial discretization of the Euler equation\ \eqref{eq:euler} is based on our previous work \cite{maxian2021integral}. To discretize the fiber centerline, we introduce an $N$-point type 1 (not including the endpoints) Chebyshev collocation grid $\{s_p\}$ for $p=1, \dots N$. In the discrete setting, we use $\V{X}(s)$ for the Chebyshev interpolant approximating the fiber centerline.

\change{In this paper, we modify the discretization of $\Lop{K}$ and $\Lop{K}^*$ as matrices $\M{K}$ and $\M{K}^*$ from \cite{maxian2021integral} to robustly simulate more curved fibers (see example in Section\ \ref{sec:whirl}). In \cite{maxian2021integral}, we introduced an orthonormal frame $\left(\Xs(s),\V{n}_1(s),\V{n}_2(s)\right)$ at each collocation point, then parameterized $\V{\Omega}^\perp$ in terms of $\V{n}_1$ and $\V{n}_2$. This removes the null space from $\M{K}$, but for curved fibers it leads to aliasing problems because products of Chebyshev polynomials of degree $N$ are aliased on an $N$-point grid. In Appendix\ \ref{sec:DiscK}, we describe a discretization of $\Lop{K}$ and $\Lop{K}^*$ which is based on removing the null space of $\M{K}$ \emph{numerically} by doing all computations on a grid of size $2N$, then downsampling. That is the discretization we use in this paper.}



\subsection{Slender-body quadrature for RPY mobilities \label{sec:exactmob}}
We will impose\ \eqref{eq:euler} in the strong sense on the Chebyshev collocation grid. Since $\M{K}$ gives the inextensible velocity at each point on the grid, we need to evaluate the action of the mobility operators $\Lop{M}_\text{tt},\Lop{M}_\text{rt},\Lop{M}_\text{tr}$, and $\Lop{M}_\text{rr}$ at each point on the Chebyshev grid, giving discrete mobility matrices $\Mtt, \Mrt, \Mtr$, and $\Mrr$. We discuss how to do this in this section.

Our approach is to use the translational mobility operator $\Lop{M}_\text{tt}$ given by the RPY integral\ \eqref{eq:UIBdef}. If the integral is computed to sufficient accuracy, then the eigenvalues of the discrete mobility matrix $\Mtt$ are guaranteed to be positive since the RPY kernel is SPD \cite{wajnryb2013generalization}. \change{Note, however, that while the RPY kernel acting on \emph{forces} is symmetric, the matrix $\Mtt$ that acts on force densities is not symmetric, since quadrature is applied on one side but not the other.} While\ \eqref{eq:UIBdef} is still a ``first kind'' integral operator, and therefore has eigenvalues that cluster around zero, unlike the ``second kind'' integral operator proposed in \cite{andersson2021integral} to regularize SBT, we will demonstrate that the ill-conditioning of $\Mtt$ is \emph{not} a problem in practice. Furthermore, unlike \cite{andersson2021integral}, we develop a nearly-singular quadrature scheme so that the RPY integrals for smooth fibers can be evaluated to 3 digits of accuracy  with $10-40$ collocation points, regardless of the slenderness $\epsRS$. 

\subsubsection{Translational mobility \label{sec:singsub}}
We define the translational velocity of the fiber centerline due to a force density $\V{f}(s)$ using\ \eqref{eq:UIBdef}, 
\begin{align}
 \label{eq:transmob}
\tt{\V{U}}(s) &= \left(\Lop{M}_\text{tt}\V{f}\right)(s)=\int_0^L \Mbtt\left(\V{X}(s),\V{X}(s^\prime)\right)\V{f}(s^\prime) \, ds'\\ \nonumber
& = \int_{D(s)} \left(\Slet{\V{X}(s),\V{X}(s')} +\frac{2\eps^2}{3}\Dlet{\V{X}(s),\V{X}(s')} \right)\V{f}\left(s'\right) \, ds' \\[2 pt] \nonumber &+ \EPMI \int_{D^c(s)} \left(\left(\dfrac{4}{3\eps}-\dfrac{3R\left(\V{X}(s),\V{X}(s')\right) }{8\eps^2}\right)\M{I}+\dfrac{\V{R} \Rhat\left(\V{X}(s),\V{X}(s')\right) }{8\eps^2} \right)\V{f}(s') \, ds'.
\end{align}
Here we have used the definition\ \eqref{eq:MbttRPY} to split the integral into a region $D(s)$ for $R\left(\V{X}(s),\V{X}(s')\right) > 2\eps$ and $D^c(s)$, which uses the RPY kernel for $R \leq 2\eps$. We use the fiber inextensibility to make the approximation $R \approx |s'-s|$ when $R \lesssim 2\eps$, so that 
\begin{equation}
\label{eq:Ddom}
D(s) =\begin{cases} \left(0,s-2\eps\right) \cup \left(s+2\eps,L\right) & 2\eps \leq s \leq L-2\eps \\
\left(s+2\eps,L\right) & s < 2\eps \\
 \left(0,s-2\eps\right) & s > L-2\eps
\end{cases},
\end{equation}
with the complement $D^c(s) = [0,L] \setminus D(s)$. \rev{In making the approximation $R \approx |s'-s|$, we also assume that the fiber never re-encroaches itself.}

For the integral of the Stokeslet in\ \eqref{eq:transmob}, we use a singularity subtraction technique which is closely tied with the asymptotics of the Stokeslet. In particular, we subtract from the integrand the leading order singular behavior and perform that integral separately, which gives
\begin{align}
\label{eq:Stsubtr}
\tt{\V{U}}^{(S)}=& \int_{D(s)} \Slet{\V{X}(s),\V{X}(s')} \V{f}\left(s'\right)\, ds'\\[2 pt] \nonumber
= \EPMI &\int_{D(s)}\left(\frac{\M{I}+\Xs(s)\Xs(s)}{|s-s'|}\right)\V{f}(s) \, ds'\\ \nonumber+ & \int_{D(s)} \left(\Slet{\V{X}(s),\V{X}(s')} \V{f}\left(s'\right) -  \EPMI \left(\frac{\M{I}+\Xs(s)\Xs(s)}{|s-s'|}\right)\V{f}(s)\right) \, ds' \\[2 pt] \nonumber
:=&\tt{\V{U}}^{(\text{inner, S})}(s)+\tt{\V{U}}^\text{(int, S)}(s)\\ 
\label{eq:UinnerSt}
\text{where} \quad 
8 \pi \mu \tt{\V{U}}^{(\text{inner,S})}(s) =&
\left(\M{I}+\Xs(s)\Xs(s)\right)\V{f}(s)  
\begin{cases}
\log{\left(\dfrac{(L-s)s}{4\eps^2}\right)} & 2\eps < s < L-2\eps \\ 
\log{\left(\dfrac{(L-s)}{2\eps}\right)} & s \leq 2\eps \\ 
\log{\left(\dfrac{s}{2\eps}\right)} & s \geq L-2\eps
\end{cases}
\end{align}
 \delete{The term $\tt{\V{U}}^{(\text{inner, S})}$ is then given for all $a_L(s)\left(\M{I}+\Xs\Xs\right)$, where $a_L~\ln{\epsRS}$ is defined in\ \cite[Eq.~(C.2)]{TwistSBT}. Thus the leading order \emph{asymptotic} behavior is used as a singularity in the integral.Taking the domain $D(s)$ when evaluating $\tt{\V{U}}^\text{(int, S)}$ to be $[0,L]$ gives the finite part integral from SBT. Combining this with (analytical) evaluation of $\tt{\V{U}}^{(\text{inner, S})}$ and asymptotic evaluation of the doublet and $D^c(s)$ integrals in\ \eqref{eq:transmob} yields exactly the classical SBT of Keller and Rubinow \cite{krub}, assuming the choice of $\eps$ given in\ \eqref{eq:ahat}. In this work, we will take a ``slender body \emph{quadrature}'' approach of evaluating\ \eqref{eq:transmob} \emph{numerically, without any asymptotic approximations}, as in slender body \emph{theory}. For the Stokeslet integral, we perform the singularity subtraction technique\ \eqref{eq:Stsubtr} to give the integral $\tt{\V{U}}^\text{(int, S)}(s)$.}A set of corresponding singularity subtraction steps is given for the doublet integral in\ \eqref{eq:DbSS}. In both cases, we obtain nearly singular integrals which are by design smoother than the integral of the Stokeslet/doublet by itself. We then apply specialized, ``slender-body'' quadrature schemes to these integrals, as discussed in Appendix\ \ref{sec:StNS}. The main framework of these schemes was first proposed by Tornberg, Barnett, and af Klinteberg \cite{tornquad,barLud} and used by us previously for the finite part integral in SBT \cite[Sec.~4.2.1]{maxian2021integral}. The idea is to write the integral as the product of a smooth function times a singular function, expand the smooth function in a set of basis functions, and precompute the integrals of the basis functions times the singularity analytically. Using an adjoint method, the calculation of the Stokeslet and doublet integrals can then be reduced to $N$ inner products of two $N$ dimensional vectors, one of which is precomputed. The total cost of these integrals at each time step is therefore $\mathcal{O}(N^2)$. Since all of these calculations are similar to previously-published methods, we discuss them in Appendix\ \ref{sec:quads}. 

After discretizing the Stokeslet and doublet integrals, we are left with the integral over $D^c(s)$ in the third line of\ \eqref{eq:transmob}. The integrand is nonsingular, but behaves like $|s-s'|$, and so for each $s$ we split the domain $D^c$ into $(s-2\eps, s)$ and $(s,s+2\eps)$ (with appropriate modifications at the endpoints), and use $N_2/2$ Gauss-Legendre quadrature points to sample the fiber and force density (i.e., sample the Chebyshev interpolant of each) and evaluate the integral on each of the two subdomains. We use $N_2/2$ points for each of these integrals so that there are a total of $N_2$ additional (local) quadrature nodes per collocation point.

\subsubsection{Rot-trans mobilities}
\delete{Since we use the RPY integral\ \eqref{eq:transmob} to evaluate the translational velocity from the force density on the fiber centerline, a consistent formulation requires that we also treat the \emph{rotational} velocity from force using the RPY integral, instead of asymptotically. Specifically, if we use\ \eqref{eq:transmob}, we must also set.}
For the rotational velocity from force, we have the RPY integrals
\begin{align}
 \label{eq:rottransmob}
\rt{\Omega}^\parallel(s) :&=\left(\Lop{M}_\text{rt}\V{f}\right)(s)= \Xs(s) \cdot \int_0^L \Mbrt\left(\V{X}(s),\V{X}(s^\prime)\right)\V{f}(s^\prime) \, ds' \\
& = \EPMI \int_{D(s)} \dfrac{\left(\V{f}(s') \times \V{R}\left(\V{X}(s'),\V{X}(s)\right)\right) \cdot \Xs(s)}{R\left(\V{X}(s'),\V{X}(s)\right)^3} \, ds' \\ \nonumber 
& + \EPMI \int_{D^c(s)}
\dfrac{1}{2\eps^2}\left(\dfrac{1}{\eps}-\dfrac{3R\left(\V{X}(s'),\V{X}(s)\right)}{8\eps^2}\right)\left(\V{f}(s') \times \V{R}\left(\V{X}(s'),\V{X}(s)\right) \cdot \Xs(s)\right) \, ds',
\end{align}
which we evaluate to spectral accuracy using singularity subtraction and slender-body quadrature for the rotlet integral over $D(s)$ (see Appendix\ \ref{sec:RTNS}), and direct Gauss-Legendre quadrature for the integral over $D^c(s)$ (split into two pieces). 

The trans-rot mobility is the adjoint of the rot-trans mobility,
\begin{align}
 \label{eq:transrotmob}
\tr{\V{U}}(s):&=\left(\Lop{M}_\text{tr}n^\parallel\right)(s) = \int_0^L \Mbrt\left(\V{X}(s),\V{X}(s^\prime)\right)\Xs(s^\prime) n^\parallel(s')\, ds' \\
& = \EPMI \int_{D(s)} \dfrac{\left(\Xs(s') \times \V{R}\left(\V{X}(s'),\V{X}(s)\right)\right) n^\parallel(s')}{R\left(\V{X}(s'),\V{X}(s)\right)^3} \, ds' \\ \nonumber 
& + \EPMI \int_{D^c(s)}
\dfrac{1}{2\eps^2}\left(\dfrac{1}{\eps}-\dfrac{3R\left(\V{X}(s'),\V{X}(s)\right)}{8\eps^2}\right)\left(\Xs(s') \times \V{R}\left(\V{X}(s'),\V{X}(s)\right) n^\parallel(s')\right) \, ds' .
\end{align}
As in the rot-trans mobility, we compute the first (rotlet) integral using slender-body quadrature (see Appendix\ \ref{sec:TRNS}) and the second using direct Gauss-Legendre quadrature.

\subsubsection{Rot-rot mobility \label{sec:rotrotmob}}
For the rot-rot mobility, we will use always use the asymptotic result \cite[Sec.~3.4]{TwistSBT}. In the fiber interior, this reduces to\ \eqref{eq:UinRRpar},
\begin{equation}
\label{eq:rotrotmob}
\rr{\Omega}^\parallel(s) = \left(\Lop{M}_\text{rr}n^\parallel\right)(s)=\EPMI \frac{9n^\parallel(s)}{4 \eps^2}.
\end{equation}
We use the asymptotic result here because the $1/\eps^2$ term is so dominant as to render calculation of the full integral\ \eqref{eq:PsiMobs}, which involves the rapidly-decaying doublet kernel, unnecessary. More importantly, the local operator $\Lop{M}_\text{rr}$ (and the diagonal discrete matrix $\Mrr$) are positive definite and well-behaved up to the end points.

\subsection{Boundary conditions \label{sec:BCs}}
In this paper, we will consider two kinds of boundary conditions (BCs): free ends and clamped ends. In our recent work that did \emph{not} account for twist elasticity \cite{maxian2021integral}, we used rectangular spectral collocation \cite{dhale15, tref17} to impose the boundary conditions. We review that formulation here, and then extend it to the PDE\ \eqref{eq:twODE} for the twist density $\psi=\ds{\theta}$. In order to demonstrate how we treat both types of boundary conditions, we will consider a fiber with the $s=0$ end clamped, and the $s=L$ end free. We assume that the $s=0$ end is spinning with rate $\omega$ (see Section\ \ref{sec:whirl}), so that the boundary conditions for $\V{X}$ and $\psi$ are 
\begin{gather}
\V{X}(0)=\V{X}_0, \quad \Xs(0)=\Xs_0, \quad \ds^2 \V{X}(L)=\V{0}, \quad \ds^3 \V{X} (L)=\V{0}, \\
\label{eq:thetaBCs}
\Omega^\parallel(0)=\rr{\Omega}^\parallel(0)+\rt{\Omega}^\parallel(0)=\omega \rightarrow \left(\Mrr \gamma \ds{\psi}\right)(0) = \omega-\rt{\Omega}^\parallel(0), \\ \nonumber
\psi(L)=0,
\end{gather}
where we recall from Section\ \ref{sec:gmob} that $\rr{\Omega}^\parallel$ and $\rt{\Omega}^\parallel$ represent the angular velocity due to torque and force, respectively.

Another possibility is to prescribe the parallel torque at the clamped end,
\begin{gather}
\label{eq:torqBC}
\V{N}(s=0) \cdot \Xs(0) = \twmod \psi(0)=N_0
\end{gather}
which gives a boundary condition $\psi(0)=N_0/\twmod$ that is easy to handle; we therefore focus on the more challenging\ \eqref{eq:thetaBCs}. See \cite[Sec~II(A)]{nguyen2018impacts} for a discussion of the difference between the constant torque vs.\ constant angular velocity boundary conditions.

\subsubsection{Boundary conditions on $\V{X}$ \label{sec:XBCF}}
The idea of rectangular spectral collocation as proposed in \cite{dhale15, tref17}, and slightly adapted by us in \cite[Sec.~4.1.3]{maxian2021integral}, is to compute an upsampled representation $\widetilde{\V{X}}$ of $\V{X}$ on a type 2 (including the endpoints) Chebyshev grid that satisfies the boundary conditions exactly, and gives $\V{X}$ when downsampled to the original (type-1) Chebyshev grid. Since there are four boundary conditions, the upsampled representation is on a $N+4$ point type 2 Chebyshev grid that includes the endpoints. The unique configuration $\widetilde{\V{X}}$ can be obtained by solving
\begin{equation}
\label{eq:deftilde}
\begin{pmatrix} \M{R}_{\V{X}} \\[2 pt] \M{B} \end{pmatrix} \widetilde{\V{X}} =:\begin{pmatrix} \M{R}_{\V{X}} \\[2 pt] \M{I}(0) \\\M{D}_{N+4}(0)\\\M{D}_{N+4}^2(L)\\\M{D}_{N+4}^3(L) \end{pmatrix} \widetilde{\V{X}} = \begin{pmatrix} \V{X}\\[2 pt] \V{X}_0 \\ \Xs_0 \\ \V 0 \\\V 0 \end{pmatrix}:=\begin{pmatrix} \V{X}\\[2 pt] \V{b}_{\V{X}} \end{pmatrix},
\end{equation}
where $\M{R}_{\V{X}} $ is the downsampling matrix that evaluates the interpolant $\widetilde{\V X}(s)$ on the original $N$ point grid, $\M{I}(0)$ encodes evaluation at $s=0$, and $\M{D}_M^n(s)$ represents evaluation of the $n$th derivative at the point $s$ on a grid of size $M$. In\ \eqref{eq:deftilde}, $\M{B}$ is the $4 \times \left(N+4\right)$ matrix that encodes the boundary conditions and $\V{b}_{\V{X}}$ is the values of the BCs. We now write the affine transformation from $\V{X}$ to $\widetilde{\V{X}}$ as
\begin{equation}
\label{eq:getX}
\widetilde{\V{X}} = \begin{pmatrix} \M{R}_{\V{X}}  \\[2 pt] \M{B} \end{pmatrix}^{-1} \left(\begin{pmatrix} \V{X} \\[2 pt] \M{0} \end{pmatrix} +\begin{pmatrix} \V{0} \\[2 pt] \V{b} \end{pmatrix}\right):= \M{E}_{\V{X}}\V{X}+\V{\beta}_{\V{X}}.
\end{equation}
Note that $\V{\beta}_{\V{X}}$ is independent of time in this case. 

We then use the configuration $\widetilde{\V{X}}$ to compute any quantities with more than one $s$ derivative in the Euler equations. Importantly, the elastic force density is computed on the upsampled grid as $\widetilde{\V{f}}^\kappa=-\kappa \M{D}_{N+4}^4 \widetilde{\V{X}}$ and is then downsampled to the original $N$ point type 1 grid to give the bending elasticity discretization
\begin{equation}
\label{eq:fkapd}
\V{f}^\kappa\left(\V{X}\right)=\M{R}_{\V{X}} \widetilde{\V{f}}^\kappa = -\kappa \M{R}_{\V{X}} \M{D}_{N+4}^4 \left(\M{E}_{\V{X}}\V{X}+\V{\beta}_{\V{X}}\right):=\M{F}_\kappa \V{X}+\V{\beta}_\kappa
\end{equation}

In all cases in which derivatives of $\Xs$ enter in the equations, we compute the derivatives on the $N+4$ point grid from the configuration $\widetilde{\V{X}}$, then downsample them to evaluate the derivatives on the original $N$ point grid, e.g., we use $\M{R}_{\V{X}}\M{D}_{N+4}^2 \widetilde{\V{X}}$ as the discretization of $\ds{\Xs}$. This helps by smoothing the higher order derivatives due to the imposition of smoothing BCs at the endpoints. Note however that $\widetilde{\Xs} = \ds \widetilde{\V{X}}$ is not of unit length at the collocation points, so we use $\Xs$ on the original $N$ point grid instead when a vector of unit length is required (e.g., in the mobility and kinematic matrices).

\subsubsection{Twist angle}
Let us now consider using rectangular spectral collocation to solve for the twist density $\psi$. We define an affine operator which gives an upsampled representation $\tilde{\psi}$ on an $N+2$ point type 2 Chebyshev grid that satisfies the boundary conditions\ \eqref{eq:thetaBCs},
\begin{gather}
\label{eq:Rpsi}
\begin{pmatrix} 
\M{R}_{\psi} \\ \rr{\widetilde{\Omega}^\parallel}(0) \\ \M{I}(L) 
\end{pmatrix}
\tilde{\psi}= 
\begin{pmatrix}
\psi \\ \omega-\rt{\Omega}^\parallel(0) \\ 0 
\end{pmatrix} = \begin{pmatrix} \psi \\ \V{b}_\psi \end{pmatrix}, \\
\label{eq:OmparNp2}
\change{\textrm{where} \qquad \rr{\widetilde{\Omega}^\parallel}\tilde{\psi} = \Mrr \twmod \M{D}_{N+2}\tilde{\psi}},
\end{gather}
\change{is the angular velocity due to torque on the $N+2$ point grid (throughout this paper, the matrix $\Mrr$ is a diagonal matrix on the $N+2$ point grid).} The analogue of\ \eqref{eq:getX} in this case is the affine operator that gives $\widetilde{\psi}$ from $\psi$ on the $N$ point grid, which we define as $\tilde{\psi}= \M{E}_\psi \psi+\V{\beta}_\psi$. 

To calculate the torque on an $N$ point grid, we take derivatives of $\tilde{\psi}$ on the $N+2$ point grid, and evaluate the result on the $N$ point grid using $\M{R}_\psi$, 
\begin{gather}
\label{eq:nparBC}
n^\parallel(\psi) = \gamma \M{R}_{\psi} \M{D}_{N+2} \left( \M{E}_\psi \psi+\V{\beta}_\psi\right):=\M{R}_\psi\left(\widetilde{\V{N}}_\psi \psi +\widetilde{\V{\beta}}_n\right),\end{gather}
where $\widetilde{\V{\beta}}_n$ is time-dependent in general. The twist force $\V{f}^{(\twmod)}$ as defined in\ \eqref{eq:euler} is slightly more complex, given that it is a nonlinear function of $\V{X}$ and involves interaction between $\V{X}$ and $\psi=\ds{\theta}$. To reduce aliasing on the spectral grid, we use the product rule to separate $\V{f}^{(\twmod)}$ into two terms, then use a common $2N$ point grid to compute each term with proper boundary conditions, 
\begin{equation}
\label{eq:ftwmodC}
\V{f}^{(\twmod)}\left(\V{X},\psi\right)= \left(\twmod \left(\ds{\tilde{\psi}} \left(\ds{\widetilde{\V{X}}} \times \ds^2{\widetilde{\V{X}}}\right)+\tilde{\psi} \left(\ds{\widetilde{\V{X}}} \times \ds^3{\widetilde{\V{X}}}\right)\right)\right)
\end{equation}
We perform the multiplications of $\psi$ and derivatives of $\V{X}$ on a $2N$ point grid, then downsample the result to the $N$ point grid to obtain the final force density $\V{f}^{(\twmod)}$ at the collocation points. 

\subsection{\rev{Well-posedness and convergence} of the static problem \label{sec:staticconv}}
In this section, we study the convergence of our numerical method for the case of a static filament (solving for velocity at $t=0$) in three steps. First, we show that our discrete mobility matrix $\Mtt$ in Section\ \ref{sec:exactmob} is indeed positive definite for larger $\epsRS$, in contrast to $\Mtt$ in SBT, which has negative eigenvalues for as few as 15 collocations points on the fiber centerline when $\epsRS=10^{-2}$. In Appendix\ \ref{sec:convquad}, we show that our quadratures for the forward evaluations of $\Mtt \V{f}$, $\Mrt \V{f}$, $\Mtr n^\parallel$, and $\Mrr n^\parallel$ give spectral accuracy in the velocity due to the bending and twisting elastic force/torque, even for a fiber with relatively high curvature. In particular, we find that we can obtain 3 digits of accuracy in the velocities with about 32 collocation points. Finally, we solve for the velocity $\V{U}$ and constraint force $\V{\lambda}$ in the Euler method and examine their convergence. This is where our method hits a snag: if we maintain a cylindrical fiber shape, there is a physical singularity that develops at the fiber endpoints, where the force required to produce a uniform motion becomes very oscillatory \cite[Sec.~4.3.2]{walker2020regularised}. This means that $\V{\lambda}$ fails to converge pointwise, which calls into question the convergence of $\V{U}$ as well. Nevertheless, we show empirically that the spectral method \emph{does} give a convergent velocity and a weakly-convergent constraint force $\V{\lambda}$ (i.e., the moments of $\V{\lambda}$ converge). Furthermore, the accuracy we obtain in the spectral method with just $16-30$ points on the fiber centerline is better than with $1/\epsRS$ points in a second-order method.

\begin{figure}
\centering
\includegraphics[width=0.3\textwidth]{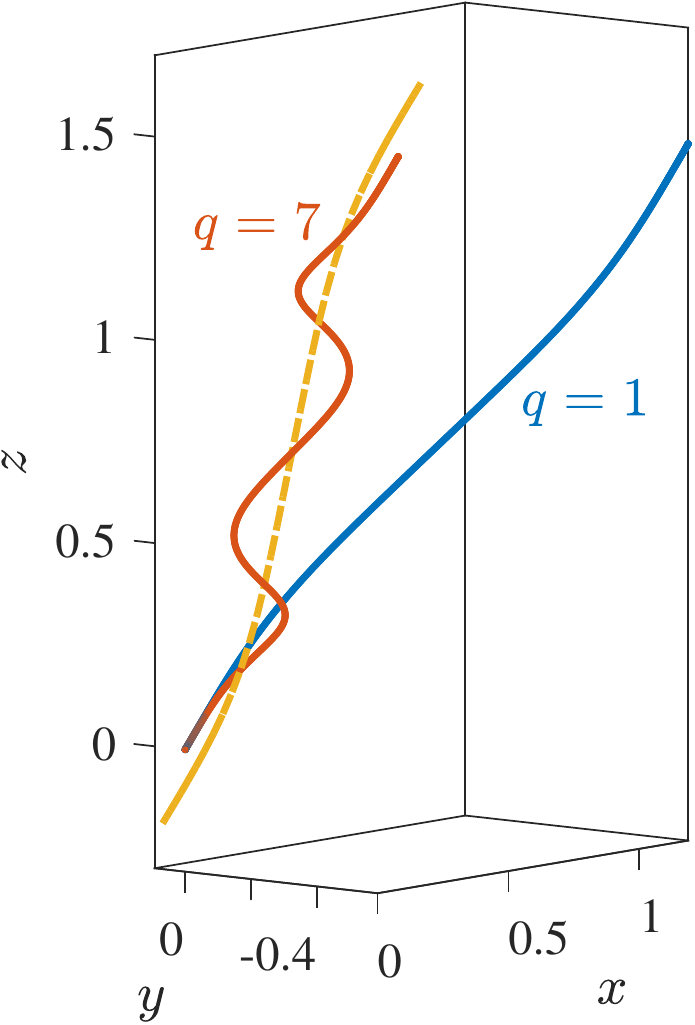}
\includegraphics[width=0.6\textwidth]{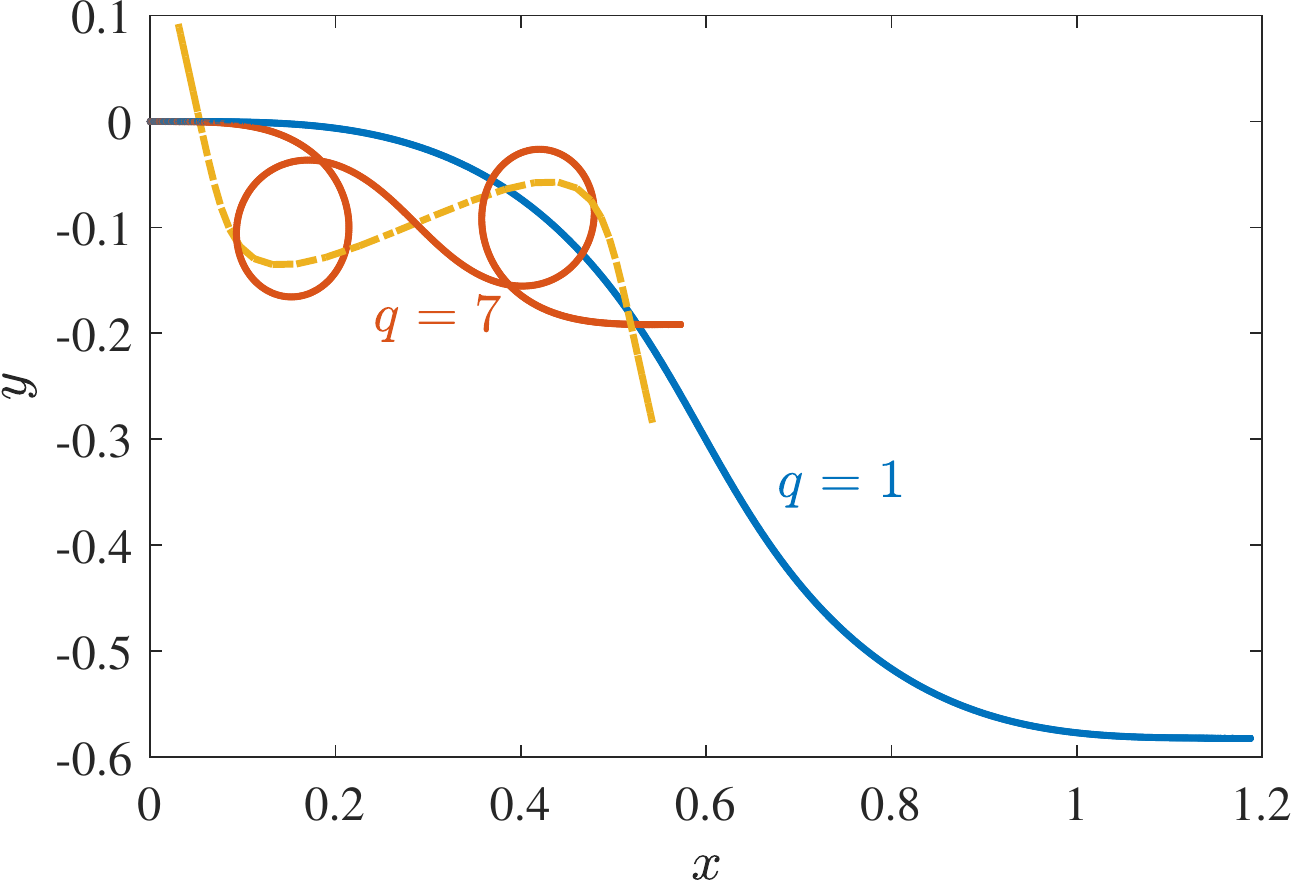}
\caption{\label{fig:FibShapes} Fibers we use for the static convergence study, shown from two different views. Here $L=2$ and the fiber tangent vector is given by\ \eqref{eq:Xsfibq}, with $q=1$ for the blue curve and $q=7$ for the red curve. The dotted red curve shows the $q=7$ fiber after $t=0.01$ seconds of relaxation (when $\epsRS=10^{-2}$, $\mu=1$, $\kappa=1$ and $\twmod=0$; see Section\ \ref{sec:fibrelax}).}
\end{figure}

Throughout this section, we will consider free fibers with tangent vectors of the form 
\begin{equation}
\label{eq:Xsfibq}
\Xs(s) = \frac{1}{\sqrt{2}}\left(\cos{\left(qs^3(s-L)^3\right)}, \sin{\left(qs^3 (s-L)^3\right)}, 1\right),
\end{equation}
where we set $L=2$, and $q$ is a parameter that determines the number of helical turns and fiber curvature and smoothness. Fibers with $q=1$ and $q=7$ are shown in blue and red, respectively, in Fig.\ \ref{fig:FibShapes}. \delete{For our study of the convergence of the forward mobility action, we concern ourselves with robustness, so we use $q=7$, but for the study of the saddle point system convergence, we choose a smooth problem with $q=1$ to reduce the number of collocation points required to achieve high accuracy.}

\subsubsection{Eigenvalues of $\Mtt$}
To illustrate that the RPY mobility\ \eqref{eq:transmob} alleviates the problem of negative eigenvalues that plagues SBT, we compute the eigenvalues of the translation-translation matrix $\Mtt$ numerically. We use $q=7$ in the free fiber\ \eqref{eq:Xsfibq}, although the mobility does not change substantially with the fiber shape (specifically, a straight fiber and four-turn helix give similar results). 

\begin{figure}
\centering
\includegraphics[width=0.49\textwidth]{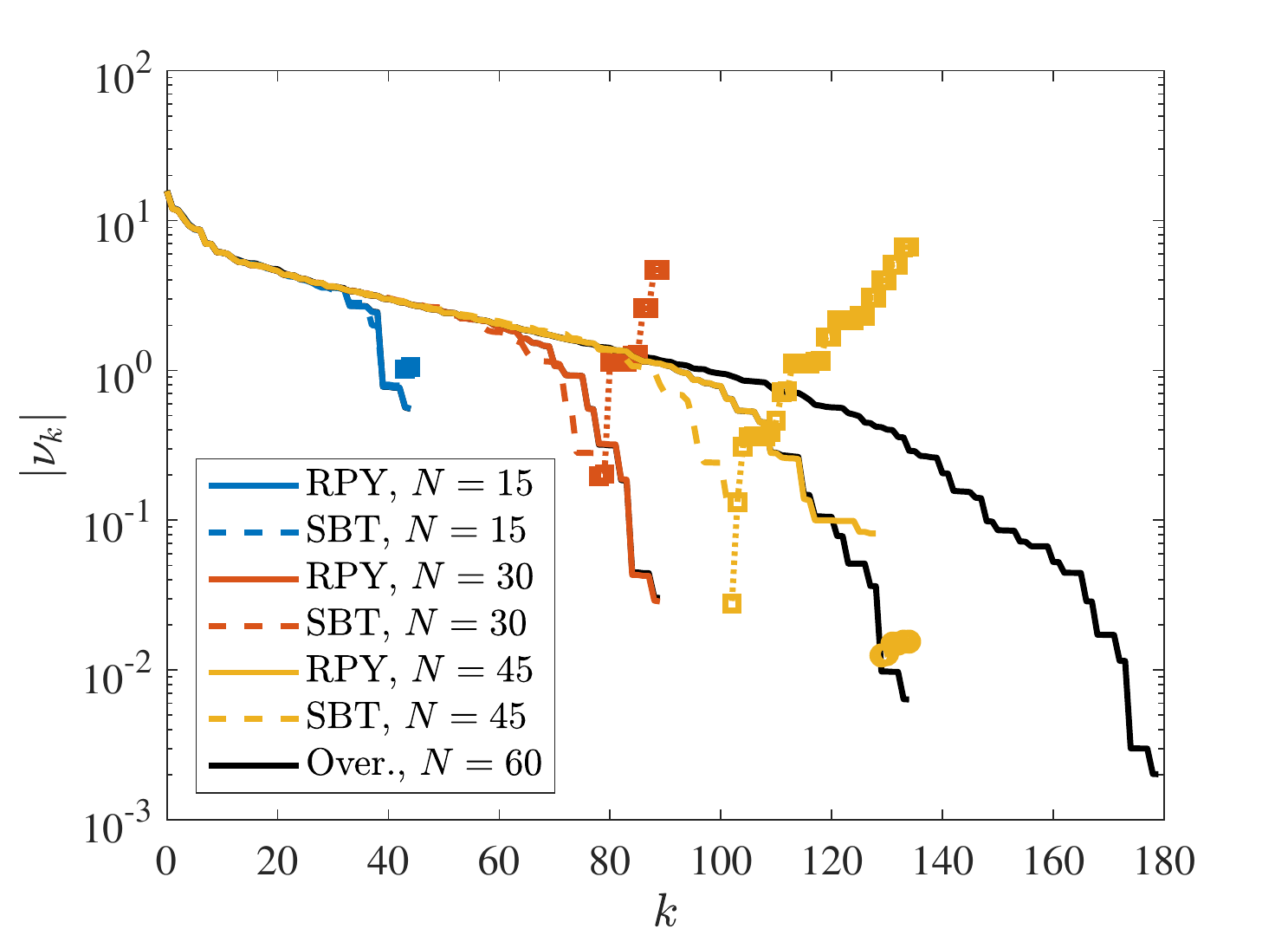}
\includegraphics[width=0.49\textwidth]{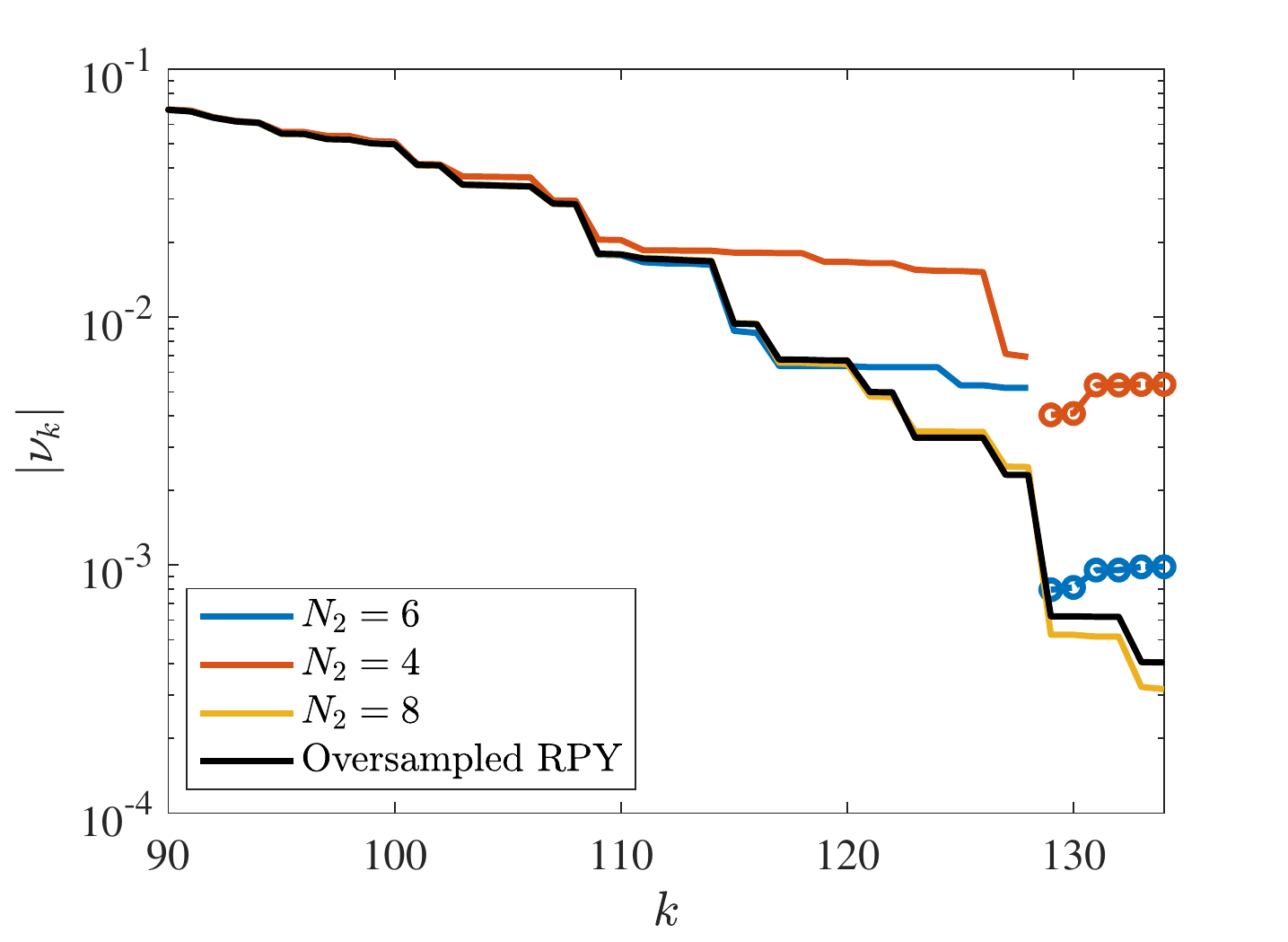}
\caption{\label{fig:MobEigs}Absolute value of the eigenvalues of the discrete mobility matrix $\Mtt$ for $\epsRS=10^{-2}$ with slender-body quadrature for the RPY integrals (positive eigenvalues shown using a solid line and negative using circles), as well as the RPY-based slender body theory of \cite{TwistSBT} (positive eigenvalues shown using a dotted line and negative using squares). Left: the eigenvalues as a function of the number of collocation points $N$ with $N_2=6$ fixed. A result for the RPY integrals\ \eqref{eq:transmob} computed with oversampled quadrature is shown as a black line for each $N=15,30,45$, and 60. There are a few small negative eigenvalues (yellow circles) for $N=45$ due to numerical errors. Right: eigenvalues for $k \geq 90$ for fixed $N=45$ and changing $N_2$  (eigenvalues for $k < 90$ are visually indistinguishable for the lines plotted). Using $N_2=4$ and $N_2=6$ still gives some negative eigenvalues (circles), while $N_2=8$ gives all positive eigenvalues that accurately approximate the true eigenvalues (solid black line).}
\end{figure}

There are two parameters associated with the quadrature for\ \eqref{eq:transmob}: the number of collocation points $N$ that we use for the fiber centerline, which we also use for the Stokeslet and doublet integrals, and the total number of additional points $N_2$ that we use for the integrals in the region $|s-s'| \leq 2\eps$ (recall that we use $N_2/2$ Gauss-Legendre points on either side of $s'=s$ for these integrals). In Fig.\ \ref{fig:MobEigs}, we report the eigenvalues of the mobility matrix for $\epsRS=10^{-2}$ as a function of $N$ (left plot) and $N_2$ (right plot). In the left plot, we also show the eigenvalues of $\Mtt$ when we use the SBT of \cite{TwistSBT} for translation. We see that the eigenvalues with SBT quadrature become negative even using $N=15$ collocation points, while the eigenvalues with the RPY slender-body quadrature remain positive up to $N \approx 35$. When we increase $N_2$ from 6 to 8 (right plot), we match the true eigenvalues of the RPY mobility almost perfectly, and all eigenvalues remain positive. We emphasize that the negative eigenvalues for the RPY quadrature are the result of numerical quadrature errors, while for SBT negative eigenvalues are inherent to the continuum operators. The discrepancy gets better as $\epsRS$ decreases; for $\epsRS=10^{-3}$ there is close agreement in the eigenvalues between SBT and RPY for $N \lesssim 50$, and, even for $N_2=2$, there are no negative eigenvalues until $N \gtrsim 60$ for SBT quadrature and $N \gtrsim 80 $ for RPY quadrature (results not shown). 

\subsubsection{Strong convergence of velocity}
Before considering time-dependent problems, we examine the convergence of the solution to the static problem\ \eqref{eq:euler}. The main object of study here is the constraint force $\V{\lambda}(s)$. Even though our discrete translational mobility $\Mtt$ (with sufficiently large $N_2$) does not have negative eigenvalues (see Fig.\ \ref{fig:MobEigs}), it still has very small positive eigenvalues and the resulting $\V{\lambda}$ changes rapidly near the endpoints. Therefore, we cannot expect \emph{pointwise} spectral convergence for the constraint force $\V{\lambda}$.

The key result of this section is that, while the constraint force $\V{\lambda}$ is not smooth, the velocity for $\epsRS=10^{-3}$ to $10^{-2}$ \emph{is} sufficiently smooth  to be resolved by the spectral method. Furthermore, while the constraint force does not converge pointwise at the fiber endpoints, moments of it, which are the physical observables in the problem (e.g., stress $\sim \int \V{\lambda}(s) \V{X}(s) \, ds$), do converge, and can be captured by the spectral method to reasonably high accuracy. 

To define a reference solution, we utilize the second-order discretization described in Appendix\ \ref{sec:Euler2}. This discretization is more robust because the values of functions at the fiber endpoints do not affect derivatives at the fiber interior, which ensures that $\V{\lambda}$ does not contain spurious oscillations in the fiber interior. Since the mobility in the second-order method converges to the integrals we compute in Section\ \ref{sec:exactmob} as the number of blobs goes to infinity, we will utilize Richardson extrapolation to form a reference solution and compute error with respect to that solution.  \rev{Then, to compare the accuracy of our method with that of the second-order method, we also compute a solution using $1/\epsRS$ points. According to \cite{kallemov2016immersed}, this is (approximately) the minimum required resolution such that the filament is treated hydrodynamically as a solid cylinder, rather than a series of disconnected blobs.}\footnote{\rev{It might be possible to reduce the required resolution of the second-order method by designing second-order quadrature schemes for the integrals in the RPY mobility instead of evaluating direct sums of the RPY kernel. However, since our slender-body quadrature presented in Appendix\ \ref{sec:quads} is based on \emph{global} interpolants of the fiber centerline, it is not straightforward to extend these to second-order discretizations.}} In the second-order method we compute the forces and torques due to bending and twisting analytically, thus focusing our investigation on the saddle point solve\ \eqref{eq:euler}, and not the accuracy of computing the right hand side.

\begin{figure}
\centering
\subfigure[$\epsRS=10^{-2}$]{
\includegraphics[width=0.48\textwidth]{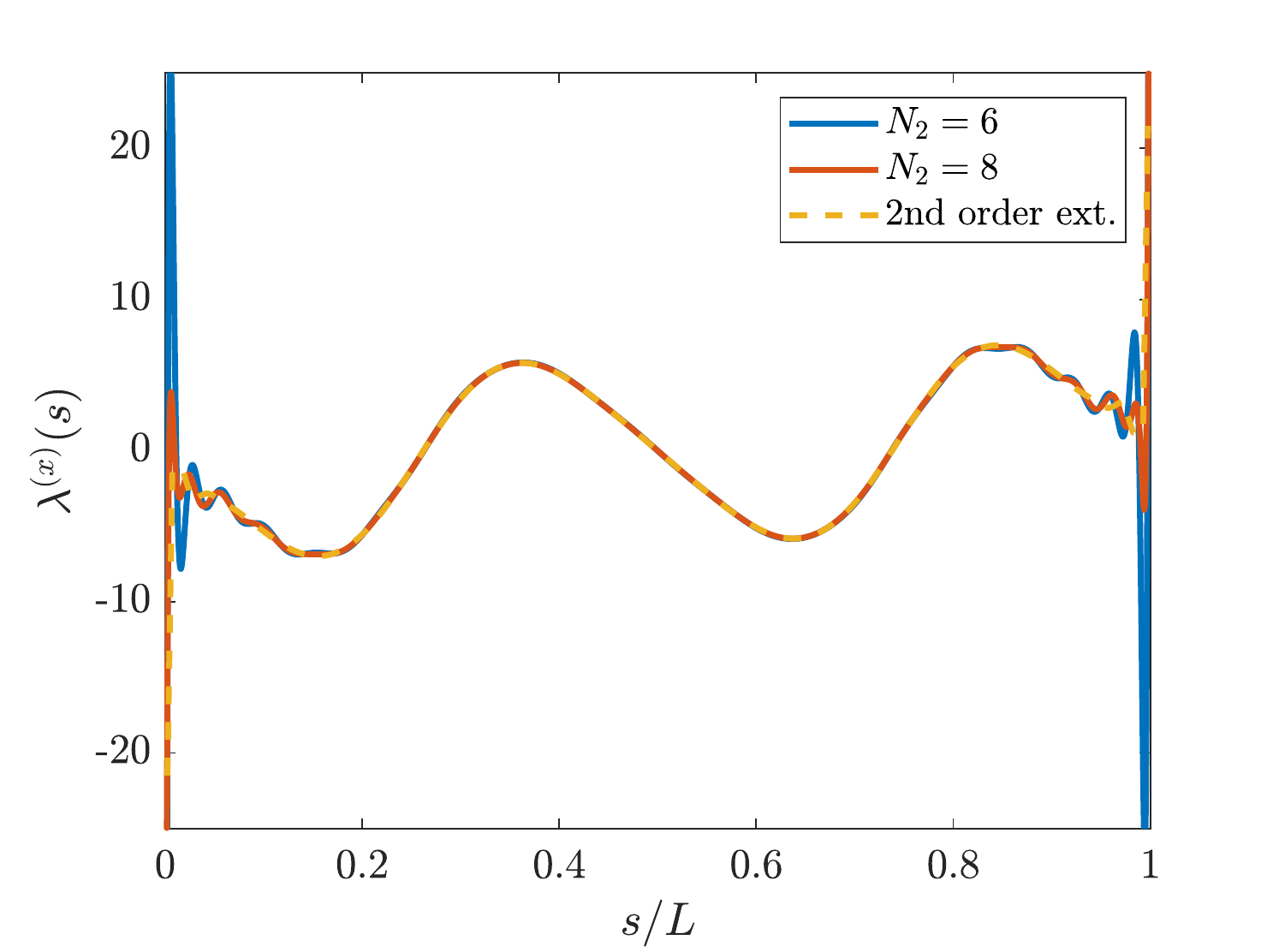}
\includegraphics[width=0.48\textwidth]{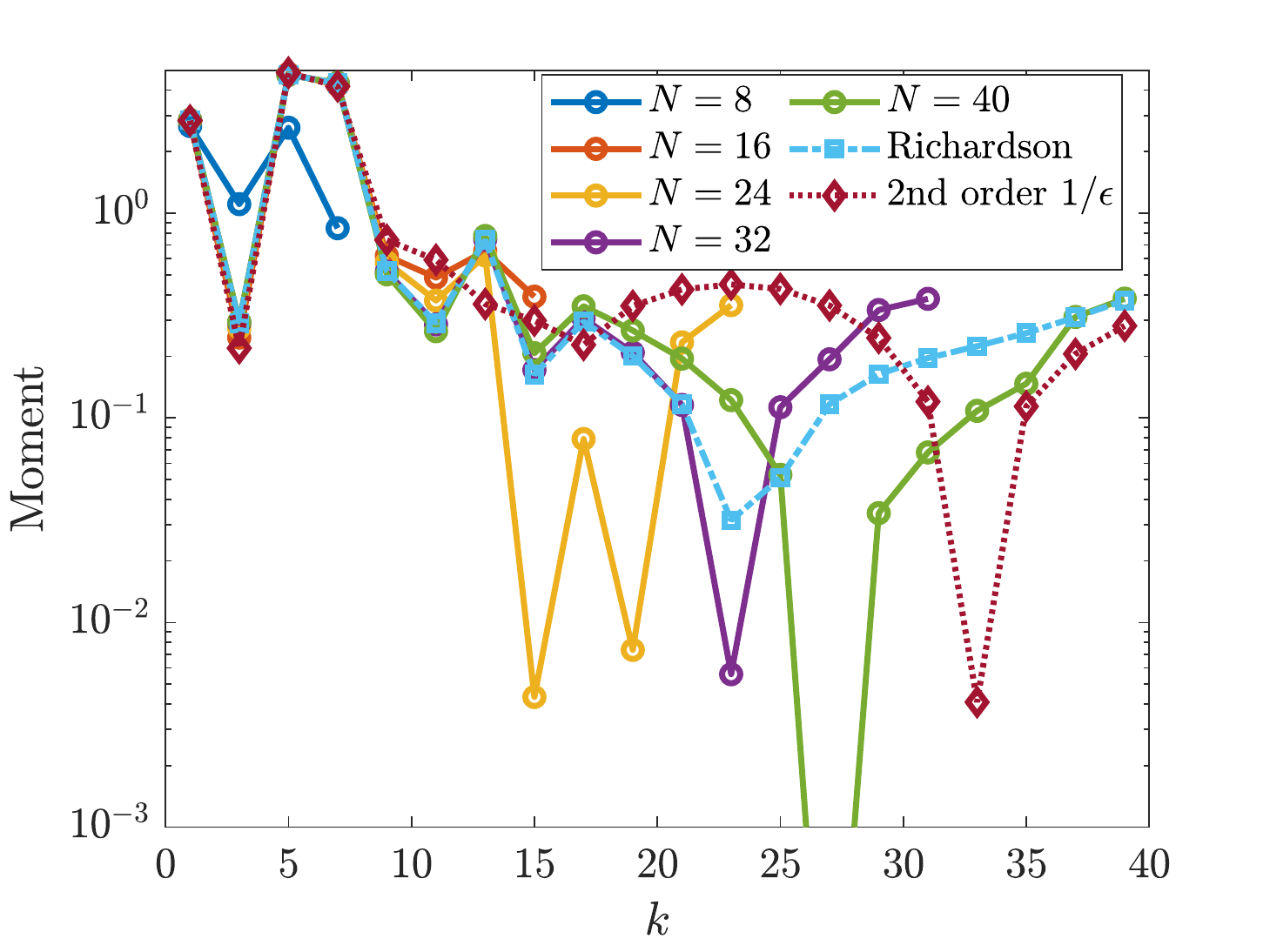}}
\subfigure[$\epsRS=10^{-3}$]{
\includegraphics[width=0.48\textwidth]{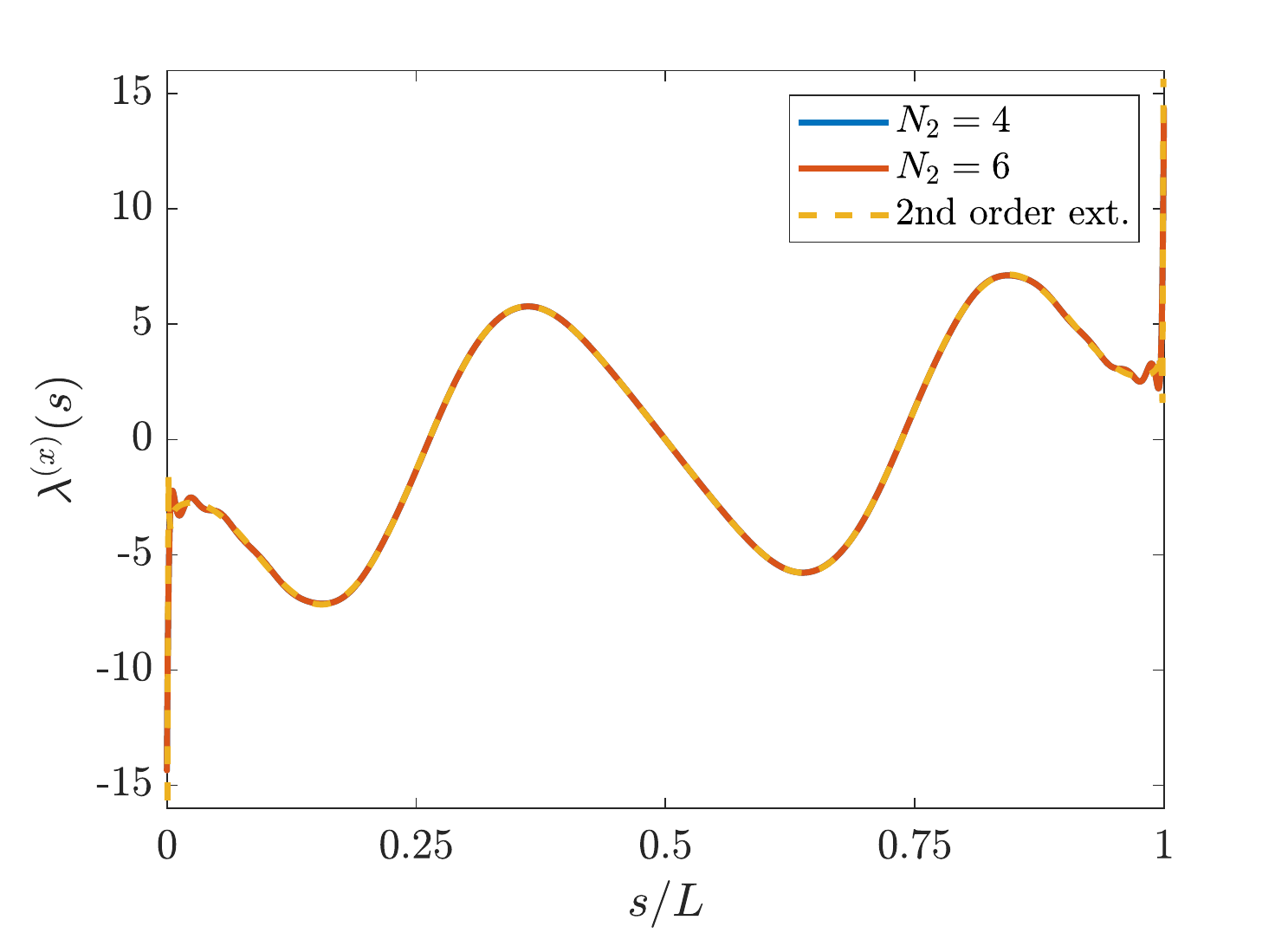}
\includegraphics[width=0.48\textwidth]{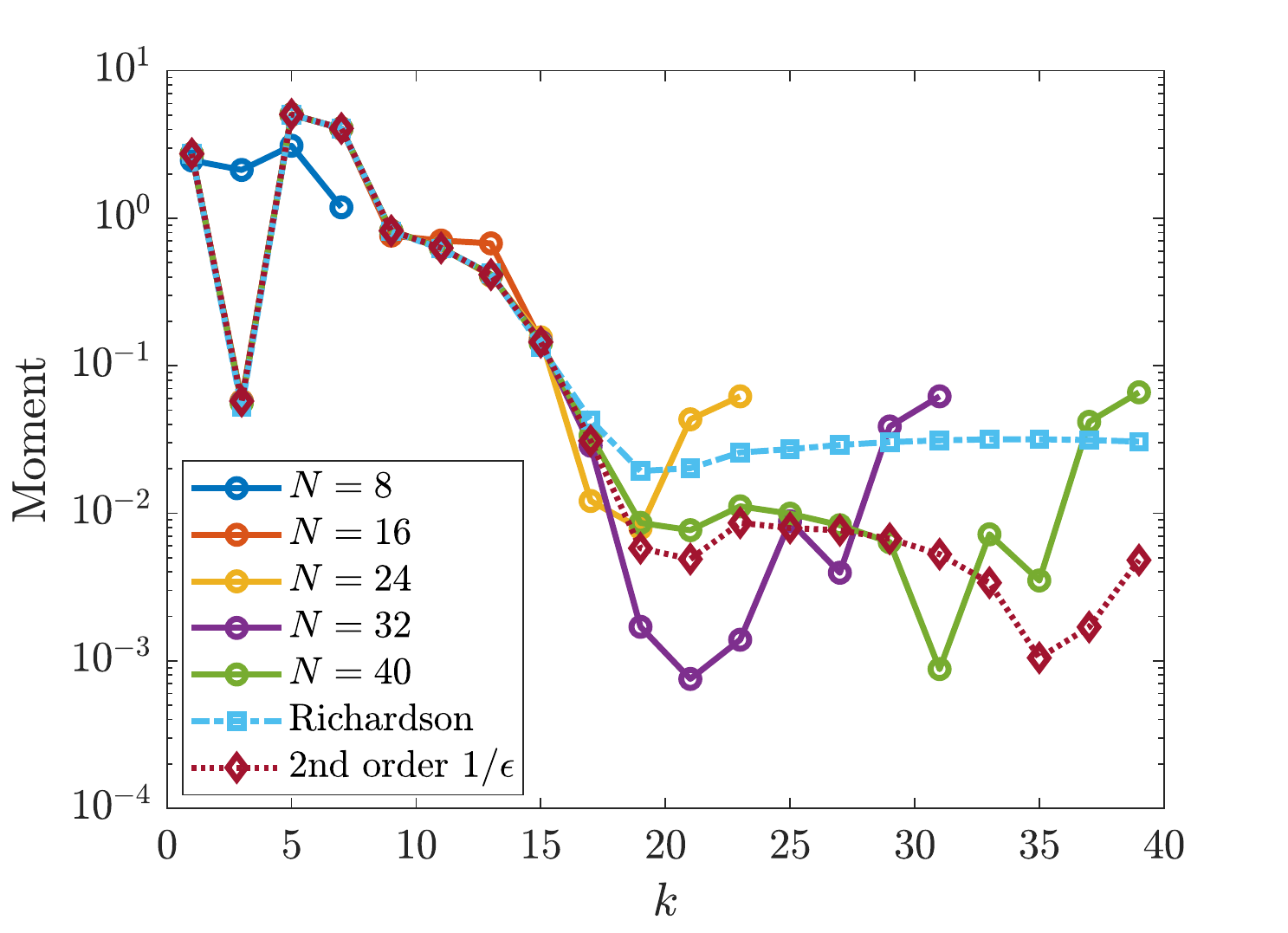}}
\caption{\label{fig:LamMoments}  Comparison of the pointwise values (left) and moments (right) of $\V{\lambda}$ in the spectral and second-order methods using the fiber\ \eqref{eq:Xsfibq} with $q=1$ and slenderness (a) $\epsRS=10^{-2}$ and (b) $\epsRS=10^{-3}$. In the panels on the left, we show the $x$ component, which is the least smooth, for $N=40$ and various $N_2$, and compare to the reference solution obtained by Richardson extrapolation of the second-order discretization. In the panels on the right, we show the magnitude of the moments of $\V \lambda^{(x)}(s)$ against the Chebyshev polynomials $T_k(s)$ with $N_2=8$ for $\epsRS=10^{-2}$ and $N_2=4$ for $\epsRS=10^{-3}$. The spectral method with $N \approx 16$ gives errors in the moments of $\V{\lambda}$ that are smaller than or comparable to the second-order method with $1/\epsRS$ points for both values of $\epsRS$. }
\end{figure}

Let us begin by examining the function $\V{\lambda}(s)$ in the solution of\ \eqref{eq:euler}. \rev{For this test, we begin with} the fiber\ \eqref{eq:Xsfibq} with $q=1$, which is a smoother problem for which we could get high accuracy with on the order of 10 points (and enter the asymptotic regime of the second-order method). We perform Richardson extrapolation on the second-order solutions, and in Fig.\ \ref{fig:LamMoments} compare the result to the spectral method with $N=40$ points. We consider two different values of $N_2$ in both cases. For $\epsRS=10^{-2}$, we start with $N_2=6$, for which we observe uncontrolled oscillations in $\V{\lambda}$ that grow near the fiber endpoints. This is consistent with our observation in Fig.\ \ref{fig:MobEigs} that negative eigenvalues exist when $N_2$ is too small. When we increase to $N_2=8$, we again see oscillations at the endpoints, but the magnitude of these is closer to that of the reference constraint forces we obtain from the second-order method. 

\begin{figure}
\centering
\includegraphics[width=\textwidth]{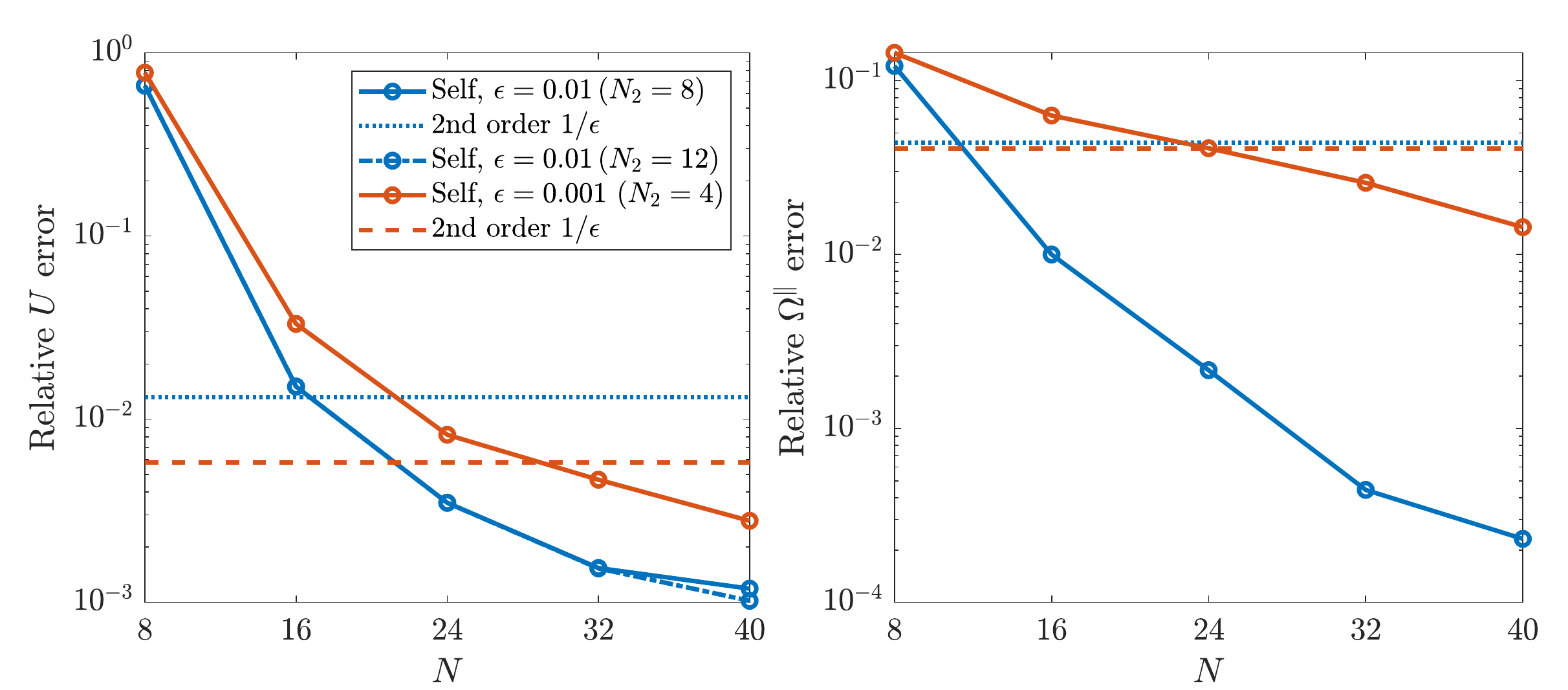}
\caption{\label{fig:EulerSols2} Convergence of the translational and (parallel) rotational velocities for the static Euler model equations\ \eqref{eq:euler} with $\epsRS=10^{-2}$ ($N_2=8$ and $N_2=12$, blue) and $\epsRS=10^{-3}$ ($N_2=4$, orange). For each $\epsRS$, the solid lines (and circles)  show the self-convergence of the spectral method (error is the relative $L^2$ difference between the solution with $N$ points and the solution with $N+8$ points). The dotted (for $\epsRS=10^{-2}$) and dashed (for $\epsRS=10^{-3}$) show the error in the second-order method using $1/\epsRS$ points, relative to Richardson extrapolation. For $\V{U}$, the difference between the Richardson-extrapolated solutions and spectral solutions with $N=40$ is similar to the estimated error when $N=40$. For $\Omega^\parallel$, the Richardson-extrapolated solution differs from the spectral solution with $N=40$ by 0.02 for $\epsRS=10^{-2}$ and 0.01 for $\epsRS=10^{-3}$. \change{Here $\Omega^\parallel$ is computed on an $N$ point grid by applying $\M{R}_\psi$ to\ \eqref{eq:OmparNp2}.}}
\end{figure}

Quite surprisingly, despite the misbehavior of $\V{\lambda}$, the translational velocity $\V{U}$ and parallel rotational velocity $\Omega^\parallel$ appear to converge pointwise, albeit with rapidly-changing behavior at the fiber endpoints (see Fig.\ \ref{fig:PertEffects}, top left plot for an example). To quantify the errors in $\V{U}$ and $\Omega^\parallel$, we perform a self-convergence study, in which the errors are computed by successive refinements, then confirm that the converged spectral solution is close to that of the second-order method (see Fig.\ \ref{fig:EulerSols2} caption). The error as a function of $N$ is shown in Fig.\ \ref{fig:EulerSols2} for both $\epsRS=10^{-2}$ and $\epsRS=10^{-3}$. We obtain about 3 digits of accuracy in both velocities when $N = 32$ and $\epsRS=10^{-2}$. When we increase to $N=40$, the endpoints become more resolved, and the error in the translational velocity using $N_2=8$ points \change{stagnates}. By increasing to $N_2=12$, we \change{obtain a smaller} error when $N=40$ (dashed-dotted squares in Fig.\ \ref{fig:EulerSols2}). In theory, this tells us that $N_2$ should increase with $N$, especially for large $\epsRS$, which is when the contributions from the integrals on $|s-s'| < 2\eps$ are nontrivial. However, for $\epsRS=10^{-2}$, using $N\gtrsim 50$ makes the number of degrees of freedom in the spectral method comparable to the second-order one with $1/\epsRS=100$ points, and in that case we might as well use the more robust low-order method.

When we decrease the slenderness to $\epsRS=10^{-3}$, Fig.\ \ref{fig:EulerSols2} shows that we obtain between $2-3$ digits in both $\V{U}$ and $\Omega^\parallel$ when $N=40$ (the errors saturate at $N\approx 64$). The error in this case is larger than for $\epsRS=10^{-2}$ since the problem is less smooth at the fiber endpoints. While the convergence is slow due to this lack of smoothness, our spectral method still represents a great improvement over the lower-order discretization with $1/\epsRS$ points, as we show in Fig.\ \ref{fig:EulerSols2} with dashed-dotted lines. Thus, if we define the limit of infinitely many blobs as the reference solution, the spectral method provides a cheap, efficient way to approximate that result. In particular, we can obtain the same accuracy with 20 Chebyshev points on the fiber as we do with $1/\epsRS$ blobs. \rev{Unsurprisingly, increasing the fiber shape to $q=7$ in\ \eqref{eq:Xsfibq} makes the required number of Chebyshev points larger; in that case we find that 40 points is sufficient to give a translational velocity with lower error than $1/\epsRS$ blobs, while about 24 points is sufficient for angular velocity (not shown). The translational, but not angular velocity, requires more points as curvature increases because $\Mctt$ is the most nonlocal operator.}

That said, it is still clear from the left panel in Fig.\ \ref{fig:LamMoments} that $\V{\lambda}$ at the endpoints is not converging \emph{pointwise}. Indeed, even the second-order method shows large jumps in $\V{\lambda}$ near the endpoints, which suggests the problem is with the model physics, and not the spectral numerical method. Nevertheless, we see from Fig.\ \ref{fig:LamMoments}(b) that there is still a surprisingly good match between the spectral and second-order solutions for $\V{\lambda}$ at the endpoints for $\epsRS=10^{-3}$, despite the lack of smoothness in $\V{\lambda}$. We also wish to emphasize that, while $\V{\lambda}$ is not smooth at the endpoints, when multiplied by $\Mtt$ it gives \emph{smooth enough} (to resolve with a spectral method to 2--3 digits) contributions to the translational and rotational velocities, which are the important quantities for the evolution of the fiber centerline. The smooth velocity maintains the smooth fiber shape, which makes our spectral method viable.  

\subsubsection{Weak convergence of $\V{\lambda}$ \label{sec:lamweak}}
Given that $\V{\lambda}$ is not smooth at the endpoints, we do not expect pointwise convergence of $\V{\lambda}$, and we cannot a priori expect spectral accuracy from the Chebyshev collocation discretization. We can, however, hope for \emph{weak} convergence of $\V{\lambda}$, which we study by computing the moments $\int_0^L \V{\lambda}(s) T_k(s) \, ds$ for increasing $k$. 

We study this in the following way: for each $N$, we solve for $\V{\lambda}$ and compute the integral against the Chebyshev polynomial $T_k(s)$ for each $k$ on a fine, upsampled grid. We do the same in the second-order method, except we compute the integral directly on the second-order grid without any upsampling. We use Richardson extrapolation of the second-order solutions to get a reference solution for the moments, then compare to the spectral moments in Fig.\ \ref{fig:LamMoments}. The goal for our spectral discretization is then to obtain the first $N$ moments of $\V{\lambda}$ with greater accuracy than the second-order method with $1/\epsRS$ points.  

If we define the extrapolated moments from the second-order moments as the ``true'' moments of $\V{\lambda}$, we see that the spectral method performs better than the second-order method with $1/\epsRS$ points in approximating those moments, especially for $\epsRS=10^{-2}$. For $\epsRS=10^{-2}$, we obtain about $2-3$ digits of accuracy in the first 20 moments using $N \geq 24$, whereas the second-order discretization with 100 points gives only $1-2$ digits in each moment, which is comparable (but worse) than the accuracy of the spectral method with only 16 points. For $\epsRS=10^{-3}$, the error in the moments is comparable for $N \geq 16$ in the spectral method and the second-order code with $1/\epsRS=10^3$ points. For both values of $\epsRS$, the moments of the spectral approximation to $\V{\lambda}$ start to increase for larger $N$ when $k \gtrsim 24$. This comes from the non-smoothness at the endpoints, but has little effect in practice on the accuracy of physical quantities of interest like stress. 

\section{Dynamics \label{sec:dyn}}
In this section we use the numerical method presented in Section\ \ref{sec:specEul} to study the dynamics of bent, twisting filaments. In Section\ \ref{sec:tint}, we first describe a temporal integration scheme in which we treat the bending force and twisting torque implicitly using a linearized backward Euler discretization. We use this scheme to study two examples involving twist elasticity. In Section\ \ref{sec:fibrelax}, we consider the role of twist elasticity in the dynamics of a relaxing bent filament. We verify that twist elasticity contributes a small $\mathcal{O}(\epsc^2)$ correction to the position of the fiber in this case. In Section\ \ref{sec:whirl}, we consider the instability of a twirling filament spinning about its axis at a clamped end. In this case we attempt to bridge the gap between previous experimental \cite{bruss2019twirling} and theoretical \cite{wolgemuth2000twirling, powers2010dynamics} work by showing that rot-trans coupling reduces the critical frequency required to trigger the transition from twirling to whirling by 20\%. \change{We proceed to study the stable ``overwhirling'' dynamics that result when twirling is unstable.}  Matlab code and input files for the examples in this section are available at \url{https://github.com/stochasticHydroTools/SlenderBody}, and supplementary animations are available in the supplementary material. 

\subsection{Temporal integration \label{sec:tint}}
We now develop a temporal integrator for the Euler model\ \eqref{eq:euler}$-$\eqref{eq:PsiEul}. In our previous work without twist \cite[Sec.~4.5]{maxian2021integral}, we introduced a second-order temporal integrator that required only a single saddle-point solve at each time step. This integrator, which is based on using extrapolated values of $\V{X}$ from previous time steps, breaks down in simulations of bundled actin filaments \cite{maxian2021simulations}, in which the hydrodynamic interactions between filaments are important. This led us to switch to a backward Euler discretization for our work on bundled actin filaments \cite{maxian2021bundling}. 

With this in mind, in this work we will also use a linearly-implicit, first-order, backward Euler discretization of the equations. While using higher-order schemes is certainly possible \cite{pareschi2005implicit, delong2013temporal}, this is complicated by the fact that any scheme we use must be $L$-stable since twist equilibrates on a much faster timescale than bending; the higher-order bending modes also equilibrate on a fast timescale. Since $\Mtt$ is a nonlinear function of $\V{X}$, linearizing the mobility as in \cite[Eq.~(25)]{delong2013temporal} does not guarantee stability, which is why prior studies \cite{keavRPY} have used nonlinearly implicit BDF formulas to obtain higher-order accuracy. Designing a temporal integrator for twist suitable for dense suspensions of fibers is a question we leave open for future work; here we only consider a single fiber to first order in time. 

At the $k$th time step, we solve the \emph{linear} system
\begin{gather}
\label{eq:eulerdisc1}
\tdisc{\M{K}}{k} \tdisc{\V{\alpha}}{k} = \tdisc{\Mtt}{k}\left(\V{f}^{(\kappa)} \left(\tdisc{\V{X}}{k+1,*}\right)+\V{f}^{(\gamma)} \left(\tdisc{\V{X}}{k},\tdisc{\psi}{k}\right)+\tdisc{\V{\lambda}}{k}\right)+\tdisc{\Mrt}{k}n^\parallel \left(\tdisc{\psi}{k}\right), \\ \nonumber
\tdisc{\left(\M{K}^*\right)}{k} \tdisc{\V{\lambda}}{k}=\V{0},
\end{gather}
where $\tdisc{\M{K}}{k}=\M{K}\left(\tdisc{\V{X}}{k}\right)$, and likewise for all other matrices. Here $\tdisc{\V{X}}{k+1,*}=\tdisc{\V{X}}{k}+\Delta t \tdisc{\M{K}}{k} \tdisc{\V{\alpha}}{k}$ is an approximation to the position at the next time step. The elastic and twist forces are obtained using\ \eqref{eq:fkapd} and\ \eqref{eq:ftwmodC}, respectively, and the parallel torque is calculated from $\psi$ using\ \eqref{eq:nparBC}. Substituting the formula for $\tdisc{\V{X}}{k+1,*}$ and the discretization of $\V{f}^{(\kappa)}$ from\ \eqref{eq:fkapd} into\ \eqref{eq:eulerdisc1} gives the saddle point system
\begin{gather}
\begin{pmatrix}
-\tdisc{\Mtt}{k} & \tdisc{\M{K}}{k} -\Delta t \tdisc{\Mtt}{k} \M{F}_\kappa \tdisc{\M{K}}{k} \\
\tdisc{\left(\M{K}^*\right)}{k} & \M{0}
\end{pmatrix}
\begin{pmatrix} \tdisc{\V{\lambda}}{k} \\ \tdisc{\V{\alpha}}{k} \end{pmatrix}
= \begin{pmatrix}
\tdisc{\Mtt}{k}\left(\V{f}^{(\kappa)} \left(\tdisc{\V{X}}{k}\right)+\V{f}^{(\gamma)} \left(\tdisc{\V{X}}{k},\tdisc{\psi}{k}\right)\right)+\tdisc{\Mrt}{k}n^\parallel \left(\tdisc{\psi}{k}\right) \\ \V{0}
\end{pmatrix}
\end{gather}
which we solve for $\V{\alpha}^{k}$ and $\V{\lambda}^{k}$ using dense pseudoinverses via the Schur complement \cite[Sec.~4.5.1]{maxian2021integral}. 

We then compute the perpendicular and parallel rates of rotation from the force via the discretization of\ \eqref{eq:OmperpFromK} and\ \eqref{eq:OmParAgain}, 
\begin{gather}
\V{\Omega}^{\perp,k} = \tdisc{\Xs}{k} \times \M{D}_N\tdisc{\M{K}}{k} \tdisc{\V{\alpha}}{k}, \\
\rt{{\Omega}}^{\parallel,k} = \tdisc{\Mrt}{k} \left(\V{f}^{(\kappa)} \left(\tdisc{\V{X}}{k+1,*}\right)+\V{f}^{(\gamma)} \left(\tdisc{\V{X}}{k},\tdisc{\psi}{k}\right)+\tdisc{\V{\lambda}}{k}\right),
\end{gather}
and evolve the twist density $\psi = \ds{\theta}$ via a backward Euler discretization of\ \eqref{eq:PsiEul},
\begin{gather}
\label{eq:psiwrite}
\psi^{k+1} = \psi^k + \Delta t \left(-\V{\Omega}^{\perp,k} \cdot \tdisc{\ds{\Xs}}{k} + \M{D}_N \rt{{\Omega}}^{\parallel,k} + \ds \rr{\Omega^\parallel}\left(\tdisc{\psi}{k+1}\right)\right).
\end{gather}
Here $\tdisc{\ds{\Xs}}{k}$ is computed from $\tdisc{\V{X}}{k}$ by downsampling the second derivative of the upsampled $\tdisc{\widetilde{\V{X}}}{k}$, as discussed in Section\ \ref{sec:BCs}. \change{To compute $\ds{\rr{\Omega}^\parallel}$, we differentiate\ \eqref{eq:OmparNp2} on the $N+2$ point grid, then downsample by applying $\M{R}_\psi$.} This converts\ \eqref{eq:psiwrite} into the linear system
\begin{gather}
\label{eq:psifinal}
\change{\left(\M{I}-\Delta t \M{R}_\psi\M{D}_{N+2} \Mrr\widetilde{\M{N}}_\psi\right)
\tdisc{\psi}{k+1}
= \psi^k + \Delta t \left(-\V{\Omega}^{\perp,k} \cdot \tdisc{\ds{\Xs}}{k} + \M{D}  \rt{{\Omega}}^{\parallel,k} + \M{R}_\psi\M{D}_{N+2}\Mrr \widetilde{\V{\beta}}_n\right),}
\end{gather}
where $\widetilde{\M{N}}_\psi$ and $\M{R}_\psi$ are defined in\ \eqref{eq:Rpsi} and\ \eqref{eq:nparBC}. 

To complete the time step, we rotate the fiber tangent vectors $\Xs$ by $\Delta t \V{\Omega}^{\perp,k}$ to obtain $\tdisc{\Xs}{k+1}$ , then obtain $\tdisc{\V{X}}{k+1}$ by integration, as discussed in detail in \cite[Sec.~4.5.3]{maxian2021integral}. This procedure keeps the fiber strictly inextensible, but is a nonlinear update, which is why we used the approximation $\tdisc{\V{X}}{k+1,*}$ in the saddle-point system\ \eqref{eq:eulerdisc1}.

\subsection{Role of twist in the relaxation of a bent filament \label{sec:fibrelax}}
The first example we consider is a relaxing fiber that is initially bent. Previous works \cite{bergou2008discrete, ehssan17, maxian2021integral}, have neglected the twist elasticity of free filaments, with the justification that twist equilibrates much faster than bend, and therefore is always in a quasi-equilibrium state as the filament deforms. Our goal here is to test this assumption by simulating the relaxation of a filament with and without twist elasticity. 

The fiber has length $L=2$ and the initial tangent vector is given by\ \eqref{eq:Xsfibq} with $q=7$. We start the fiber with no twist, $\psi(s,t=0)=0$, and simulate until $t=0.01$ with $\mu=1$, $\kappa=1$, and $\twmod=0$. This is long enough for the fiber to make a non-negligible shape change, as shown in Fig.\ \ref{fig:FibShapes}. We then repeat the same simulations using a varying twist modulus $\twmod$ and in Fig.\ \ref{fig:RelaxFibDiff} report the $L^2$ error in the fiber position when twist elasticity is ignored (equivalently, $\gamma=0$). We report results using \change{$N=36$} and $\Delta t = 10^{-4}$ for $\epsRS=10^{-2}$ and $N=40$ and $\Delta t = 1.25 \times 10^{-5}$ for $\epsRS=10^{-3}$, having verified that increasing $N$ and decreasing $\Delta t$ makes little difference in the error curves shown in Fig.\ \ref{fig:RelaxFibDiff}. 

\begin{figure}
\centering
\includegraphics[width=0.6\textwidth]{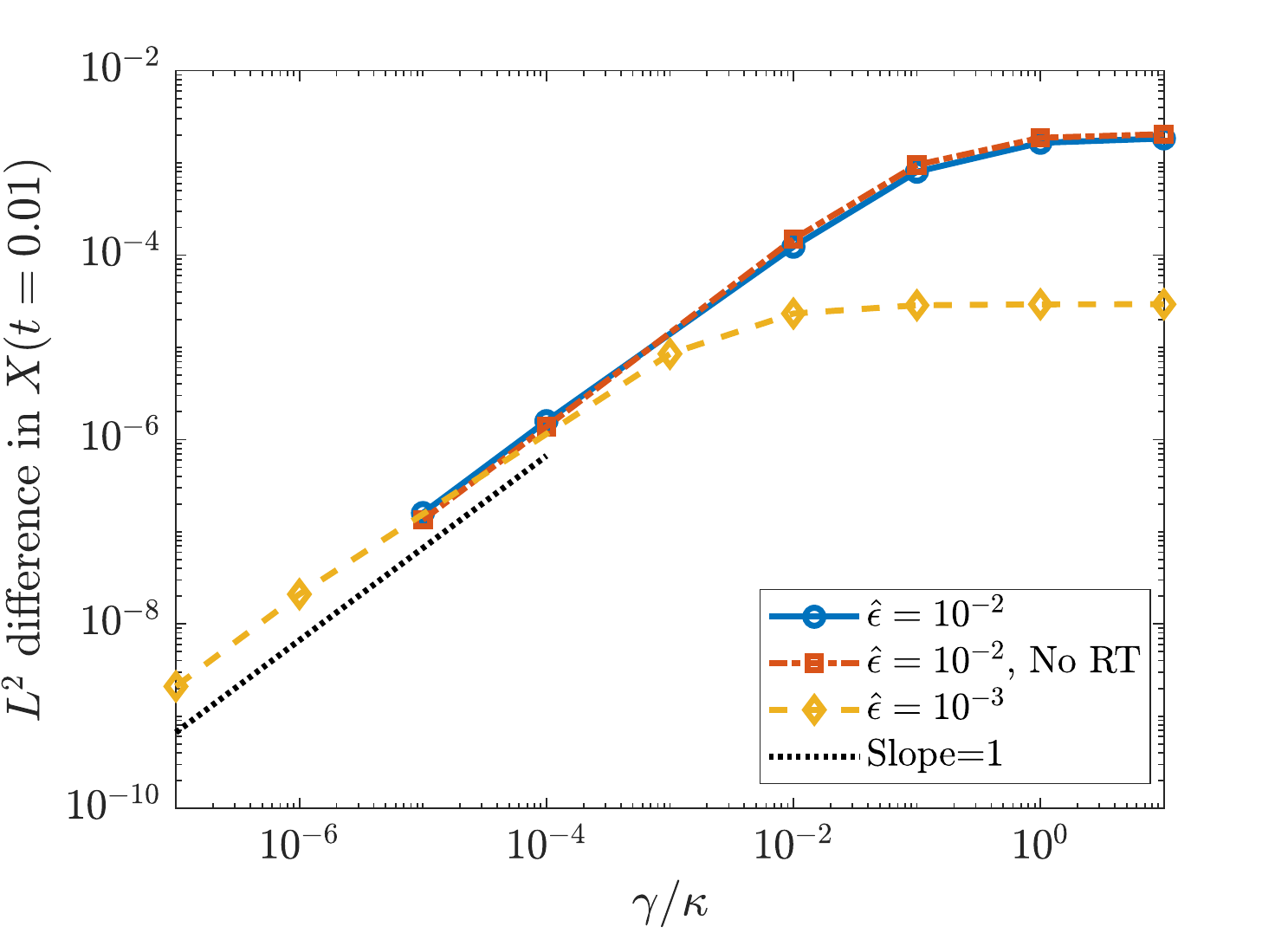}
\caption{\label{fig:RelaxFibDiff} $L^2$ error in the final position of a relaxing bent filament when we do not consider twist elasticity. We measure the error as a function of $\twmod$ and consider $\epsRS=10^{-2}$ (blue circles), $\epsRS=10^{-2}$ without rot-trans coupling (orange squares), and $\epsRS=10^{-3}$ (yellow diamonds). }
\end{figure}

Figure\ \ref{fig:RelaxFibDiff} shows how the position of the filament changes when we add twist elasticity. We recall that the timescale of twist equilibration is $1/\twmod$, and the forcing that results from twist is proportional to $\twmod$. Thus, when $\gamma$ is small, the timescale of twist equilibration is larger than that considered here, and there is an $\mathcal{O}(\twmod)$ correction to the position of the fiber. As $\twmod$ increases, the $1/\gamma$ timescale goes to zero, but the contribution of $\V{f}^{(\twmod)} \propto \twmod$ goes to infinity, so in total there is a constant difference in the position which scales as $\epsRS^2$ in the limit as $\twmod \rightarrow \infty$ (in which twist is always in quasi-equilibrium). 

The scaling of the curves with $\epsRS$ suggests that the rot-rot mobility controls the equilibration of twist, which feeds back to the fiber shape through the force $\V{f}^{(\twmod)}$ (see\ \eqref{eq:euler}) and trans-trans mobility. Indeed, as shown in Fig.\ \ref{fig:RelaxFibDiff}, dropping all rot-trans coupling from the dynamics, so that the \emph{only} coupling between the position and twist angle comes through the term $\Mtt \V{f}^{(\twmod)}$ in\ \eqref{eq:euler}, gives the same behavior for $\epsRS=0.01$. Thus, in the case of a relaxing fiber, the rot-trans coupling is negligible. Furthermore, in most materials, $\gamma/\kappa \approx 1$ \cite{powers2010dynamics}, which according to Fig.\ \ref{fig:RelaxFibDiff} falls within the plateau regime of equilibrated twist, especially for $\epsRS \lesssim 10^{-3}$.


\subsection{Twirling to whirling instability \label{sec:whirl}}
Having considered an example for which twist elasticity is negligible, we now consider an example in which it is vital: the instability in a twirling fiber \cite{wolgemuth2000twirling, lim2004simulations, wada2006non, bruss2019twirling}. In this case, the instability results when the torque due to twist becomes larger than the filament buckling torque, and the filament transitions from a straight twirling state to a curved whirling state \cite{wolgemuth2000twirling}. 

To initialize a small perturbation of a straight fiber, we set
\begin{equation} 
\label{eq:Xswhirl}
\Xs(s) =\frac{1}{\sqrt{1+\delta^2}}\begin{pmatrix} \cos{\left(\tan^{-1}(\delta)\right)} & -\sin{\left(\tan^{-1}(\delta)\right)} & 0 \\
\sin{\left(\tan^{-1}(\delta)\right)} & \cos{\left(\tan^{-1}(\delta)\right)} & 0 \\
0 &0& 1 \end{pmatrix}
\begin{pmatrix}
\delta \cos{\left(s(s-L)^3\right)}\\1 \\\delta \sin{\left(s (s-L)^3\right)}
\end{pmatrix}
\end{equation} 
at $t=0$. Integrating\ \eqref{eq:Xswhirl} numerically (with the pseudo-inverse of the Chebyshev differentiation matrix) and setting $\V{X}(0)=\V{0}$ then gives a filament that satisfies free boundary conditions at $s=L$ and the clamped boundary conditions $\V{X}(0)=\V{0}$ and $\Xs(0)=(0,1,0)$. At the clamped end, there is a delta-like singularity in the perpendicular component of $\V{\lambda}(s)$, which enforces the clamped boundary condition. As a result of this, the velocity has a boundary layer that requires more collocation points to resolve for smaller $\epsilon$. We will therefore use $\epsRS=10^{-2}$ throughout this section, in addition to $\delta=0.01$, $L=2$, $\mu=1$, $\kappa=1$, and $\twmod=1$. 

To trigger the instability, we prescribe the rate of twist at the $s=0$ end as $\Omega^\parallel(0)=\omega$ (see Section\ \ref{sec:BCs}). Note that we could also simulate with the constant torque BC\ \eqref{eq:torqBC}, in which case the applied torque $N_0$ is equal to the torque required to spin a straight fiber at rate $\omega$, 
\begin{gather}
N_0:=-\int_0^L \frac{\omega}{m_\text{rr}(s)} \, ds,
\end{gather}
where $m_\text{rr}$ is the local rot-rot mobility mobility coefficient, $ \rr{\Omega}^\parallel(s)=m_\text{rr}(s)n^\parallel(s)$. The results in this section are unchanged regardless of the BC we use. 

We start the twist profile in a state that satisfies the boundary conditions with no rot-trans coupling, $\psi(s)= \omega/\twmod \int_L^s m_\text{rr}^{-1}(s') \, ds'$. In the case of an ellipsoidally-tapered filament we use\ \eqref{eq:UinRRpar} along the whole filament, so $m_\text{rr}(s) \equiv 9/\left(32\pi \mu \eps^2\right)$, and the steady state twist profile reduces to $\psi(s)=\omega(s-L)/(\psi m_\text{rr})$ \cite{wolgemuth2000twirling}. When the fiber is cylindrical, we use the formulas in \cite{TwistSBT} to obtain $m_\text{rr}(s)$. 

In the case of an ellipsoidal fiber with local drag and no rot-trans coupling, Wolgemuth et.\ al \cite{wolgemuth2000twirling} performed a linear stability analysis to show that the spinning of the fiber about its axis is unstable at critical frequency \cite[Eq.~(74)]{powers2010dynamics}
\begin{equation}
\label{eq:omegac}
\omega^\text{(ELD)}_c \approx 8.9 \frac{m_\text{rr}\kappa}{L^2} = 8.9\frac{9\kappa }{32 \pi \mu L^2 \eps^2}
\end{equation}
with the centerline oscillating at a significantly smaller frequency \cite[Eq.~(75)]{powers2010dynamics}
\begin{equation}
\label{eq:chic}
\chi^\text{(ELD)}_c \approx 22.9 \frac{m_\text{tt}^\perp\kappa}{ L^4}  = 22.9 \frac{\kappa \left(\log{\left(\epsRS^{-2}/16\right)}+4\right)}{8\pi\mu L^4}.
\end{equation}
Here $m_\text{tt}^\perp$ is the mobility coefficient for force perpendicular to the fiber centerline for an ellipsoidally-tapered filament. We will therefore report frequency in units $\bar{\omega}=\omega/\omega^\text{(ELD)}_c$ and time in units $\bar{t}=t\chi^\text{(ELD)}_c/(2\pi)$.

Since the publication of \cite{wolgemuth2000twirling}, there have been a number of studies on the whirling instability that leave a few questions open. A study by Lim and Peskin, for instance \cite{lim2004simulations}, found a critical frequency of $\bar{\omega}_c=0.3$, while recent experimental work \cite{bruss2019twirling} has put the critical frequency at about $\bar{\omega}_c=0.6$. The discrepancy in Lim and Peskin's work \cite{lim2004simulations} could be due to a number of factors, including their use of nonlinear fluid equations, a periodic domain, and a larger initial deflection. It is not, however, due to nonlocal trans-trans interactions, as Wada and Netz \cite{wada2006non} included those and obtained $\bar{\omega}_c=1$. They also showed, however, that the frequency required to induce the overwhirling state drops as the initital perturbation increases (and the nonlinear regime is entered). It could be possible, therefore, that the simulations of \cite{lim2004simulations} fall in this nonlinear regime. This observation was used in \cite{bruss2019twirling} to explain the deviations from the theory as well.

An alternative explanation for the different values of $\omega_c$ is the influence of rot-trans coupling. This is neglected in computational studies which use the RPY tensor or SBT \cite{wolgemuth2000twirling, wada2006non}, but is included by necessity in simulations that use the IB method \cite{lim2004simulations}, and of course in experimental systems \cite{bruss2019twirling}. Here we will study the role of rot-trans coupling to see if we can understand the discrepancies in the literature. We will approach the problem in two steps. First, we will confirm the relationships\ \eqref{eq:omegac} and\ \eqref{eq:chic} hold for an ellipsoidal filament with local drag and a cylindrical filament with all hydrodynamic terms included except rot-trans. Then, we will examine the influence of the rot-trans coupling in the linear regime. 

\subsubsection{Critical frequency as a function of hydrodynamic model \label{sec:wcrit}}

\begin{figure}
\centering
\subfigure[Growth at $\omega=\omega_c^\text{(ELD)}$]{
\includegraphics[width=0.6\textwidth]{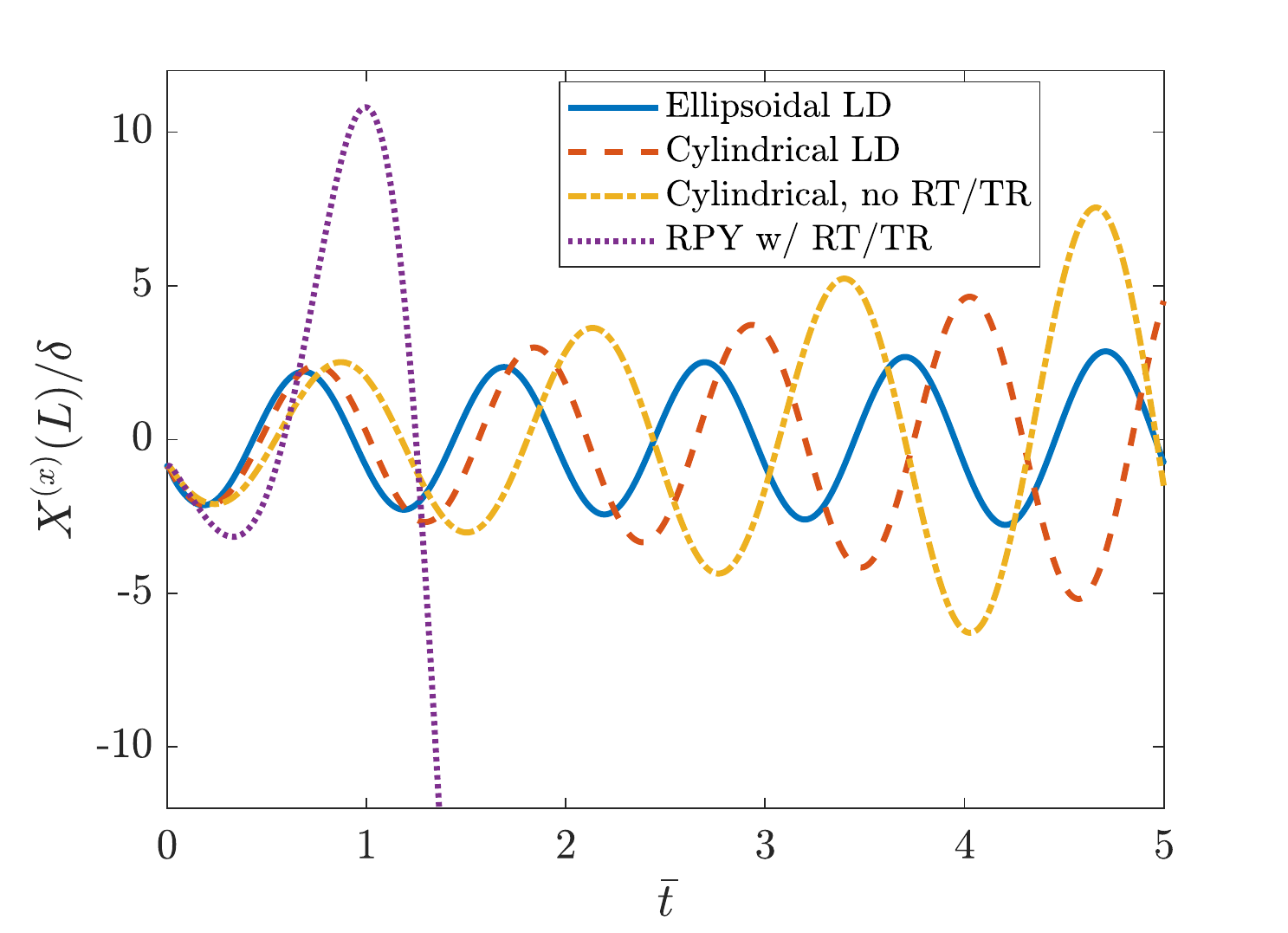}}
\subfigure[Values of $\omega$ with no growth]{
\includegraphics[width=0.6\textwidth]{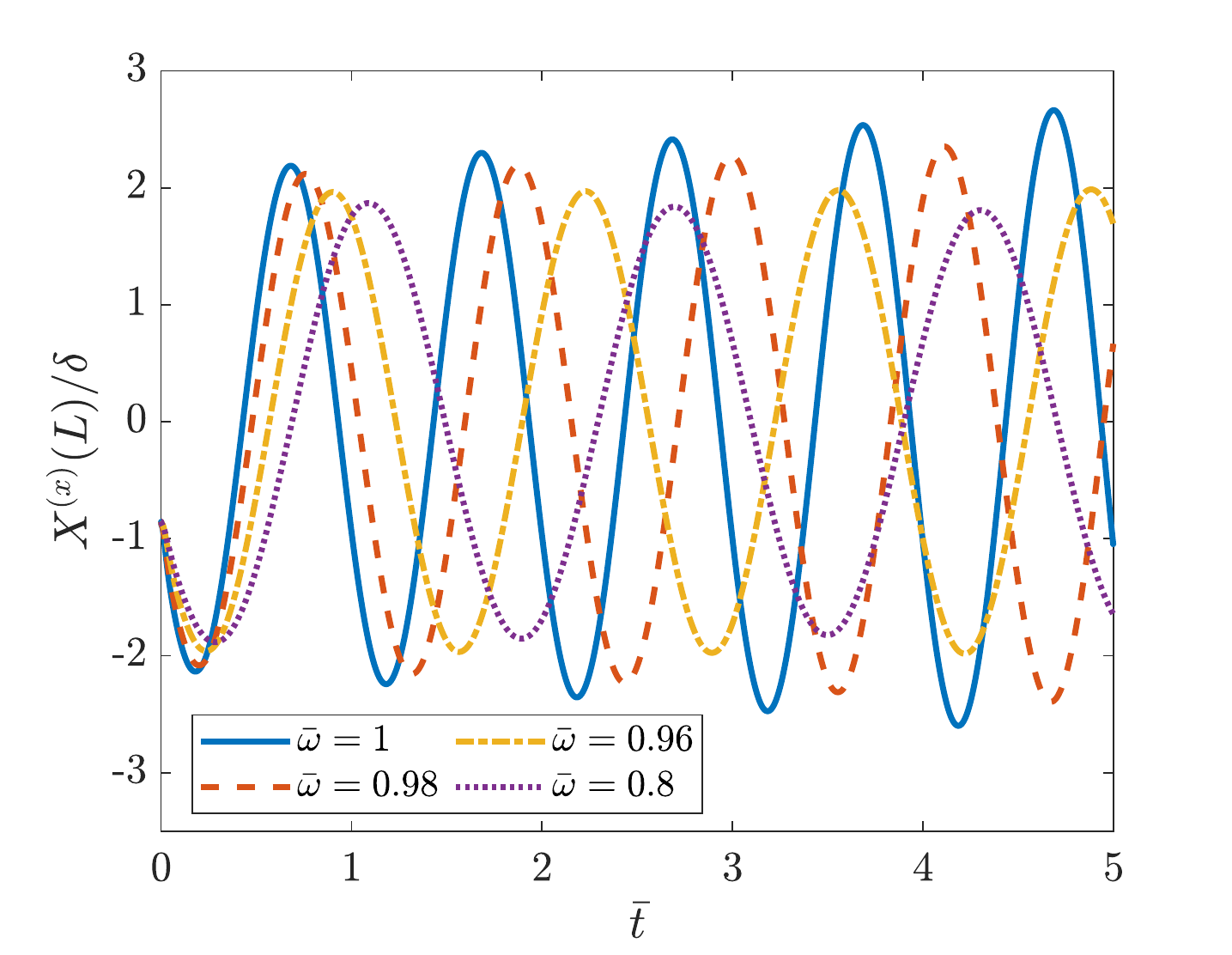}}
\caption{\label{fig:WhirlingLinear} Simulating the whirling instability in the linear regime, $\delta=0.01$ in\ \eqref{eq:Xswhirl}. We show the $x$ coordinate of the $s=L$ endpoint as a function of dimensionless time $t^*=t\chi^\text{(ELD)}_c/(2\pi)$ for $\epsRS=10^{-2}$ (results for $\epsRS=10^{-3}$ are very similar, but require larger $N$ to capture the changes in velocity near the endpoints, and are not shown). We consider ellipsoidal (blue) and cylindrical (dashed red) fibers with local drag, and cylindrical fibers (dashed-dotted yellow), where in all three of these cases we neglect rot-trans and trans-rot coupling. We then add rot-trans and trans-rot coupling and show the results as dotted purple lines. In (a), we show the endpoint $x$ coordinate at $\omega=\omega^\text{(ELD)}_c$ using $1000$ time steps per period $2\pi/\chi_c^\text{(ELD)}$, $N=16$ for ellipsoidal fibers, and $N=32$ for all other fibers. In (b) we tune $\bar{\omega}=\omega/\omega^\text{(ELD)}_c$ to be as close at possible to the frequency at which the perturbation neither grows nor decays (each line represents the same mobility as the top plot), and use $2000$ time steps per period with $N=40$ for all fibers.}
\end{figure}

Figure\ \ref{fig:WhirlingLinear} shows our study of the behavior close to the critical frequency $\omega_c$. In Fig.\ \ref{fig:WhirlingLinear}(a), we simulate at $\bar{\omega}=1$ with four different hydrodynamic models: local drag with ellipsoidal fibers (no rot-trans coupling), local drag with cylindrical fibers (endpoint formulas in \cite[Appendix~C.1]{TwistSBT}, no rot-trans coupling), and the RPY integral mobility with and without rot-trans coupling. We see that using ellipsoidal local drag gives a trajectory which neither grows nor decays over five periods, and the period is exactly equal to that predicted from\ \eqref{eq:chic}. When we switch to the cylindrical local drag formulas, the mobility at the fiber endpoints decreases, and so the period $T_c \sim \chi_c^{-1} \sim m_\text{tt}^{-1}$ increases (see\ \eqref{eq:chic}). Likewise, the rot-rot mobility $m_\text{rr}$ at the fiber endpoints is also smaller when we use cylindrical fibers, and so the critical frequency $\omega_c$ ought to be smaller (see\ \eqref{eq:omegac}). As shown in Fig.\ \ref{fig:WhirlingLinear}(b), this is indeed the case, with the critical frequency for cylindrical filaments coming in lower than that for elliposidal filaments, but with a small difference of only 2\%. The same can be said for the RPY mobility without rot-trans coupling, since in this case the difference in the critical frequencies is at most 4\%. Thus, nonlocal trans-trans hydrodynamics, as well as the shape of the fiber radius function, make little difference for the critical frequency, and cannot explain the previously-observed large deviations from $\omega_c^\text{(ELD)}$. 

It is only when we account for rot-trans coupling that we see a substantial difference in the critical frequency. Indeed, as shown in Fig.\ \ref{fig:WhirlingLinear}, including trans-rot dynamics gives a large growth rate of the oscillations in the fiber endpoints, which eventually leads to the ``overwhirling'' behavior that has been documented previously \cite{lim2004simulations, powers2010dynamics} (see Fig.\ \ref{fig:WhirlingPic}, but note that the same overwhirling behavior can be obtained without rot-trans and trans-rot dynamics for $\omega > \omega_c$ \cite{wada2006non}). Figure\ \ref{fig:WhirlingLinear} shows that the critical frequency with rot-trans coupling is about $\bar{\omega}_c=0.8$, which explains about half of the observed experimental value of 0.6. The critical frequency of $\bar{\omega}_c=0.8$ is also unchanged within 2\% when we place the clamped end $2\eps=0.04$ above a bottom wall and use the wall-corrected RPY mobility derived in \cite{swan2007simulation}, \rev{which we discretize using oversampled integrals at every collocation point}. Thus, confinement, at least in the direction perpendicular to the centerline, is not the cause of the other 20\%, since the wall corrections decay too fast away from the wall to have a noticeable effect (the total effect is $\mathcal{O}(\epsRS)$ \cite[Sec.~2.3--4]{de1975low}). 

The rest of the deviation could be due to the initial fiber shape, as is suggested in \cite{bruss2019twirling}; however, our observations are the same up to a deflection of $\delta \approx 0.2$ in\ \eqref{eq:Xswhirl}, and Wada and Netz showed that starting the filament in a circular configuration gives about 25\% reduction in $\omega_c$ \cite[Fig.~1(g)]{wada2006non}. Seeing as the initial fiber shapes in \cite{bruss2019twirling} are far more straight than circular, it seems unlikely than a reduction as large as 20\% comes from the fiber shape. It is possible that the confinement inside a PVC pipe could reduce the critical frequency, although the authors of \cite{bruss2019twirling} say that the flow is fully developed around the filament. \rev{Other possibilities include the experimental set-up; for instance, the supplementary videos of \cite{bruss2019twirling} show that the shaft also slightly translates the fiber endpoint in addition to spinning it. This provides additional perturbations that could lower the critical frequency.}

\subsubsection{\change{Nonlinear overwhirling dynamics}}

\begin{figure}
\centering
\includegraphics[width=0.48\textwidth]{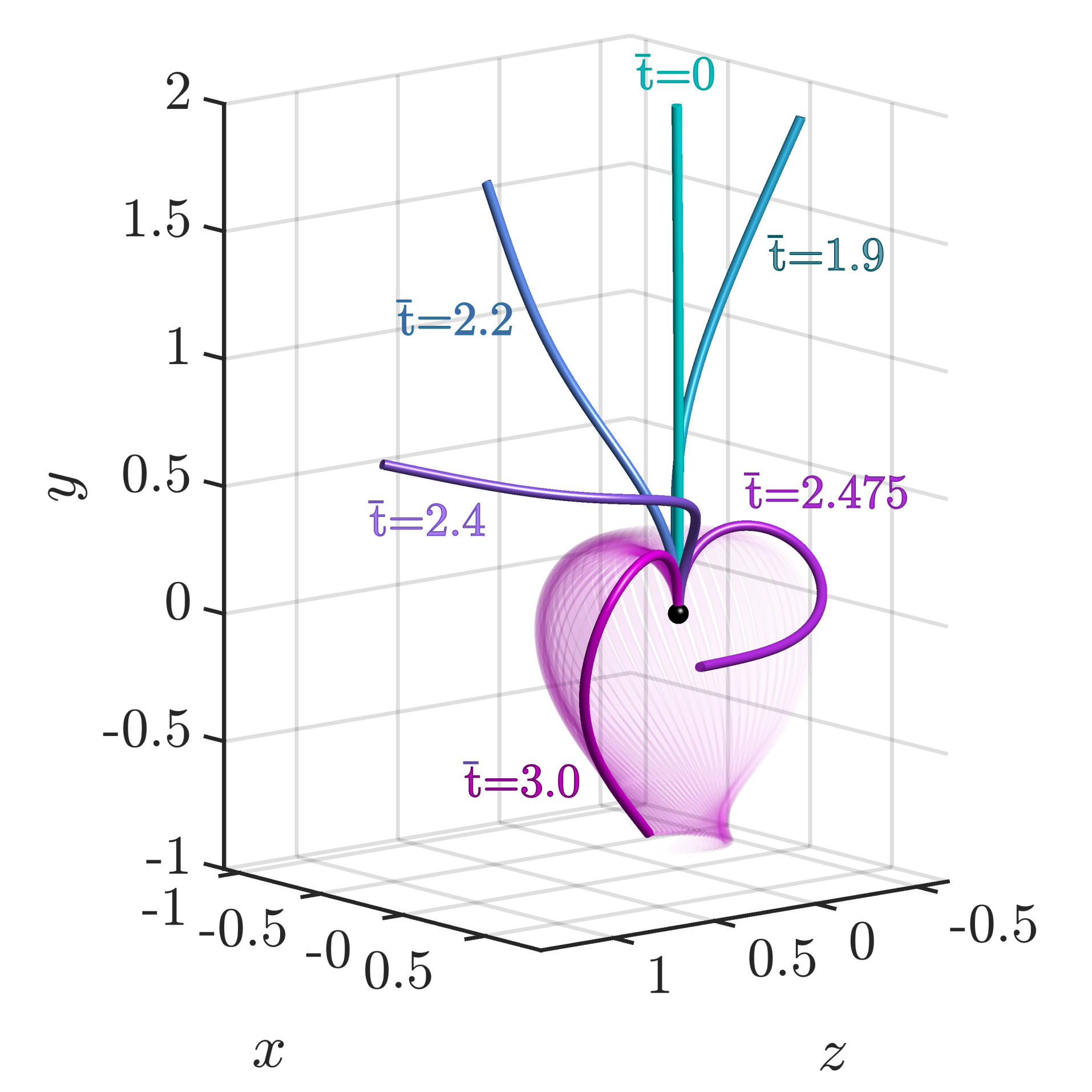}
\caption{\label{fig:WhirlingPic} \change{Nonlinear overwhirling dynamics: Fiber trajectory when $\omega=\omega_c^{\text{(ELD)}}=497$ rad/s $\approx 1.25 \omega_c$ and we include rot-trans coupling in the RPY mobility (dotted purple curve in Fig.\ \ref{fig:WhirlingLinear}). The fiber endpoint drops below $y=0$ and begins periodic crankshafting motion around $\bar{t} \approx 2.5$. }}
\end{figure}

\change{As $\omega$ increases, the twisting torque overcomes the bending moment of the filament, and the filament dynamics eventually reach a steady state in which the fiber takes the shape of a shepherd's crook and spins in a steady motion, as illustrated in Fig.\ \ref{fig:WhirlingPic}. These dynamics, which are called ``overwhirling'' \cite{lim2004simulations}, have been studied in more detail in \cite{wada2006non, lee2014nonlinear}. In particular, it was shown in \cite{lee2014nonlinear} that there is a subcritical Hopf bifurcation: if $\omega$ is increased from below, the overwhirling frequency is $\omega=\omega_c$, but if $\omega$ is decreased once the fiber is already in the overwhirling state, overwhirling is maintained for $\omega < \omega_c$. In addition to determining the boundaries of the bistable region, in this section we also study the amplitude and frequency of the periodic overwhirling motion as a function of the motor twirling frequency $\omega$. 

To study the stability and behavior of overwhirling at various frequencies, we begin by getting the fiber to the steady overwhirling state using $\omega=1.1\omega_c$ (by steady, we mean the crankshafting frequency is unchanged to $0.2\%$ over the last 5 periods). Following this, we use this steady state configuration to initialize another simulation with $\omega=1.2\omega_c$, then use that steady state to initialize $\omega=1.3 \omega_c$, and so forth. We do the same for decreasing frequency: starting from $\omega=1.1\omega_c$, we decrease the frequency until the fiber comes out of the overwhirling state. As discussed in \cite{lee2014nonlinear}, this happens at a \emph{lower} frequency than $\omega_c$, which we find to be somewhere between $0.8\omega_c$ (stable overwhirling) and $0.72\omega_c$ (stable twirling). Thus, the bistable region in the bifurcation diagram is between $\omega/\omega_c=0.8$ and $\omega/\omega_c=1$, which is (approximately) consistent with the results of \cite[Fig.~5]{lee2014nonlinear}, where the lower boundary appears to be $0.85-0.9\omega_c$. }

\begin{figure}
\centering
\subfigure[Dynamics of fiber endpoint]{
\includegraphics[width=0.8\textwidth]{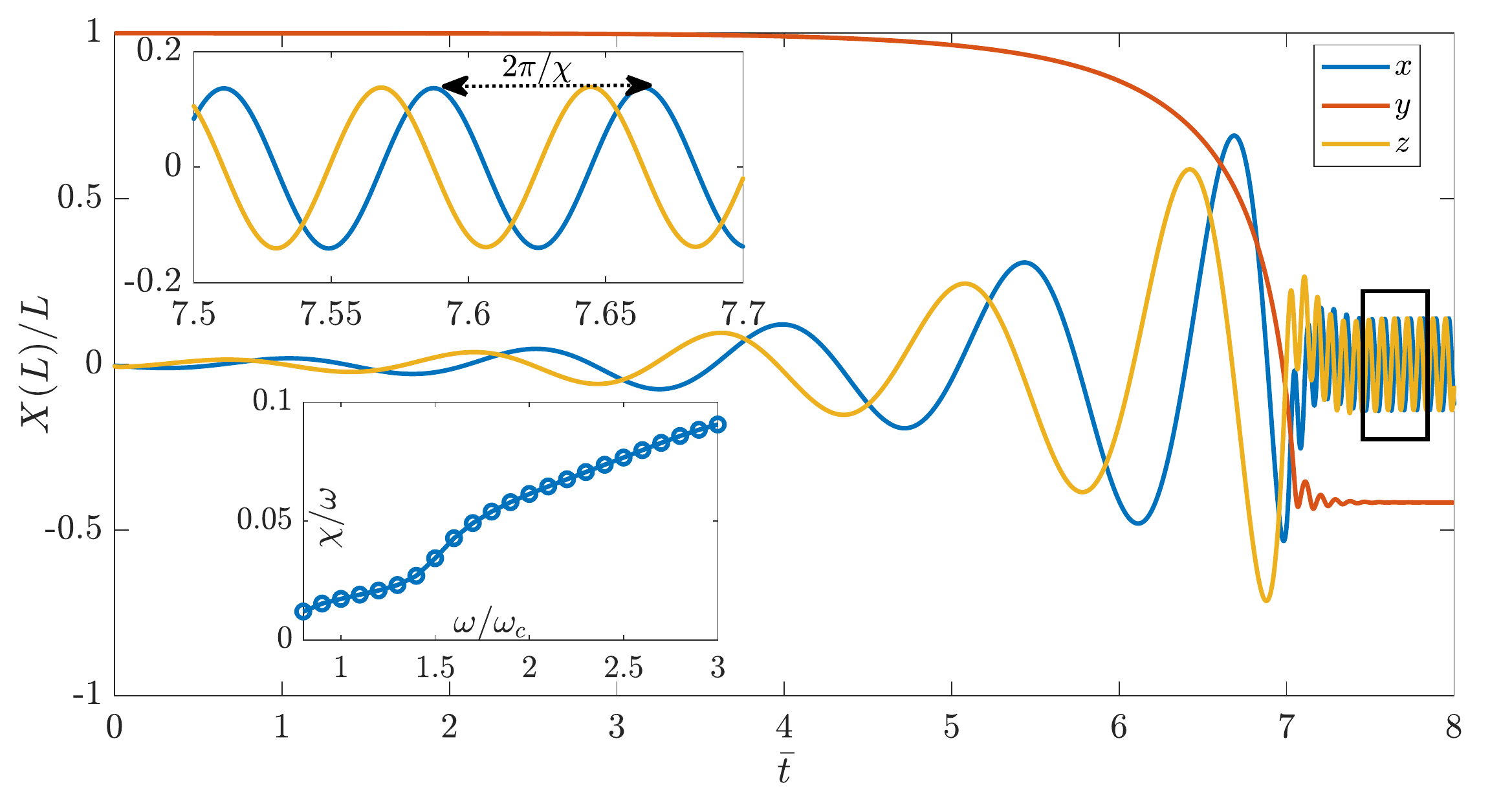}}
\subfigure[Steady state amplitude]{
\includegraphics[width=0.48\textwidth]{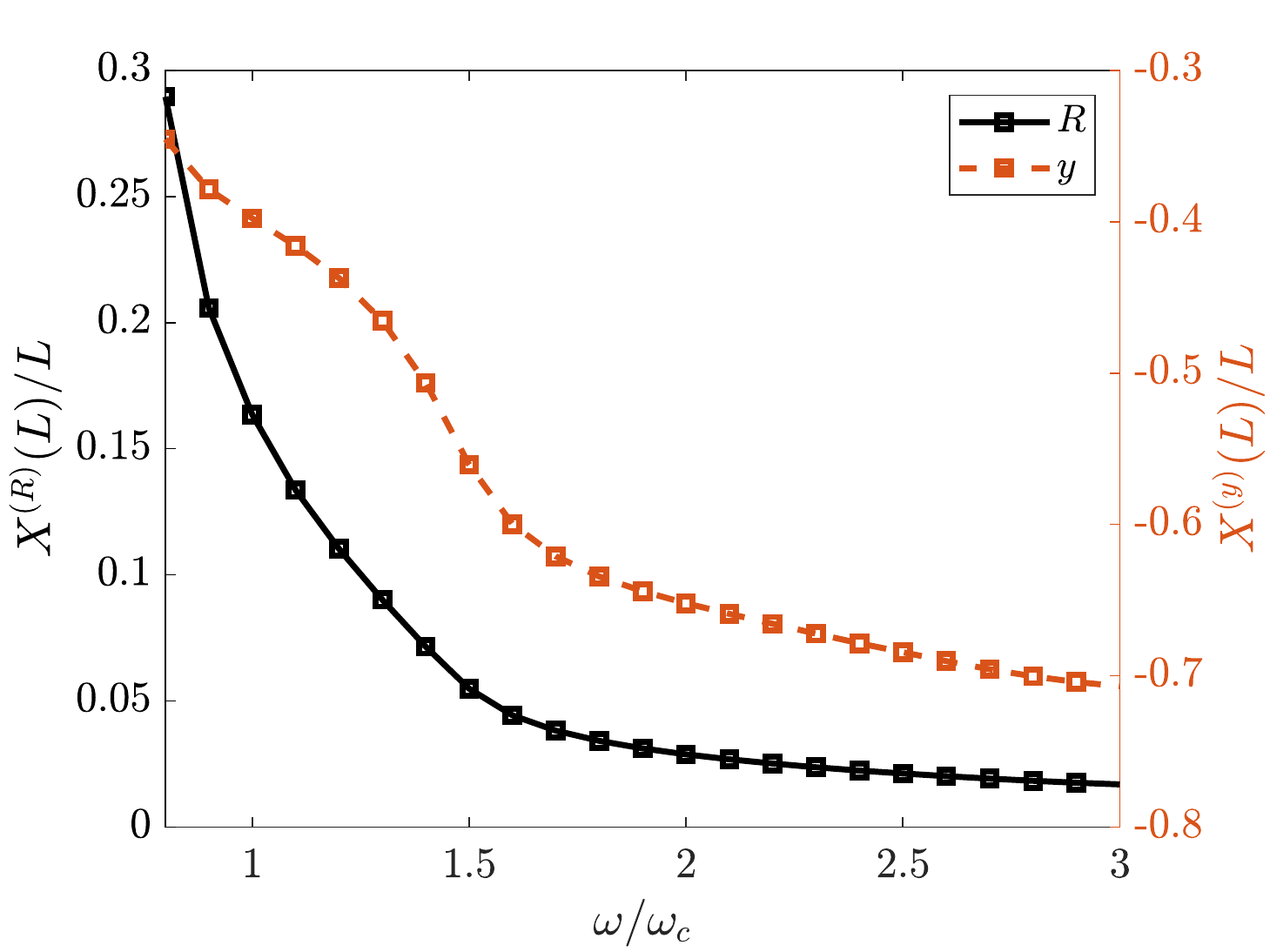}}
\subfigure[Steady state shapes]{
\includegraphics[width=0.35\textwidth]{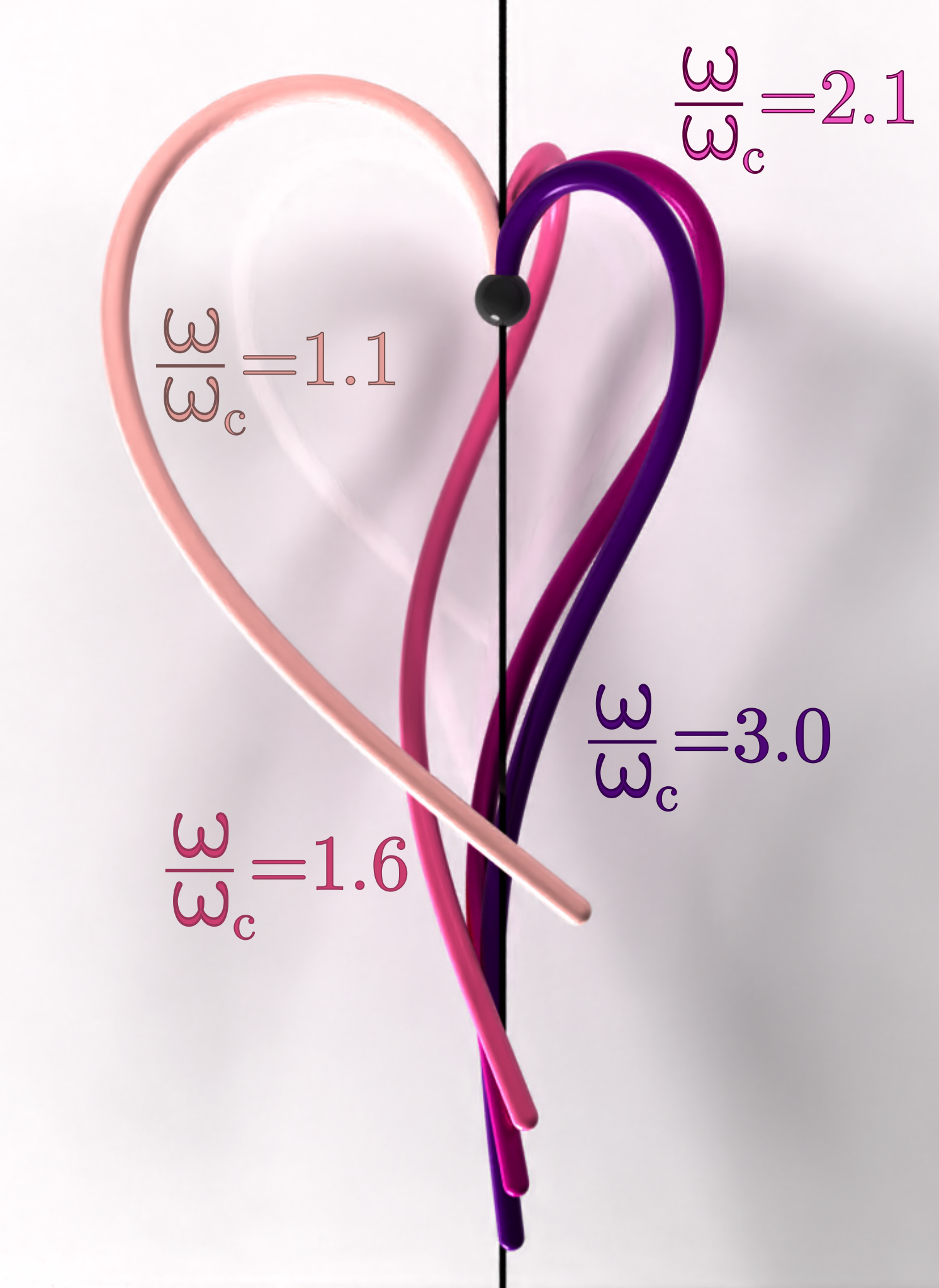}}
\caption{\label{fig:Whirling2} \change{Steady state overwhirling dynamics. (a) Transition to the overwhirling state for $\omega=1.1\omega_c$; see Fig.\ \ref{fig:WhirlingPic} for a 3D illustration for a different $\omega/\omega_c=1.25$. The fiber endpoint dips below $y=0$ and begins crankshafting motion with frequency $\chi(\omega)$. The top inset zooms in on the black box, while the bottom inset shows $\chi(\omega)$. (b) Steady state amplitudes for varying $\omega/\omega_c$. The left axis and black line shows the amplitude $R=\sqrt{X(s)^2+Z(s)^2}$, while the right axis and red line show the steady state $y$ coordinate. (c) Steady state overwhirling configurations for given $\omega/\omega_c$.}}
\end{figure}

\change{The overwhirling steady state dynamics are analyzed in Fig.\ \ref{fig:Whirling2}. In Fig.\ \ref{fig:Whirling2}(a), we show the trajectory of the $s=L$ endpoint when we use $\omega/\omega_c=1.1$ and initialize the fiber to be nearly straight using\ \eqref{eq:Xswhirl}. As depicted in Fig.\ \ref{fig:WhirlingPic}, the endpoint dips below $y=0$, and the whole fiber begins a steady crankshafting rigid rotation around the $y$ axis, such that $Y(s)$ and $\sqrt{X(s)^2+Z(s)^2}$ are constant for all $s$. In Fig.\ \ref{fig:Whirling2}, we show how the amplitude $R=\sqrt{X(L)^2+Z(L)^2}$ (Fig.\ \ref{fig:Whirling2}(b)) and frequency $\chi$ (inset of Fig.\ \ref{fig:Whirling2}(a)) of the crankshafting motion depend on the driving motor frequency. For $\omega/\omega_c \lesssim 1.6$, the amplitude contracts rapidly to about 5\% of the fiber length, while the crankshafting frequency is roughly constant at $1-2$\% of the motor frequency. On $\omega/\omega_c \gtrsim 1.6$, however, the fiber amplitude function is roughly constant; indeed, Fig.\ \ref{fig:Whirling2}(c) shows that the steady state overwhirling shapes approach a fiber that is essentially twirling upside down in this case. Instead of bringing the endpoint closer to the $y$ axis, the extra work done by the motor is used to increase the frequency of the crankshafting motion at a higher rate, so that it transitions to about 10\% of the motor frequency. The mean $\chi/\omega$ of about 0.05 is similar to what Wada and Netz found for a thermally fluctuating fiber \cite[Fig.~1(c)]{wada2006non}. Our frequency results are also similar to those presented previously in \cite[Fig.~6(a)]{lee2014nonlinear}, where nonlinear behavior is observed for $\omega/\omega_c \gtrsim 1.6$.

Interestingly, while we find that the overwhirling state can be stable for frequencies as large as $3 \omega_c$, we also find that the ability to reach it depends on the initial fiber configuration. For instance, when $\omega/\omega_c$ is as small as 2.1, the fiber intersects itself on its way to the overwhirling state. An animation of this intersection, which is more severe for larger $\omega/\omega_c$, is shown when $\omega/\omega_c=2.8$ in one of the supplementary movies. This suggests a mechanism for the formation of plectonemes that were observed at high $\omega$ in \cite{bruss2019twirling}, and that the full nonlinear dynamics are likely more complex than previously thought. }

\section{Conclusions}
In this paper, we studied a single, inextensible filament that resists bending and twisting while immersed in a viscous fluid. Our focus was on adding twist elasticity to the model that we used in previous work \cite{maxian2021integral} for inextensible filaments with bending resistance. We did this in two parts: first, we compared two possible models for twist elasticity in the presence of a fluid, which we called the Kirchhoff and Euler models. While the two models have the same elastic constituitive laws, in the Kirchhoff model the angular velocity of the material frame is constrained to be equal to the angular fluid velocity, while in the Euler model the perpendicular components of the angular fluid velocity are discarded. The contribution of the first part of this paper was to reformulate the Kirchhoff model in a form similar to the Euler model, and use this to show that these two models produce dynamics (parallel and perpendicular angular velocities) that differ by $\mathcal{O}(\epsc^2)$ in the fiber interior, where $\epsc$ is the fiber aspect ratio. We chose to utilize the Euler model in the remainder of the paper, since the extra degree of rotational freedom (parallel torque and parallel rotational velocity) can be easily incorporated into the formulation we introduced in \cite{maxian2021integral}. To compute the additional force and parallel torque due to twist, we used the Bishop (twist-free) frame and evolved the angle of twist relative to that frame \cite{bergou2008discrete}. 

We defined the hydrodynamic mobility operators for translation-translation, translation-rotation, rotation-translation, and rotation-rotation coupling in terms of integrals of the corresponding Rotne-Prager-Yamakawa (RPY) kernels, then evaluated these integrals both asymptotically and numerically. Our preference for trans-trans, trans-rot, and rot-trans was to use the slender body \emph{quadrature} schemes described in Appendix\ \ref{sec:quads} to evaluate the integrals; we showed that this removes the spurious negative eigenvalues that result when the trans-trans integral is evaluated asymptotically (as in slender body \emph{theory}). The only remaining uncontrolled approximation in our mobility evaluation is to use local drag, and not the full RPY integral, for the rot-rot mobility. We did this because the $\mathcal{O}(\epsc^2)$ accuracy of this approximation was sufficient for our purposes. If necessary, the special quadrature scheme we developed for the doublet (in the context of trans-trans mobility) could also be applied to the rot-rot kernel to account for the nonlocal rot-rot hydrodynamics to higher accuracy. Using the RPY kernels to account for the fiber thickness undoubtedly gives some error relative to the true three dimensional geometry. Nevertheless, the equations we wrote are intended as a mathematical \emph{model} of a thin elastic rod in fluid, in much the same spirit as the Kirchhoff rod model itself, which is not derived by solving any three-dimensional elasticity problem but which aims to capture the essence of a thin rod's elastic behavior in a one-dimensional self-consistent continuum model. In the same way, our formulation of the fluid-structure interaction equations for a thin rod in fluid does not exactly correspond to the solution of any 3D Stokes problem with an internal moving tubular boundary \cite{koens2022tubular}, but instead we work with the variables of the Kirchhoff rod model and couple them to the fluid in a physically reasonable way. 

We applied our numerical method to two common single-filament problems in which twist elasticity has been traditionally neglected and included, respectively. We considered a relaxing bent filament and showed that twist elasticity has a negligible $\mathcal{O}(\epsc^2)$ effect, which justifies the assumption in many studies \cite{ehssan17, maxian2021integral} to neglect it. By contrast, the spinning filament problem requires finite twist elasticity to drive the filament into an unstable trajectory for spinning frequencies $\omega > \omega_c$. We showed that our numerical method reproduces exactly the semi-numerical prediction of \cite{wolgemuth2000twirling} for $\omega_c$ when using ellipsoidally-tapered filaments with local drag, and that the shape of the filament and the addition of nonlocal trans-trans hydrodynamics modify the critical frequency by at most 5\% \cite{wada2006non}. We showed that rot-trans coupling is more significant, lowering the frequency by another 15\% to about 80\% of the value predicted for ellipsoidally-tapered local drag, which explains about half of the deviation observed experimentally \cite{bruss2019twirling}. 

\change{Above the critical frequency $\omega_c$, the fiber transitions to a steady crankshafting rotation in the shape of a shepherd's crook (overwhirling) \cite{lim2004simulations}. While we found this state to be stable up to the highest frequencies we investigated, getting to that state from an initially straight filament requires complex dynamics that can lead to self-intersections, suggesting that in practice plectonemes may be observed instead \cite{bruss2019twirling}.} \rev{Our numerical method cannot handle self intersections, both because of the quadrature scheme and because we do not include lubrication, frictional, or steric forces \cite{brady1988stokesian}. }

Despite the length and complexity of this paper, we only considered a single filament! Yet it is quite interesting how many modeling and numerical choices we had to make for this seemingly simple problem. For instance, we chose \emph{not} to regularize the shape of the filament, keeping the radius constant all the way up to the endpoints. Because we use integrals of the RPY tensor, the mobility for constant regularized radius $\eps$ \emph{is} well-defined at the fiber endpoints, although its lack of smoothness causes slow convergence in our spectral method. This is in contrast to SBT, in which the local drag coefficient is not well-defined for cylindrical fibers, and some kind of ellipsoidally-tapered radius function has to be assumed \cite[Sec.~2.1]{maxian2021integral}. Future work could develop a hybrid approach that uses the RPY kernel for \emph{unequally-sized spheres} \cite{RPY_Shear_Wall} to simulate fibers with radius functions that decay smoothly to zero at the endpoints. This should both smooth out the mobility and keep the eigenvalues of the trans-trans mobility away from zero. The integrals involved are more complex than those for equally sized spheres because the integrand includes the radius function $\eps(s)$. Still, because the singularities (Stokeslet, doublet, rotlet, etc.) in the polydisperse integral kernels are unchanged, it should be possible to account for non-constant $\eps(s)$ using the slender body quadrature schemes developed here.

Another choice we made for the single filament was to stick with a first kind integral equation for the force-velocity relationship, i.e., the trans-trans mobility. While this yielded sufficient accuracy for our purposes, the eigenvalues of the first kind mobility cluster around zero, leading to ill-defined solutions near the fiber endpoints. Of course, one argument around this (which we made in \cite{maxian2021integral}) is that both the ill-conditioning and negative eigenvalues of the SBT mobility occur when $N \gtrsim 1/\epsRS$, which is too large for slender body approaches to be efficient. Still, a better-conditioned formulation (for trans-trans mobility) would be to use a second-kind integral equation of the form given in \cite[Eq.~(4)]{andersson2021integral}. To enable rapid dynamic simulation of the model in \cite{andersson2021integral}, it is necessary to develop slender body quadratures to evaluate integrals of regularized Stokeslets and doublets along the centerline, as we did for the first kind RPY mobility in Appendix\ \ref{sec:quads}. 

An important future extension of this work is to consider the dynamics of multiple filaments in a suspension with twist, as we did without twist in \cite{maxian2021integral, maxian2021simulations}. The mobility of Section\ \ref{sec:two} directly generalizes to multiple filaments, since the same RPY kernel can be used between distinct filaments and a filament and itself. The challenge lies in efficiently numerically computing the solution. Focusing on the static case first, we need to develop an efficient numerical method to sum the rot-trans and trans-rot hydrodynamic interactions over many bodies, since our previous work \cite[Sec.~4]{maxian2021integral} already treated the trans-trans kernel with a combination of the spectral Ewald method \cite{PSRPY} (for well-separated fibers)\footnote{We chose an Ewald approach because we are interested in triply-periodic systems; in free space a fast multipole method \cite{gimbutas2015computational} is likely the best choice, and systems with mixed periodicity are an avenue of active research \cite{lindbo2011fast, yan2018universal}.} and corrections for nearly-touching fibers \cite{barLud}. Extending our many-body-summation approach to twist is straightforward, but arduous. In particular, spectral Ewald methods have been developed for the rot-trans and trans-rot singularities of Stokes flow \cite{fiore2018rapid, af2017fast}, and the nearly-singular quadratures of \cite{barLud} readily generalize to the other singularities. For rot-rot coupling, the fast decay of the doublet singularity makes it a good approximation to only include the flows generated by nearly contacting fibers.

Then there is the question of computing the dynamics of many twisting filaments. In our previous SBT-based work \cite{maxian2021integral}, we treated terms involving nonlocal hydrodynamic interactions explicitly, and only treated the local drag terms implicitly. This allowed us to solve essentially diagonal linear systems, and to evaluate long-ranged hydrodynamics only once or at most a few times per time step. In the slender body quadrature approach followed here, there is not, at first sight, an explicit splitting into a ``local drag'' contribution and a nonlocal contribution, as there is in traditional SBT. But in the trans-trans mobility, we used singularity subtraction to factor out the leading order singularity in the RPY integral. This, combined with the integrals for $R < 2\eps$, give a ``local drag'' contribution that can be treated implicitly. The rot-rot mobility is even more dominated by local behavior, and so we can treat the angular velocity coming from torque on any other fibers explicitly and only treat the self contribution implicitly. A challenge for future exploration is how to efficiently handle pairs of filaments \rev{(or pieces of a single filament)} in near contact in a semi-implicit temporal integrator, without expensive GMRES iterations \cite[Secs.~4.5.2, 5.2]{maxian2021integral} at each time step. 


The spectral numerical method developed in this work is efficient, but lacks the simplicity, flexibility, and robustness of second-order ``blob-based'' methods that use $\mathcal{O}(\epsRS^{-1})$ discrete regularized singularities along the centerline. Spectral methods excel especially when blob-based discretizations are expensive (because of small aspect ratios) or unnecessary (because of smooth fiber shapes). They are more challenging, however, when the fiber shapes are nonsmooth, as it was only possible to model overwhirling with a great effort in anti-aliasing. A related even greater challenge for the future is to account for thermal fluctuations of elastic filaments in spectral approaches. 

\section*{Acknowledgments}
We thank Bjorn Dunkel and Vishal Prakash Patil for helpful discussions involving the Bishop frame. This work was supported by the National Science
Foundation through Research Training
Group in Modeling and Simulation under award RTG/DMS-1646339 and through the Division of Mathematical Sciences award DMS-2052515.

\begin{appendices}
\setcounter{equation}{0}
\renewcommand{\theequation}{\thesection.\arabic{equation}}
\section{RPY mobility \label{sec:rpysing}}
\delete{This appendix is devoted to the Rotne-Prager-Yamakawa (RPY) mobility that we use throughout this paper. We give the formulas for the $3 \times 3$ tensors $\Mbtt$, $\Mbtr$, $\Mbrt$, and $\Mbrr$ that define the mobility between two spheres \cite{wajnryb2013generalization}; these formulas form the integral kernels that we use the the fiber mobilities in\ \eqref{eq:UIBdef} and\ \eqref{eq:PsiMobs}. We then perform asymptotic evaluation of the mobilities $\Mctt$, $\Mctr$, $\Mcrt$, and $\Mcrr$ in Appendix\ \ref{sec:EPs}, which includes asymptotic expressions for the mobility near the fiber endpoints. These formulas yield SBT-type expressions which break the velocity into a local drag part (velocity at $s$ from force $\mathcal{O}(\eps)$ away from point $s$) and a finite part integral for the remaining (nonlocal) velocity. }

When the regularized delta function $\delta_\eps$ is a surface Delta function on a sphere of radius $\eps$ (i.e., $\delta_\eps(\V r)=1/(4\pi\eps^2)\delta_\eps (r-\eps)$), the kernels $\Mbtt, \Mbtr, \Mbrt$, and $\Mbrr$ in\ \eqref{eq:Ukernels} and\ \eqref{eq:Psikernels} can be written analytically, and are collectively known as the Rotne-Prager-Yamakawa (RPY) tensor. For the translational velocity\ \eqref{eq:Ukernels} at position $\V{x}$, they take the form \cite{wajnryb2013generalization}
\begin{align}
\label{eq:MbttRPY}
\left(8\pi \mu\right) \, \Mbtt\left(\V{x},\V{y}\right) \V{F}= \begin{cases}
 \left(\dfrac{\M{I}+\Rhat \Rhat}{R}+\dfrac{2\eps^2}{3} \dfrac{\M{I}-3\Rhat \Rhat}{R^3}\right)\V{F} & R > 2\eps \\[8 pt]
\left(\left(\dfrac{4}{3\eps}-\dfrac{3R}{8\eps^2}\right)\M{I}+\dfrac{R}{8\eps^2} \Rhat \Rhat \right)\V{F} & R \leq 2\eps
\end{cases}
\end{align}
\begin{equation}
\label{eq:MbtrRPY}
\left(8\pi \mu\right)\, \Mbtr\left(\V{x},\V{y}\right)\V{N} =\begin{cases}  \dfrac{\V{N}\times \Rhat}{R^2} & R > 2\eps \\[8 pt]
\dfrac{1}{2\eps^2}\left(\dfrac{R}{\eps}-\dfrac{3R^2}{8\eps^2}\right)\left(\V{N}\times \Rhat \right) & R \leq 2\eps 
\end{cases} 
\end{equation}
where $\V{F}$ and $\V{N}$ are the strength of a regularized force and torque, respectively, applied at position $\V{y}$, the displacement vector $\V{R}=\V{x}-\V{y}$, $R=\norm{\V{R}}$, $\Rhat=\V{R}/\norm{\V{R}}$, and $\Rhat \Rhat$ denotes the outer product of $\Rhat$ with itself. Note that the RPY tensor changes form for the case of overlapping spheres ($R \leq 2\eps$). At $R=2\eps$, the mobilities are once continuously differentiable. 

Likewise, for the angular velocity, the double convolutions in\ \eqref{eq:Psikernels} can be computed analytically, which yields the formulas \cite{wajnryb2013generalization}
\begin{gather}
\label{eq:MbrtRPY}
\left(8\pi \mu\right) \, \Mbrt\left(\V{x},\V{y}\right)\V{F} = \begin{cases}  \dfrac{\V{F} \times \Rhat}{R^2} & R > 2\eps \\[8 pt]
\dfrac{1}{2\eps^2}\left(\dfrac{R}{\eps}-\dfrac{3R^2}{8\eps^2}\right)\left(\V{F} \times \Rhat \right) & R \leq 2\eps 
\end{cases} \\[4 pt]
\label{eq:MbrrRPY}
\left(8\pi \mu\right) \, \Mbrr\left(\V{x},\V{y}\right) \V{N}= \begin{cases} -\dfrac{1}{2} \left(\dfrac{\M{I}-3\Rhat \Rhat }{R^3} \right) \V{N} & R > 2\eps\\[2 pt] 
\dfrac{1}{\eps^3} \left(\left(1-\dfrac{27R}{32\eps}+\dfrac{5R^3}{64 \eps^3}\right)\M{I}+\left(\dfrac{9 R}{32\eps}-\dfrac{3R^3}{64\eps^3}\right)\Rhat \Rhat\right)\V{N} & R \leq 2\eps.
\end{cases}
\end{gather}

\setcounter{equation}{0}
\section{Derivation of twist evolution equation \label{sec:twderiv}}
In this appendix, we derive the evolution equation\ \eqref{eq:twODE} for the twist $\psi$. We present this derivation because it uses straightforward calculus, and not the variational approach in the review article \cite{powers2010dynamics}.

Because the material frame $\left(\DO(s),\DT(s),\DR(s)\right)$ is orthonormal, it evolves by the parallel transport equation
\begin{equation}
\label{eq:bishEv}
\ds{\ind{\V{D}}{i}} = \V{\omega} \times \ind{\V{D}}{i},
\end{equation}
which also gives 
\begin{align}
\nonumber 
\psi &= \ds{\DO} \cdot \DT\\ 
\label{eq:thetaXs}
& =\left(\V{\omega} \times \DO \right) \cdot \DT= \left(\DO \times \DT\right) \cdot \V{\omega} = \Xs \cdot \V{\omega}
\end{align}
So $\omega^\parallel = \psi$. Likewise, setting $i=3$ in\ \eqref{eq:bishEv} and using the constraint $\Xs=\DR$ from\ \eqref{eq:DrXs} gives $\ds{\Xs} = \V{\omega} \times \Xs$. Crossing both sides with $\Xs$, we obtain the total representation of $\V{\omega}$ as \cite{bergou2008discrete}
\begin{equation}
\label{eq:omrep}
\V{\omega}= \Xs \times \ds{\Xs} +\Xs \psi
\end{equation}

We now equate the $t$ derivative of\ \eqref{eq:bishEv} with the $s$ derivative of\ \eqref{eq:Devolve}, 
\begin{gather}
\nonumber
\dt{\left(\V{\omega} \times \ind{\V{D}}{i}\right)} = \ds{\left(\V{\Omega} \times \ind{\V{D}}{i}\right)}\\ \nonumber
\dt{\V{\omega}} \times \ind{\V{D}}{i} + \V{\omega} \times \dt{\ind{\V{D}}{i}}=\ds{\V{\Omega}} \times \ind{\V{D}}{i}+\V{\Omega} \times \ds{\ind{\V{D}}{i} } \\ \label{eq:beforeJac}
\dt{\V{\omega}} \times \ind{\V{D}}{i} = \ds{\V{\Omega}}\times \ind{\V{D}}{i}+\V{\Omega}\times \left(\V{\omega} \times \ind{\V{D}}{i} \right)-\V{\omega} \times \left(\V{\Omega} \times \ind{\V{D}}{i}\right).
\end{gather}
The last equality used the evolution equations\ \eqref{eq:bishEv} and\ \eqref{eq:Devolve} for $\ind{\V{D}}{i}$. Using the triple cross product identity $\V{a} \times \left(\V b \times \V c\right)+\V{c} \times \left(\V a \times \V b\right)+\V{b} \times \left(\V c \times \V a\right)=\V 0$, we can write\ \eqref{eq:beforeJac} as 
\begin{gather}
\dt{\V{\omega}} \times \ind{\V{D}}{i} = \ds{\V{\Omega}}\times \ind{\V{D}}{i}+\left(\V{\Omega} \times \V{\omega}\right)\times \ind{\V{D}}{i},
\end{gather}
which must hold for all $i$. This gives \cite[Eq.~(54)]{powers2010dynamics}
\begin{equation}
\label{eq:domdt}
\dt{\V{\omega}} = \ds{\V{\Omega}} + \left(\V{\Omega} \times \V{\omega}\right).
\end{equation}
To derive the evolution equation for the twist angle, we take the $\Xs$ component of\ \eqref{eq:domdt} and use\ \eqref{eq:Devolve}
\begin{gather}
\dt{\left(\V{\omega} \cdot \Xs\right)} - \V{\omega} \cdot \dt{\Xs} = \ds{\left(\V{\Omega} \cdot \Xs\right)} - \V{\Omega} \cdot \ds{\Xs} + \left(\V{\Omega} \times \V{\omega}\right) \cdot \Xs \\ \nonumber
\dt{\left(\V{\omega} \cdot \Xs\right)} = \V{\omega} \cdot \left(\V{\Omega} \times \Xs\right)+ \ds{\left(\V{\Omega} \cdot \Xs\right)} - \V{\Omega} \cdot \ds{\Xs} + \left(\V{\Omega} \times \V{\omega}\right) \cdot \Xs.
\end{gather}
The triple products cancel, and using\ \eqref{eq:thetaXs}, we get the twist evolution equation\ \eqref{eq:twODE}, 
\begin{align*}
\dt{\psi}= \ds{\Omega^\parallel}-\left(\V{\Omega} \cdot \ds{\Xs}\right).
\end{align*}

\setcounter{equation}{0}
\section{Euler model for intrinsic curvature and twist \label{sec:Intrin}}
In Sections\ \ref{sec:krod} and\ \ref{sec:euler}, we saw that the Euler model can be derived from the Kirchhoff rod model by first setting $\V{F}=\V{F}^{(\kappa)}+\V{F}^{(\twmod)}+\V{\Lambda}$ in\ \eqref{eq:Fdecomp}, where $\V{F}^{(\kappa)}$ and $\V{F}^{(\twmod)}$ are defined to eliminate the perpendicular torque in $\ds{\V{N}}$, then setting $\V{\Lambda}=T\Xs$, which eliminates the perpendicular torque Lagrange multipliers that are the distinguishing feature of the Kirchhoff rod model. These two steps eliminate \emph{all} perpendicular torques from the formulation, which leads to the Euler model in Section\ \ref{sec:euler}. 

In this appendix, we repeat this procedure on a fiber with intrinsic curvature and twist (e.g., a flagellum) to derive the corresponding Euler model. In this case, we have the moment \cite{olson2013modeling}
\begin{align}
\nonumber
\V{N} &= \kappa \left(\left(\ds{\DT} \cdot \DR - \kappa_1\right) \DO + \left(\ds{\DR} \cdot \DO - \kappa_2 \right) \DT \right) + \twmod \left(\ds{\DO} \cdot \DT - \phi\right) \DR \\ 
\label{eq:Nintrin}
& = \kappa \left(\Xs \times \ds{\Xs} - \kappa _1 \DO - \kappa_2 \DT\right) + \twmod \left(\psi- \phi\right)\Xs,
\end{align}
where $\kappa_1$ and $\kappa_2$ are the preferred curvatures, and $\phi$ is the preferred twist angle.

We illustrate the key steps here assuming constant $\kappa_1$, $\kappa_2$, and $\phi$. Taking the $s$ derivative of $\V{N}$ and separating the result into tangential and perpendicular parts, we have 
\begin{align}
\ds{\V{N}} & =  \kappa \left(\Xs \times \ds^2{\Xs} - \kappa _1 \ds{\DO} - \kappa \ds{\DT}\right) + \twmod \left(\left(\ds \psi\right) \Xs+(\psi -\phi)\ds{\Xs}\right)\\
& = \Xs \times \left(\kappa \left( \ds^2{\Xs}+\kappa_1 \left(\Xs \times \ds{\DO}\right) + \kappa_2 \left(\Xs \times \ds{\DT}\right)\right) - \gamma \left(\psi -\phi\right)\left(\Xs \times \ds{\Xs}\right)\right) \\ \nonumber
&+ \left(\kappa \left(-\kappa_1 \left(\ds{\DO}\cdot \Xs\right) - \kappa_2 \left(\ds{\DT} \cdot \Xs\right)\right) + \twmod \ds{\psi}\right)\Xs
\end{align}
To derive an Euler model, we now amend\ \eqref{eq:Ftwkap} and\ \eqref{eq:Fdecomp} to 
\begin{gather}
\V{F}^{(\kappa)} = - \kappa \left(\ds^2 \Xs - \kappa_1 \psi \DO - \kappa_2 \psi \DT\right),  \qquad \V{F}^{(\twmod)} =  \twmod \left(\psi-\phi\right) \left(\Xs \times \ds{\Xs}\right), \\ \nonumber \V{F}=\V{F}^{(\kappa)}+ \V{F}^{(\twmod)} +T\Xs.
\end{gather}
Just like in the case of intrinsically straight and untwisted fibers, substituting this representation into\ \eqref{eq:fforce} gives the force and torque density applied to the fluid as
\begin{gather}
\V{f} = \ds{\V{F}}, \qquad 
\V{n}=n^\parallel \Xs, \qquad n^\parallel = \kappa \left(\kappa_1 \left(\DO \cdot \ds{\Xs}\right) + \kappa_2 \left(\DT \cdot \ds{\Xs}\right)\right) + \gamma \ds{\psi}.
\end{gather} 

Free fiber boundary conditions are modified in the case of intrinsic curvature and twist. As before, we require $\V{N}=\V{0}$ and $\V{F}=\V{0}$ at the free end, which implies that
\begin{gather}
\psi= \phi, \qquad 
\ds{\Xs}=  \kappa_2 \DO-\kappa_1 \DT , \qquad
\ds^2{\Xs}= \kappa_1 \phi \DO + \kappa_2 \phi \DT, \quad \text{and} \quad T=0
\end{gather}
at a free end. The equation for $\ds{\Xs}$ is derived by using\ \eqref{eq:Nintrin} and subsituting $\DO = \DT \times \Xs$ and $\DT = \Xs \times \DO$.

\subsection{Solving the Bishop ODE \label{sec:BishNumer}}
In the case of intrinsic curvature and twist, the material frame is involved in the forces and torques, and so it must be computed from the Bishop frame and $\theta$ using\ \eqref{eq:D1}. Computing the Bishop frame requires the solution of the Bishop ODE\ \eqref{eq:BishODE}, which is subject to the condition\ \eqref{eq:BishBC},
\begin{gather}
\label{eq:BishODEBC}
\ds{\V{b}^{(1)}} = \left(\Xs \times \ds{\Xs}\right) \times {\V{b}^{(1)}} \\ \nonumber
\V{b}^{(1)}\left(L/2\right) = \DO\left(L/2\right).
\end{gather}

\begin{figure}
\centering
\includegraphics[width=0.5\textwidth]{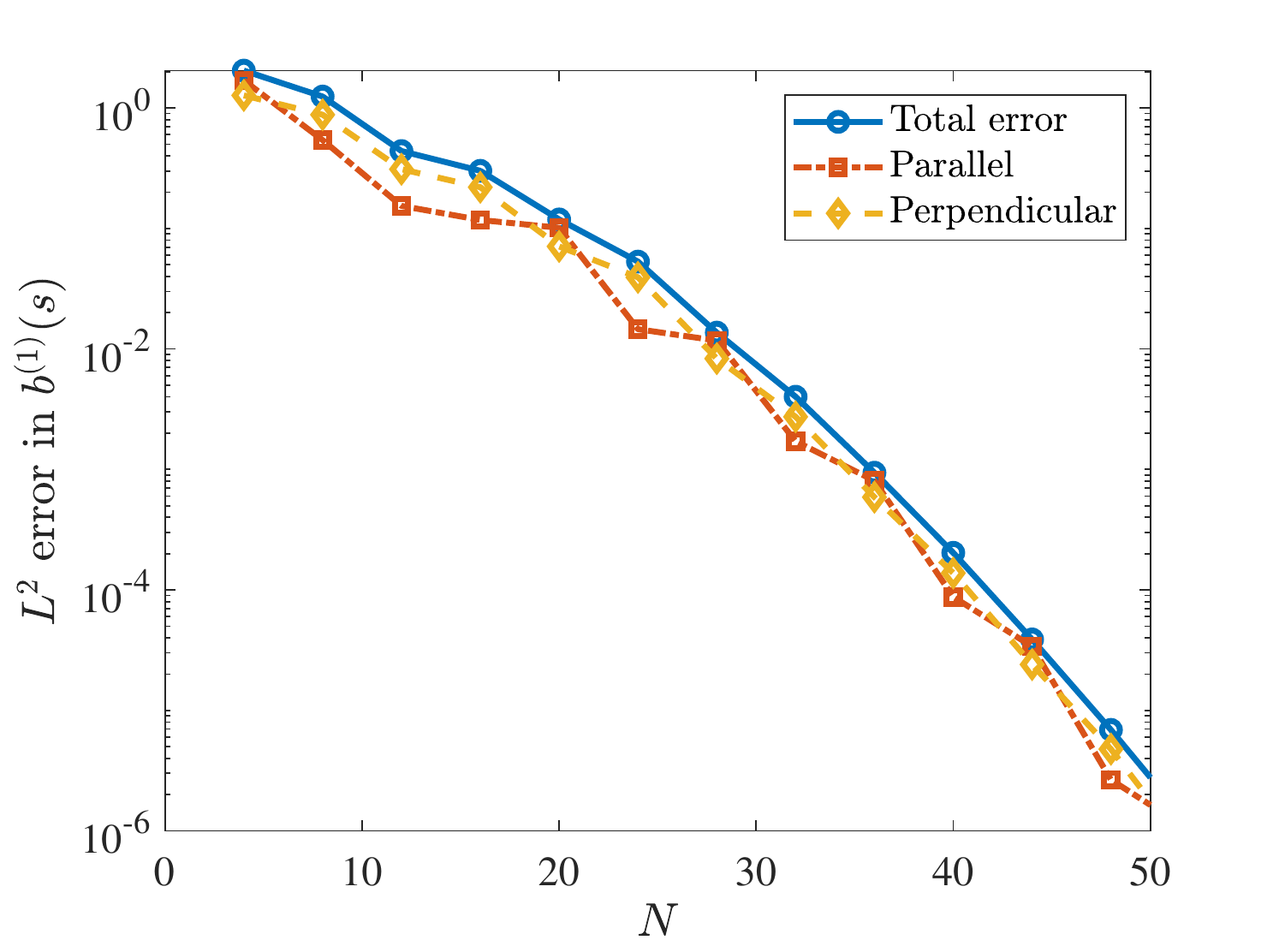}
\caption{\label{fig:BishErs}The $L^2$ function errors in the Bishop vector $\V{b}^{(1)}(s)$. We use the fiber\ \eqref{eq:Xsfibq} with $q=7$ and generate a reference solution using $100$ collocation points (these vectors have unit length to 13 digits). We then solve\ \eqref{eq:DiscBish} using $N$ points and measure the $L^2$ function error in $\V{b}^{(1)}(s)$ on a common grid. We separate the total error (blue circles) into a component parallel to $\V{b}^{(1)}(s)$ (red squares) and a component perpendicular to $\V{b}^{(1)}(s)$ (yellow diamonds).}
\end{figure}

To discretize the Bishop ODE\ \eqref{eq:BishODEBC}, we first integrate both sides of\ \eqref{eq:BishODEBC} to obtain an integral equation for $\V{b}^{(1)}$. Since each fiber has at most $40$ points on its centerline, we will use dense linear algebra to solve for $\V{b}^{(1)}$ at all of the collocation points. The discrete integral form of\ \eqref{eq:BishODEBC} is
\begin{equation}
\V{b}^{(1)}+ \V{c} = \M{D}^\dagger  \left(\Xs \times \ds{\Xs}\right) \times\V{b}^{(1)}.
\end{equation}
where $\V{c}$ is an unknown integration constant and $\M{D}^\dagger$ is the pseudo-inverse of the Chebyshev differentiation matrix. Combining with the boundary condition in\ \eqref{eq:BishODEBC}, and denoting with $\M{I}(s)$ evaluation at $s$, we get the system
\begin{gather}
\label{eq:DiscBish}
\begin{pmatrix} 
\M{I}- \M{D}^\dagger \left( \left(\Xs \times \ds{\Xs}\right) \times \right) & \M{I} \\
\M{I}(L/2) & \V{0}
\end{pmatrix}
\begin{pmatrix} \V{b}^{(1)} \\ \V c \end{pmatrix} = \begin{pmatrix} \V 0\\ \DO(L/2) \end{pmatrix},
\end{gather}
which we solve for $\V{b}^{(1)}$. 

In Fig.\ \ref{fig:BishErs}, we plot the errors in $\V{b}^{(1)}(s)$ as a function of the grid size for the fiber\ \eqref{eq:Xsfibq} with $q=7$ that satisfies the free fiber BCs (the refined solution is generated with $N=100$, in which all vectors have unit length to tolerance $10^{-13}$). We see that the errors in the parallel and perpendicular direction are the same, which means that choosing to renormalize $\V{b}^{(1)}$ will not increase the overall error (in Fig.\ \ref{fig:BishErs} we do not renormalize $\V{b}^{(1)}$).

\setcounter{equation}{0}
\section{Principle of virtual work \label{sec:ftilde}}
In this appendix, we show that the Lagrange multiplier forces $\V{\Lambda}(s)$ in both the Kirchhoff and Euler models do no work on the fluid. When we substitute\ \eqref{eq:Fdecomp} and\ \eqref{eq:ndecomp} into the corresponding Stokes equation for the pointwise fluid velocity $\V{u}(\V{x})$, we obtain 
\begin{gather}
\nonumber
\grad\pi(\V{x})=\mu \grad^{2}\V u(\V{x})+ \int_0^L \left(\ds{\left(\V{F}^{(\kappa)}(s)+\V{F}^{(\twmod)}\right)}(s)\right)\delta_{\eps}\left(\V{x}-\V X(s)\right) \, ds\\ \ + \frac{\V{\nabla}}{2} \times \int_0^L n^\parallel(s) \Xs(s) \delta_{\eps}\left(\V{x}-\V X(s)\right) \, ds + \left(\Lop{A}\V{\Lambda}\right)\left(\V{x}\right),
\end{gather}
where we have defined a linear operator
\begin{gather}
\label{eq:AF}
\left(\Lop{A}\V{\Lambda}\right)(\V{x}) =  \int_0^L \left(\ds{\V{\Lambda}}(s)\right)\delta_{\eps}\left(\V{x}-\V X(s)\right) \, ds+ \frac{\V{\nabla}}{2} \times \int_0^L \left(\Xs(s) \times \V{\Lambda}(s)\right) \delta_{\eps}\left(\V{x}-\V X(s)\right) \, ds.
\end{gather}
Starting with the Kirchhoff model, let us consider the constraint\ \eqref{eq:simplest_cons}, which can be written as 
\begin{equation}
\label{eq:ConstrD}
\ds{\V{U}}+\Xs \times \V{\Psi}=\V{0}.
\end{equation} 
Recalling that $\V{U}$ and $\V{\Psi}$ are also convolutions with the regularized function as defined in\ \eqref{eq:Uavg} and\ \eqref{eq:OmIBdef}, we write this constraint as 
\begin{gather}
\label{eq:AstarU}
\ds \left(\int \V{u}(\V{x}) \delta_\eps\left(\V{X}(s)-\V{x}\right) d\V{x}\right) + \Xs(s) \times \int \left(\frac{\V{\nabla}}{2} \times \V{u}(\V{x})\right) \delta_\eps\left(\V{X}(s)-\V{x}\right) d\V{x} =0.
\end{gather}
It is not difficult to show that\ \eqref{eq:AstarU} can be written as $\Lop{A}^* \V{u}=\V{0}$, where $*$ denotes an $L^2$ adjoint. It follows that 
\begin{gather}
\label{eq:D3}
\left\langle \Lop{A} \V{\Lambda}, \V{u}\right\rangle = \left\langle \V{\Lambda}, \Lop{A}^*\V{u}\right\rangle=\V{0},
\end{gather}
meaning the force $\V{\Lambda}$ does no work on the fluid when the velocity $\V{u}$ is constrained to satisfy\ \eqref{eq:simplest_cons}. 

The same principle continues to hold in the Euler model as well. Specifically, since in the Euler model\ \eqref{eq:LamEul} gives $\V{\Lambda}=T\Xs$, the forcing term $\Lop{A}\V{\Lambda}=\Lop{B}T$, where 
\begin{equation}
\Lop{B} T = \int_0^L \ds{\left(T(s)\Xs(s)\right)} \delta_\eps\left(\V{x}-\V{X}(s)\right) \, ds .
\end{equation}
It is straightforward to show that the principle of virtual work $\langle \Lop{B} T , \V{u} \rangle =0$ is (weakly) equivalent to the inextensibility constraint\ \eqref{eq:inex2},
\begin{equation}
\label{eq:D5}
\Lop{B}^* \V{u} = \Xs \cdot \ds{\left(\int \V{u}(\V{x}) \delta_\eps\left(\V{x}-\V{X}(s)\right) \, d\V{x}\right)}=\Xs \cdot \ds{\V{U}}=0.
\end{equation}
Note that both\ \eqref{eq:D3} and\ \eqref{eq:D5} assume the free fiber boundary conditions\ \eqref{eq:Tfree}. For clamped fibers, there may be work performed by the forcing at the endpoints (e.g., the motor spinning the fiber in Section\ \ref{sec:whirl}). \rev{The saddle-point structure of our formulation shows that the dynamical Euler equations we solve are a generalized gradient descent dynamics on the manifold of inextensible fibers.}

By combining\ \eqref{eq:AstarU} and\ \eqref{eq:D3} we can express the (vanishing) work performed by the constraint forces in terms of the fiber motion 
\begin{align}
\left\langle \Lop{A}\V{\Lambda}, \V{u} \right\rangle &= \int_0^L \V{\Lambda}(s) \cdot \left(\ds{\V{U}}(s)+\Xs(s) \times \V \Psi(s)\right) \, ds  \\
\label{eq:CombinedVW}
& = -\int_0^L \V{U}(s) \cdot \V{\lambda}(s)\, ds - \int_0^L \V{n}_\Lambda(s) \cdot \V{\Psi}(s) \, ds = 0,
\end{align}
where $\V{n}_\Lambda = \Xs \times \V{\Lambda}=\Xs \times \V{\Lambda}^\perp$ is the contribution to the torque density\ \eqref{eq:torqK} from $\V{\Lambda}$. \rev{In the Euler model, $\V{\Lambda}^\perp = \V{0}$ and\ \eqref{eq:CombinedVW} becomes\ \eqref{eq:kIP1}, which we used to enforce $\V{\Lambda}(s)$ be tangential. By contrast, for the Kirchhoff model\ \eqref{eq:CombinedVW} vanished for any $\V{\Psi}$ since\ \eqref{eq:ConstrD} holds, and therefore there is no additional constraint on $\V{\Lambda}$ from the principle of virtual work.}

\setcounter{equation}{0}
\section{Comparison of Kirchhoff and Euler models\label{sec:KENumer}}
\rev{In this appendix, we derive and confirm the asymptotic results\ \eqref{eq:dUEst} and\ \eqref{eq:dOmEst} that were given in Section\ \ref{sec:compareMeth}. After using the estimates\ \eqref{eq:LDmats} to estimate the difference between the two formulations, we develop a second order method for the Kirchhoff and Euler models which can robustly handle jumps in the forces at the endpoints. We then apply the method to a bent, twisted filament to show the equivalence of the two models (away from the endpoints) as $\epsc \rightarrow 0$.}

\subsection{Asymptotic estimates}
To estimate the strength of $\dF$, we use the scalings\ \eqref{eq:LDmats}, repeated here for clarity
\begin{gather}
\Lop{M}_\text{tt}\sim \log{\epsc}, \qquad \Lop{M}_\text{tr} \sim \log{\epsc}, \qquad
\Lop{M}_\text{rt}\sim \log{\epsc}, \qquad \Lop{M}_\text{rr} \sim \epsc^{-2},
\end{gather}
Beginning with the perpendicular Lagrange multipliers $\dF^\perp=\V{\Lambda}^\perp$, we first observe that the mismatch $\V{m}$ in\ \eqref{eq:mmT} is $\mathcal{O}(\log{\epsc})$ at most, since we don't expect $\Eul{\V{F}}$ as defined in\ \eqref{eq:F0} to change with $\epsc$. So $\dF^\perp$ cannot exceed $\mathcal{O}(\epsc^2 \log{\epsc})$, since otherwise the angular velocity from $\Mcrr$ would be too large to be consistent with the $\mathcal{O}(\log{\epsc})$ mismatch. If $\dF^\perp \sim \epsc^2 \log{\epsc}$, then the rotational velocity $\Delta \V{\Psi}$ is $\mathcal{O}(\log{\epsc})$, and therefore of the correct order to cancel the mismatch $\V{m}$, while the translational velocity $\Delta \V{U}=\Mctt \left(\ds{\dF}\right)+\Mctr \left(\Xs \times \dF\right)$, is of the order $\Delta \V{U}\sim \epsc^2 \left(\log{\epsc}\right)^2$ (if we assume that $\dF$ and its derivative are of the same order, which should be true away from the endpoints).\footnote{Recall that $\dF^\parallel$ is a Lagrange multiplier that makes the perturbation velocity $\Delta \V{U}$ inextensible. Since the interior velocity is at worst $\Delta \V{U} \sim \epsc^2 \left(\log{\epsc}\right)^2$ and the translational mobility is $\mathcal{O}(\log{\epsc})$, the Lagrange multiplier $\dF^\parallel \sim \epsc^2 \log{\epsc}$ at worst, just like $\dF^\perp$. }  This implies that in the fiber interior $\ds{\dU} \sim  \epsc^2 \left(\log{\epsc}\right)^2$, which gives\ \eqref{eq:dUEst}.

The difference in parallel angular velocity $\Delta \Omega^\parallel$ is given asymptotically by
\begin{align}
\Delta \Omega^\parallel =\Xs \cdot \dOm &=  \Xs \cdot \left(\Lop{M}_\text{rt}\left(\ds{\dF}\right) +\Lop{M}_\text{rr}\left(\Xs \times \dF\right)\right)  \approx  \Xs \cdot \left(\Lop{M}_\text{rt}\left(\ds{\dF}\right)\right),
\end{align}
where the last approximate equality is because $\Lop{M}_\text{rr}\left(\Xs \times \dF \right)$ is orthogonal to $\Xs$ to leading order (see Section\ \ref{sec:RPYasymp}). In fact, in the fiber interior $\Lop{M}_\text{rr}\left(\Xs \times \dF \right) \cdot \Xs \sim \log{\epsc} \norm{\dF}\sim \epsc^2 \left(\log{\epsc}\right)^2$. Adding the rot-trans contribution, since $\dF \sim \epsc^2 \log{\epsc}$ and $\Mcrt \sim \log{\epsc}$, we get $\Delta \Omega^\parallel \sim \epsc^2 \left(\log{\epsc}\right)^2$. To non-dimensionalize this, we need to estimate the parallel rotational velocity $\Eul{\Omega^\parallel}$  defined in\ \eqref{eq:OmE}. When there is no twist, there is no parallel torque, so $\Eul{\Omega^\parallel} = \Xs \cdot \left(\Lop{M}_\text{rt}\ds{\V{F}_E}\right) \sim \log{\epsc}$, and the relative difference should be $\norm{\Delta \Omega^\parallel} / \norm{\Omega_E^\parallel} \sim \epsc^2 \log{\epsc}$. When there is twist, the disparity is even larger, since $\Eul{\Omega^\parallel} \sim \epsc^{-2}$, which explains the two different cases in\ \eqref{eq:dOmEst}.

Summarizing, the mismatch $\V{m}$ in\ \eqref{eq:mprobT} is compensated for by $\dOm^\perp$, which grows as $\epsc \rightarrow 0$, even as $\dF$ decays. These are the components of the velocity that are discarded in the Euler model, so it does not matter if they differ from those of the Kirchhoff model. That said, because $\dF$ must be zero at the endpoints by the boundary conditions\ \eqref{eq:FBC}, the  perpendicular torques cannot generate angular velocity there, and the term $\dOm^\perp$ cannot compensate for the mismatch $\V{m}$ at the endpoints. Instead, the mismatch there has to be compensated for by terms involving $\ds{\dF}$, meaning that the derivative of $\dF$ must be large near the boundaries, as we see in Fig\ \ref{fig:PertEffects}. This makes a spectral method a poor choice for the Kirchhoff rod formulation, and is the reason why we develop a second-order method to solve the Kirchhoff equations next.

\subsection{Second-order numerical method \label{sec:SecOrder}}
\delete{This Appendix outlines the \emph{second-order} discretizations that we use for the Kirchhoff and Euler models. These discretizations are not the main focus of the paper, and are used to compare the Kirchhoff and Euler models in Section\ \ref{sec:compareMeth}, where we want to obtain an accurate solution for the nonsmooth Lagrange multipliers $\dF$ in the Kirchhoff model, and to generate robust reference solutions to verify the spectral method in Section\ \ref{sec:specEul}. }

Following the method outlined in \cite{keavRPY}, we introduce a number of links $N$ and set $\Delta s = L/N$. The discretization of the fiber is \emph{staggered}: the ``nodes'' are located at $s_j=j \Delta s$, where $j=0, \dots N$, so that there are actually $N+1$ nodes. The ``links'' are centered on $s_{j+1/2}=(j+1/2)\Delta s$, where $j=0, \dots, N-1$, so that there are $N$ links. 

\subsubsection{Euler model \label{sec:Euler2}}
In our second-order discretization of the Euler problem\ \eqref{eq:euler}, we define $\V{\lambda}$, $\V{f}^{(\kappa)}$, and $\V{f}^{(\gamma)}$  at link centers $(j+1/2)\Delta s$, for $j=0,\dots N-1$. All mobility operators in\ \eqref{eq:UIBdef} and\ \eqref{eq:PsiMobs} (including the rot-rot mobility) are then replaced with the obvious discretizations of the integrals which sum over the links with weight $\Delta s$. For example, 
\begin{equation}
\label{eq:loworderMob}
\left(\Mctt \V{f} \right)(s_{i+1/2}) \approx \sum_{j=0}^{N-1} \Mbtt\left(\V{X}_{i+1/2},\V{X}_{j+1/2}\right) \V{f}_{j+1/2}\Delta s :=\left(\Mtt\V{f}\right)_{i+1/2}.
\end{equation}
Note that this discretization does not correspond to a more common discretization of an inextensible bead-link polymer chain, which would place forces on $N+1$ blobs connected by $N$ inextensible links, e.g., \cite{nguyen2018impacts}. We use the staggered discretization instead for a more direct comparison to the Kirchhoff model.

To complete the second-order discretization of\ \eqref{eq:euler}, we have to define the matrices $\M{K}$ and $\M{K}^*$. In the second-order case, these matrices are exact transposes of each other, so we only have to define one of them. The matrix $\M{K}$ is a midpoint discretization of the continuum operator\ \eqref{eq:Kcdef}, with action given by
\begin{align}
\left(\M{K}\V{\alpha}\right)_{j+1/2} = \widebar{\V{U}}+ \sum_{k=0}^{j-1} &\Delta s \left(\alpha_{1,k+1/2} \V{n}_{1,k+1/2}  + \alpha_{2,k+1/2} \V{n}_{2,k+1/2}\right) \\ \nonumber +& \frac{\Delta s}{2} \left(\alpha_{1,j+1/2} \V{n}_{1,j+1/2} + \alpha_{2,j+1/2} \V{n}_{2,j+1/2}\right).
\end{align}
The midpoint quadrature assigns a weight of $\Delta s$ to all prior nodes $1/2, 3/2, \dots j-1/2$, and then a weight of $\Delta s/2$ to the current node $j+1/2$. The matrix $\M{K}^*$ is the transpose of $\M{K}$. 

Finally, the parallel angular velocity\ \eqref{eq:OmParAgain} is obtained at the link centers by dotting the full angular velocity, $\Lop{M}_\text{rt} \V{f}+  \Lop{M}_\text{rr}  \left(n^\parallel \Xs \right)$ (computed at the link centers using the forces computing at the link centers), with the tangent vectors $\Xs$, which are also defined at the link centers. 

\subsubsection{Kirchhoff model \label{sec:Kirch2}}
The discretization of the Kirchhoff model, which we take from \cite{keavRPY}, is more complicated because the forces $\V{\Lambda}$ are defined at the \emph{nodes}, and not the link centers. To communicate between the nodes and links, we let $\M{D}_\text{lb}$ be the $N \times (N+1)$ differentiation matrix, which gives derivates at the links from values at the nodes. We also define $\M{D}_\text{bl}$ as the $(N-1) \times N$ differentiation matrix, which gives derivates at interior nodes from values at the link centers. Likewise the $N \times (N+1)$ averaging matrix $\M{A}_\text{lb}$ gives averages at links from nodes, and the $(N-1) \times N$ averaging matrix $\M{A}_\text{bl}$ gives averages at interior nodes from links. This gives the $(N-1) \times (N+1)$ system of equations discretizing\ \eqref{eq:mprobT}
\begin{gather}
\label{eq:discmmp}
\left[\M{D}_\text{bl}\left(\M{M}_\text{tt}\M{D}_\text{lb}+\M{M}_\text{tr} \M{C} \M{A}_\text{lb} \right) + \M{A}_\text{bl}\M{C}\left(\M{M}_\text{rt}\M{D}_\text{lb}+\M{M}_\text{rr} \M{C}\M{A}_\text{lb}\right)\right] \dF = -\V{m}\\
\label{eq:discmismatch}
=-\left(\M{D}_\text{bl}\Mtt+\M{A}_\text{bl}\M{C}\Mrt\right)\left(\V{f}^{(\kappa)}+\V{f}^{(\gamma)}+\V{\lambda}_E\right)-\left(\M{D}_\text{bl}\Mtr+\M{A}_\text{bl}\M{C}\Mrr \right)\left(n^\parallel \Xs\right)
\end{gather}
which we solve for $\dF$ at the $N-1$ \emph{interior blobs} away from the fiber boundaries.  Because $\V{\Lambda}$ is defined at $N+1$ points, the other two equations are the BCs $\dF_j=\V{0}$ for  $j=0, N$, see\ \eqref{eq:FBC}. The matrix $\V{C}$ represents cross products with $\Xs_{j+1/2}=\Xs\left((j+1/2)\Delta s\right)$, where in the second-order method $\Xs=\ds{\V{X}}$ is evaluated analytically at the links and $\left(\M{C}\V{X}\right)_{j+1/2}=\Xs_{j+1/2} \times \V{X}_{j+1/2}$. Notice that the calculation of the mismatch\ \eqref{eq:discmismatch} requires calculation of the Euler quantities $\V{f}^{(\kappa)}$ and $\V{f}^{(\gamma)}$ and $\V{\lambda}_E$ at the $N$ links, which we first obtain by solving the Euler problem to second order accuracy using the method in Section\ \ref{sec:Euler2}. 

\subsection{Numerical results \label{sec:numerres}}
\delete{We now solve the mismatch problem\ \eqref{eq:mprobT} numerically to verify the estimates of Section\ \ref{sec:compareMeth}. To do this robustly, we use the second-order, blob-based numerical method described in Appendix\ \ref{sec:Kirch2} and in \cite{keavRPY}. The Euler model that we use \cite{maxian2021integral} is slightly different from the traditional tension-based form \cite{ts04}, as we discuss in Section 5. The second-order numerical method that we use for it is given in Appendix\ \ref{sec:Euler2}. }

To numerically compare the Euler and Kirchhoff constitutive models, we will apply the second-order methods to the fiber\ \eqref{eq:Xsfibq} with $q=1$ and $L=2$ (see Fig.\ \ref{fig:FibShapes} for a visual). This fiber is chosen to satisfy the free fiber boundary conditions\ \eqref{eq:EulBC}. We will also assume an instantaneous twist
\begin{equation}
\label{eq:theta}
\psi= \sin{(2\pi s)},
\end{equation}
so that all of the free fiber boundary conditions\ \eqref{eq:EulBC} are satisfied. Because $\V{X}$ satisfies the free-fiber BCs, we can compute the elastic force to spectral accuracy simply by differentiating $\Xs$ three times analytically, and then setting $\V{f}^{(\kappa)} = -\kappa \ds^3{\Xs}$. For this section, we will compute all geometric quantities, for example $\ds{\Xs}$, $\ds^2{\Xs}$ and $\ds{\psi}$, analytically. We will use a viscosity of $\mu=1$ and aspect ratios $\epsRS=\eps/L$ that are a function of the regularized singularity radius $\eps$; the corresponding value of the true fiber radius $\rc\approx \eps/1.12$ (see\ \eqref{eq:ahat}). 

\subsubsection{Bending but no twisting}
\begin{figure}
\hspace{-3em} 
\includegraphics[width=1.1\textwidth]{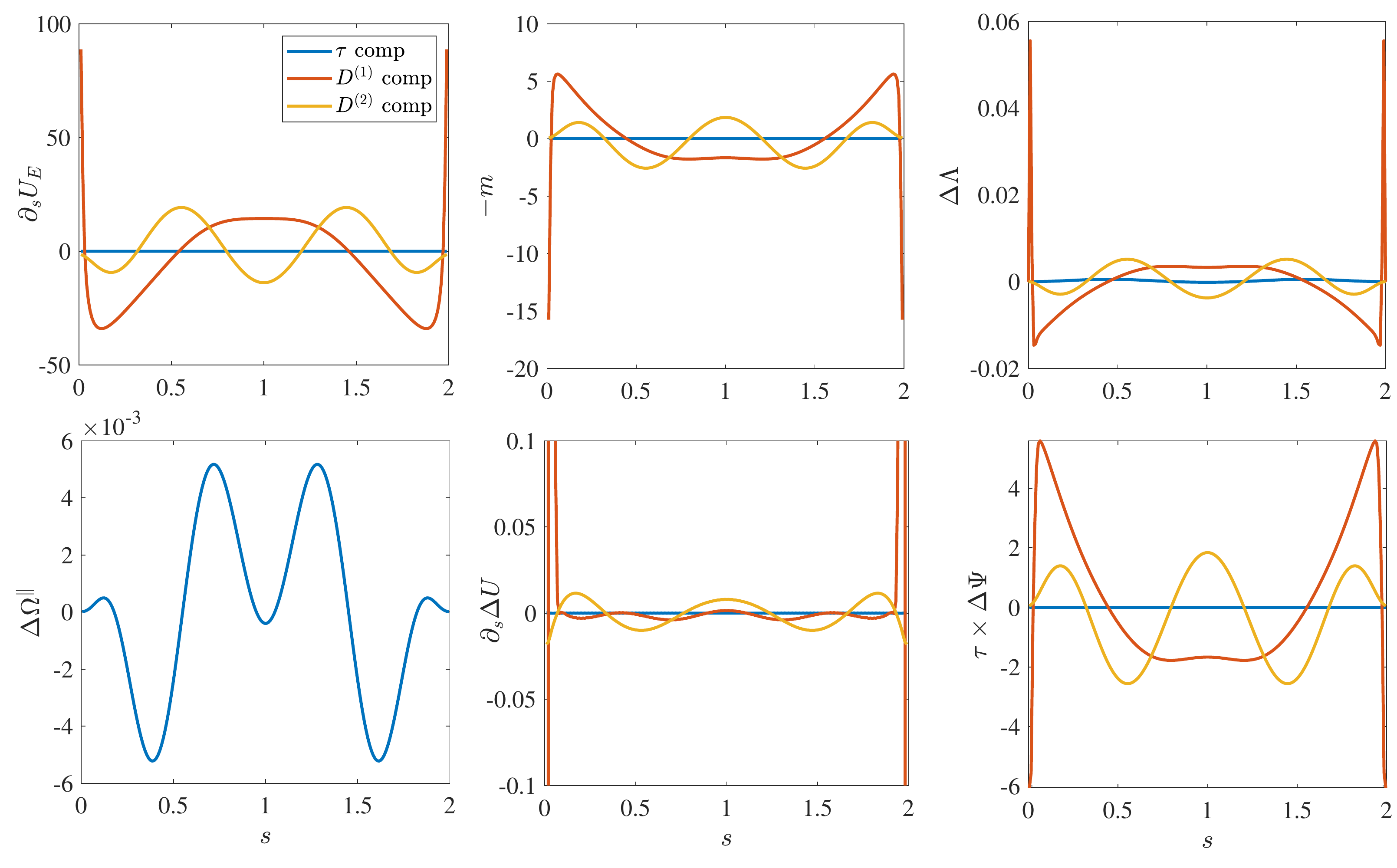}
\caption{\label{fig:PertEffects} Solving the mismatch problem\ \eqref{eq:mprobT} with bending stiffness $\kappa=1$, twist modulus $\gamma=0$, and slenderness $\epsRS=\eps/L =0.005$. All quantities are shown in the basis of the local triad, with the $\Xs$ component in blue and the $\V{D}^{(1/2)}$ component in red/yellow. In the top row, from left to right, we show $\ds{\V{U}_E}$, (negative) mismatch $-\V{m}$, and the corresponding $\dF$. Notice how the shapes of $\dF$ and $\V{m}$ match, which implies that $\dF(s) \approx c(\epsRS) \V{m}(s)$ away from the fiber boundaries. In the bottom row, from left to right, we show the resulting velocity perturbations $\Delta \Omega^\parallel$, \rev{$\partial_s \Delta \V{U}$}, and $\Xs \times \dOm$. In the fiber interior, $\Xs \times \dOm \approx -\V{m}$, and $\ds{\dU}(s) \approx \V{0}$. At the endpoints, $\partial_s \dU$ reaches values as negative as $-10$, so we truncate the $y$ axis to show the behavior in the fiber interior.}
\end{figure}

We first solve the mismatch problem\ \eqref{eq:mprobT} with $\kappa=1$ and $\gamma=0$ using our second-order method with $N=1/\epsRS$ points. In Fig.\ \ref{fig:PertEffects}, we plot $\V{m}$, $\dF$, $\Delta \V{U}$, $\Xs \times \dOm$, and $\ds{\dU}$ in the basis of the local triad for $\epsRS=\eps/L=0.005$. We observe several features: first, the mismatch $\V{m}$ is smooth until it reaches the endpoints, after which it jumps to a large value at the boundaries. In the fiber interior, the function $\dF$ is approximately a constant multiple of $\V{m}$, but at the endpoints it transitions to zero (as required by the boundary condition). Away from the fiber boundaries, the mismatch $\V{m}$ is almost entirely compensated for by the angular velocity perturbation $\Xs \times \dOm$, and $\ds{\dU}$ is small in the fiber interior, validating\ \eqref{eq:dUEst}.

Figure\ \ref{fig:PertEffects} also shows that $\Xs \times \dOm$ by itself does not compensate for the mismatch $\V{m}$ near the boundaries, since the dominant contribution from torque $\dF$ must be zero there. The mismatch $\V{m}$ near the boundaries must therefore be compensated for by $\ds{\dU}$, which means that $\Delta \V{U}$ is highly nonsmooth at the boundaries. Since $\dU$ comes from $\dF$, it is not surprising that $\dF$ is also not smooth at the boundaries, as it takes on its largest values near the endpoints, only to have to fall back to zero as required by the boundary conditions. 

\begin{figure}
\centering
\subfigure[$\kappa=1, \gamma=0$]{\label{fig:DeltaConv}
\includegraphics[width=0.48\textwidth]{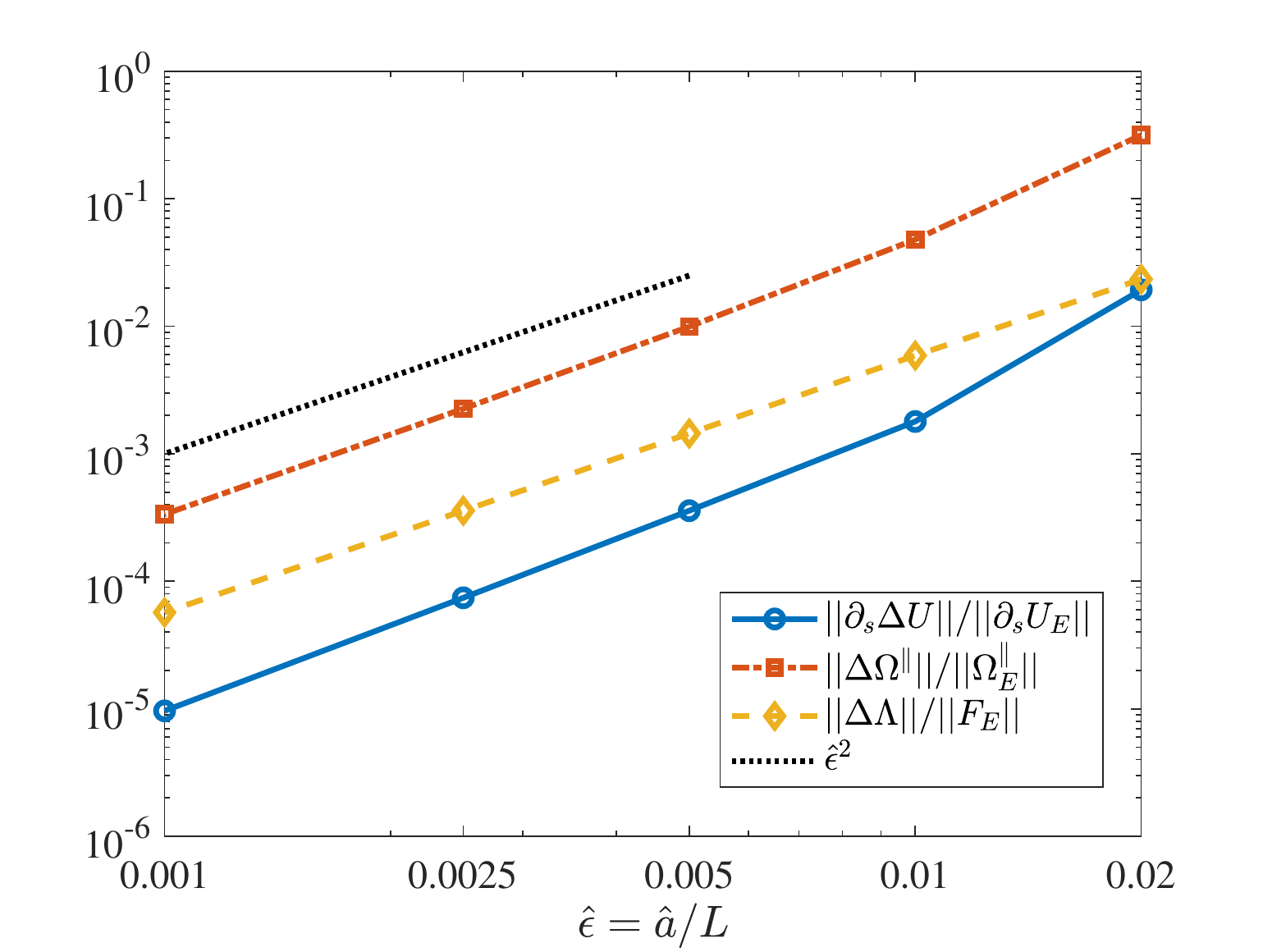}}
\subfigure[$\kappa=0, \gamma=1$]{\label{fig:DeltaConvT}
\includegraphics[width=0.48\textwidth]{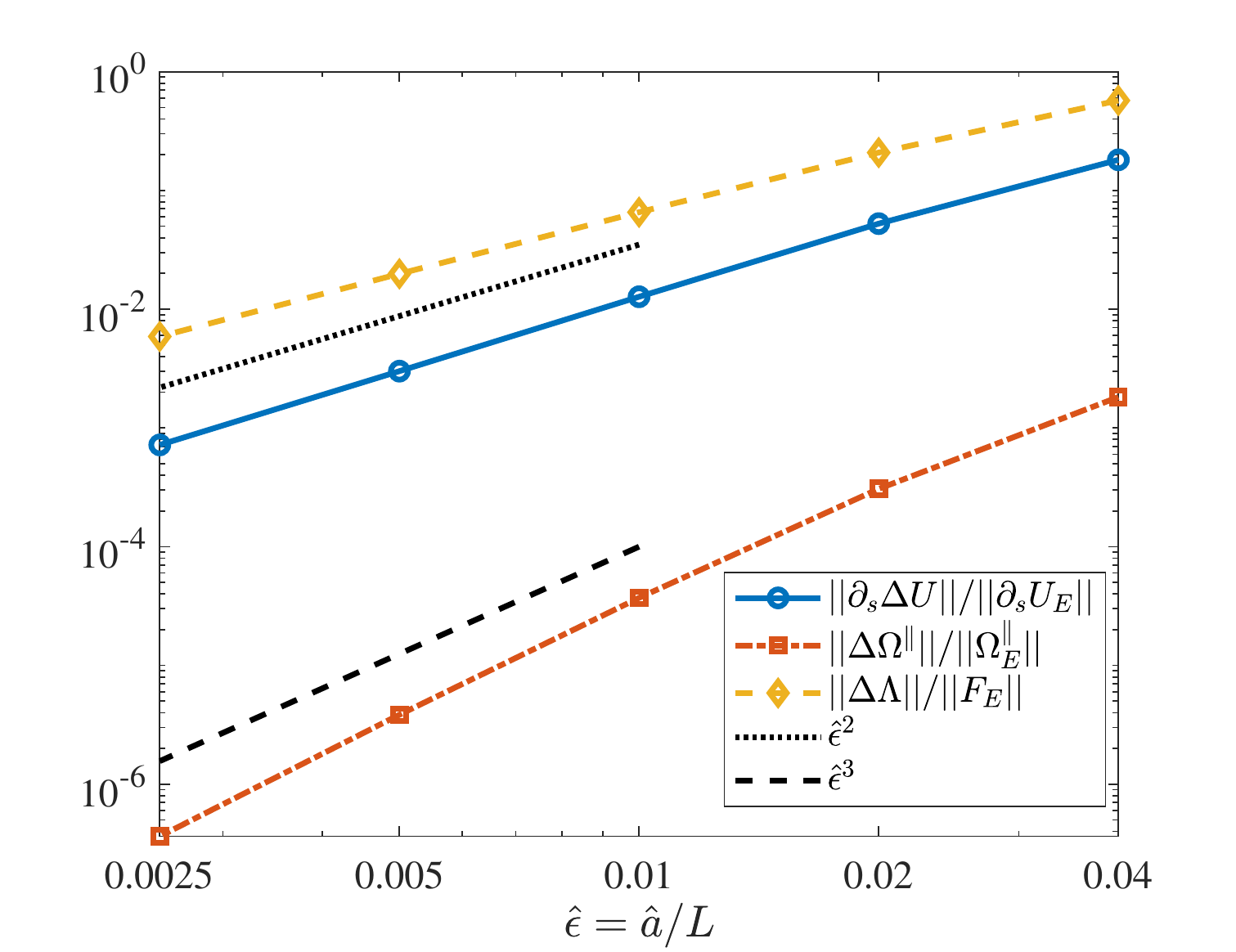}}
\caption{\label{fig:DeltaC}Decay of (properly normalized) magnitude of $\dF$ (yellow), $\ds{\Delta \V{U}}$ (blue), and $\Delta \Omega^\parallel$ (red) as slenderness $\epsRS \rightarrow 0$ for: (a) $\gamma=0$ and $\kappa=1$ with $N=1/\epsRS+1$ grid points; and, (b) $\gamma=1$ and $\kappa=0$ and  $N=4/\epsRS+1$ grid points. We use an interior ``norm'' ($L^2$ norm on $0.05 \leq s/L \leq 0.95$). All quantities scale like $\epsRS^2$ when properly normalized, with the exception of $\dOmpar$ when there is twist, which scales faster than $\epsRS^3$. This establishes\ \eqref{eq:dUEst} and\ \eqref{eq:dOmEst} numerically.  }
\end{figure}

We will focus on the equivalence of the Euler and Kirchhoff models in the fiber \emph{interior}, or points $\mathcal{O}(1)$ from the fiber endpoints. We therefore define an interior ``norm'' for all functions which  is simply the $L^2$ norm on $0.05 \leq s/L \leq 0.95$. In Fig.\ \ref{fig:DeltaConv}, we plot the relative difference in the velocity and forces using this interior norm. We observe convergence of order $\epsRS^2$ for all quantities. This verifies, and in fact exceeds, the predictions of\ \eqref{eq:dUEst} and\ \eqref{eq:dOmEst} of $\epsRS^2 \log{\epsRS}$. The reason is that the mismatch is actually dominated by $\mathcal{O}(1)$ terms for the values of $\epsRS$ we are interested in, and so our estimates are off by a factor of $\log{\epsRS}$.

\subsubsection{Twisting but no bending}
\begin{figure}
\hspace{-3em} 
\includegraphics[width=1.1\textwidth]{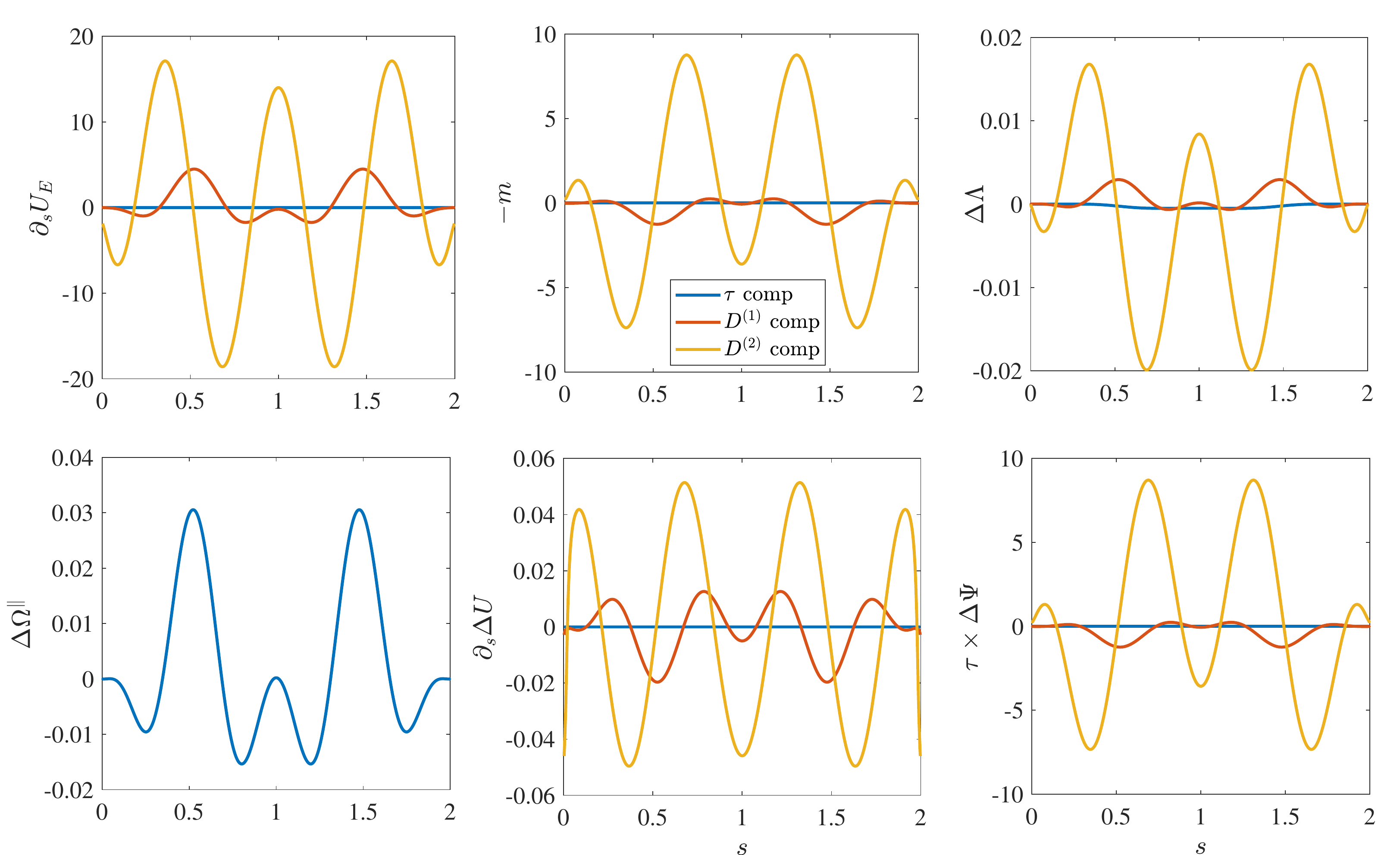}
\caption{\label{fig:MMTwist} Solutions of the mismatch problem\ \eqref{eq:mprobT} with twist modulus $\gamma=1$, bending modulus $\kappa=0$, and slenderness ratio $\epsRS=\eps/L =0.005$. The behavior is qualitatively the same as the mismatch problem without twist, with the perturbations in $\Delta \Omega^\parallel$ and $\partial_s \dU$ being small relative to the values obtained from the Euler method. Because the twist force $\V{f}^{(\gamma)}$ is zero at the endpoints, the mismatch there is small, and $\partial_s \dU$ is smoother than seen in Fig.\ \ref{fig:PertEffects}. }
\end{figure}

We now proceed to the case of finite twist modulus, $\twmod=1$, with zero bend modulus. We first plot the quantities related to the mismatch problem\ \eqref{eq:mprobT} in Fig.\ \ref{fig:MMTwist}. We see similar behavior to the case without twist and with bending force. The velocity difference $\ds{\dU}$ is again small, since most of the mismatch is compensated for by $\Xs \times \dOm$. A major difference from the non-twisted case is that $\dF$ and $\dU$ are smoother at the boundaries, since the force when $\kappa=0$, given in\ \eqref{eq:forceK}, is zero at the fiber boundaries because of the free fiber boundary conditions\ \eqref{eq:EulBC}. Because the force decays to zero there, the mismatch is much smaller at the boundaries (but is still not zero). This means that $\Xs \times \dOm$ can almost compensate for the entire mismatch everywhere. At the boundaries, whatever small mismatch there is must be compensated for by $\ds{\dU}$, and so as $\epsRS$ becomes smaller we do see some localized peaks in $\ds{\dU}$ near the endpoints.

In Figure\ \ref{fig:DeltaConvT}, we plot the convergence of the velocities and forces for $\twmod=1$ and $\kappa=0$. For $\dF$ and $\dU$, we observe a similar rate of convergence (approximately $\epsRS^2$) to Fig.\ \ref{fig:DeltaConv}. In accordance with\ \eqref{eq:dOmEst}, the parallel angular velocity $\Omega^\parallel_E$ is larger when there is twist, and so as a result the relative rate at which $\Delta \Omega^\parallel$ approaches zero is faster. In\ \eqref{eq:dOmEst} we estimated $\epsRS^4 \left(\log{\epsRS}\right)^2$; in Fig.\ \ref{fig:DeltaConvT} we simply show it is faster than $\epsRS^3$, but the line isn't straight, which indicates some kind of log correction. This verifies\ \eqref{eq:dUEst} and\ \eqref{eq:dOmEst} numerically for finite twist modulus.

 \delete{Although the fluid velocity in the Euler model \emph{does not} satisfy the velocity constraint\ \eqref{eq:simplest_cons}, it is the perpendicular components of the fluid angular velocity $\V{\Psi}$ that differ between the Euler and Kirchhoff models. The reason is that the perpendicular torques in the Kirchhoff model, $\Xs \times \V{\Lambda}^\perp$ generate $1/\epsc^2$ more angular velocity than translational velocity, and so their impact on the translational velocity $\V{U}$, which is used to compute $\V{\Omega}$ in the Euler model, is necessarily small. Thus the tangent vector rotation rates $\V{\Omega}$ become equivalent in the two models as $\epsc\rightarrow 0$. In Section\ \ref{sec:specEul} we develop a spectral numerical method based on the simpler Euler model. 
}

\setcounter{equation}{0}
\section{\change{Discretization of kinematic operators} \label{sec:DiscK}}
\change{In this appendix, we discuss the spectral discretization of the kinematic operators\ \eqref{eq:Kcdef} and\ \eqref{eq:KstPt} on a type 1 Chebyshev grid of size $N$. For simplicity of notation, we will consider $s=0$ to be a clamped end so that $\V{U}_0=\V 0$ and $\V \alpha \equiv \V \Omega$. The continuum operators $\Lop{K}$ and $\Lop{K}^*$ (c.f.\ \eqref{eq:Kcdef} and\ \eqref{eq:KstPt}) in this case are defined as
\begin{gather}
\label{eq:Kcdef2}
\V{U}(s)=\left(\Lop{K}\V{\alpha}\right)(s)=\int_0^{s} \V{\Omega}(s') \times \Xs(s') \, ds',\\
\label{eq:KstPt2}
\Lop{K}^*\V{\lambda}:=
\Xs(s) \times \int_s^L \V{\lambda}(s') \, ds'= \V 0.
\end{gather}
These operators must be discretized with care because of nonlinearities, which lead to aliasing on the spectral collocation grid. Our approach in this paper, which is somewhat different from what we have done previously \cite[Sec.~4.1.2]{maxian2021integral}, is to discretize these operators on a type 1 Chebyshev grid of size $2N$, and then appropriately upsample and downsample to a grid of size $N$. This is in contrast to \cite{maxian2021integral}, where we set $\Omega^\parallel(s)=0$ on the $N$ point grid. 

As in \cite[Sec.~4.1.2]{maxian2021integral}, we let $\M{J}$ be a matrix such that $\M{J}\V{\Omega}$ is an inextensible motion \mbox{$\int_0^s \V{\Omega}(s^\prime) \times \Xs(s^\prime) ds^\prime$} on a grid of size $2N$. Unlike \cite[Sec.~4.1.2]{maxian2021integral}, however, we will use $\M{J}$ to project off extensible motions on the $2N$ grid, rather than trying to write an explicit basis for inextensible motions on the $N$ point grid. To define $\M{J}$, we let $\M{C}(\Xs)$ be the matrix that computes cross products with $\Xs$ on the $2N$ grid, i.e., $\M{C}\V{x}=\Xs \times \V{x}$, with $\M{C}^T=-\M{C}$. We use the pseudoinverse of the Chebyshev differentiation matrix \cite[c.~6]{trefethen2000spectral} for integration, with $\M{D}^\dagger_0$ denoting integration from 0 to $s$, which has $L^2$ adjoint $\left(\M{D}^\dagger_0\right)^* = -\M{D}^\dagger_L$ (integration from $s$ to $L$). The discretization of $\M{J}$ and its $L^2$ adjoint $\M{J}^*$ on a grid of size $2N$ can then be written as
\begin{gather}
\label{eq:newJ}
\M{J} =\M{D}_0^\dagger \M{C}, \qquad \qquad \M{J}^* =\M{C}\M{D}_L^\dagger.
\end{gather}
Here $\M{J}$ is a $6N \times 6N$ matrix that discretizes\ \eqref{eq:Kcdef2}, and $\M{J}^*$ is a $6N \times 6N$ matrix that discretizes\ \eqref{eq:KstPt2}. 

To write a discretization on the $N$ point grid, we introduce an extension matrix $\M{E}$ and restriction matrix $\M{R}$. The extension matrix is well defined: to evaluate $\M{E}f_N$, where $f_N$ is a polynomial of degree $N-1$ on the $N$ point grid, we simply evaluate $f$ on the grid of size $2N$. Following \cite[Eq.~(62)]{maxian2021integral}, we let the restriction matrix $\M{R}$ be the matrix that gives the best possible (in the $L^2$ sense) polynomial approximation of degree $N-1$ on the $N$ grid to a polynomial of degree $2N-1$ on the $2N$ point grid, i.e., if $f_{2N}$ is a function on the $2N$ point grid, we want to choose $\M{R}$ such that $\norm{\M{E}\M{R}f_{2N}-f_{2N}}_{L^2}$ is minimized. This can be accomplished by setting
\begin{equation}
 \M{R} = \left(\M{E}^T \M{W}_{2N} \M{E}\right)^{-1} \M{E}^T \M{W}_{2N},
\end{equation}
where $\M{W}_{2N}$ is a diagonal matrix of Clenshaw-Curtis weights on the $2N$ point grid. We do this for both $\M{K}$ and $\M{K}^*$ to obtain the discretizations of $\Lop{K}$ and $\Lop{K}^*$,
\begin{equation}
\label{eq:Kmat}
\M{K} = \M{R}\M{J}\M{E}, \qquad \text{and} \qquad \M{K}^* = \M{R}\M{J}^*\M{E} 
\end{equation}
where $\M{K}$ and $\M{K}^*$ are both $3N \times 3N$ matrices. It is straightforward to modify $\M{K}$ and $\M{K}^*$ to allow a nonzero $\V{U}_0$ for a free end.


To update the fiber tangent vectors from the velocity $\V{U}=\M{K}\V{\alpha}$, we follow\ \eqref{eq:OmperpFromK}. This  update to the tangent vectors projects off any parallel components in $\V{U}_s$, which come from numerical discretization artifacts. Therefore, the parallel component of $\V{\Omega}$ on the $N$ point grid does not affect the dynamics; it is included merely to lessen aliasing artifacts.

As discussed in Section\ \ref{sec:specDisc}, the main difference between our discretization here and that of \cite[Sec.~4.1.2]{maxian2021integral} is the treatment of the null space of the continuum operator $\Lop{K}$. Because we project off the null space of $\Lop{K}$ on a grid of size $2N$ and downsample the result, the rank of the matrix $\M{K}$ (and $\M{K}^*$) is a function of the fiber curvature. For a rigid fiber, the rank of $\M{K}$ is $2N+1$, which is the same as when we directly write a basis for inextensible motions on an $N$ point grid. For curved fibers, the rank is larger because inextensible motions are aliased to extensible ones on the $N$ point grid. For the fiber\ \eqref{eq:Xsfibq} with $q=7$ for instance, for large $N$ we observe $2N+23$ nonzero singular values of $\M{K}$, so there are 22 additional positive singular values on top of the $2N+1$ nonzero singular values for a rigid fiber. Hence, typically the matrix $\M{K}$ does not have full rank, so the Schur complement $\M{K}^* \M{M}^{-1}\M{K}$ that we use to solve the saddle point system\ \eqref{eq:euler} is singular. We therefore use the pseudo-inverse of $\M{K}^* \M{M}^{-1}\M{K}$ to solve the saddle point system (\texttt{pinv} function in Matlab with default tolerance).}

\setcounter{equation}{0}
\section{Nearly singular quadrature \label{sec:quads}}
In this appendix, we discuss the nearly singular quadrature scheme we use to evaluate the integrals of the RPY kernels in Section\ \ref{sec:exactmob} to spectral accuracy. The idea behind this scheme has already been published in \cite{barLud, tornquad, maxian2021integral}, but in those works the nearly-singular integral is actually a singular finite part integral because the domain is $[0,L]$ instead of $D(s)$ as defined in\ \eqref{eq:Ddom}. This section presents a review of the quadrature scheme and how we modify it when the domain changes. We do this sequentially for the Stokeslet integral, doublet integral, and rotlet integral. 

\subsection{Stokeslet integral \label{sec:StNS}}
For the Stokeslet integral\ \eqref{eq:Stsubtr}, the nearly-singular integral that we must evaluate is the Stokeslet minus the leading order singularity, 
\begin{equation}
\label{eq:StNSint}
\tt{\V{U}}^\text{(int, S)}(s)=\int_{D(s)} \left(\Slet{\V{X}(s),\V{X}(s')} \V{f}\left(s'\right) -  \EPMI \left(\frac{\M{I}+\Xs(s)\Xs(s)}{|s-s'|}\right)\V{f}(s)\right) \, ds' ,
\end{equation}
which has the same integrand as the finite part integral in slender body theory, but with the different domain of integration
\begin{equation}
D(s) =\begin{cases} \left(0,s-2\eps\right) \cup \left(s+2\eps,L\right) & 2\eps \leq s \leq L-2\eps \\
\left(s+2\eps,L\right) & s < 2\eps \\
 \left(0,s-2\eps\right) & s > L-2\eps
\end{cases}.
\end{equation}
We employ a singularity subtraction scheme that ensures that the second term in the integral\ \eqref{eq:StNSint} cancels the leading order $1/|s'-s|$ singularity in the first. The next singularity is $\text{sign}(s'-s)$, which means that the near singular integral\ \eqref{eq:StNSint} can be written as
\begin{gather} 
\label{eq:FPre}
 \tt{\V{U}}^\text{(int, S)}(s) =\int_{D(s)} \V{g}_\text{tt}(s,s') \frac{s'-s}{|s'-s|} \, ds'
= \frac{L}{2}\int_{D'(\eta)} \V{g}_\text{tt}(\eta, \eta') \frac{\eta'-\eta}{|\eta'-\eta|} \, d\eta', \qquad \text{where}\\ 
\nonumber
\V{g}_\text{tt}(s,s') = \left(\Slet{\V{X}(s),\V{X}(s')} \V{f}\left(s'\right)  |s'-s| - \EPMI \left(\V{I}+\Xs(s) \Xs(s)\right) \V{f}(s)\right)\frac{1}{s'-s}
\end{gather}
is a smooth function, $\eta=-1+2s/L$ is a rescaled arclength coordinate on $[-1,1]$, and $D'(\eta)=[0,\eta_\ell] \cup [\eta_h,L]$ is defined from $D(s)$. The function $\V{g}_\text{tt}$ is nonsingular at $s=s'$ with the finite limit
\begin{equation}
\label{eq:gttlimit}
\lim_{s' \to s} \V{g}_\text{tt}(s,s') = \EPMI\left(\frac{1}{2}\left(\Xs(s)\ds{\Xs}(s)+\ds{\Xs}(s)\Xs(s)\right)\V{f}(s) + \left(\V{I}+\Xs(s)\Xs(s)\right)\ds{\V{f}}(s)\right). 
\end{equation}
Furthermore, $\V{g}_\text{tt}$ is smooth, so we can express it in a truncated Chebyshev series on $[-1,1]$, 
\begin{equation}
\label{eq:gmono}
\frac{L}{2}\V{g}_\text{tt}(\eta,\eta') \approx \sum_{k=0}^{N-1} \V{c}_k(\eta) T_k(\eta'), 
\end{equation}
where $\V{c}_k$ is a vector of 3 coefficients for each $\eta$. Substituting the Chebyshev  expansion\ \eqref{eq:gmono} into the integrand\ \eqref{eq:FPre}, we obtain
\begin{align}
\label{eq:expandmono}
 \tt{\V{U}}^\text{(int, S)}(\eta) &=\sum_{k=0}^{N-1} \V{c}_k(\eta) \int_{D'(\eta)} T_k\left(\eta^\prime \right) \frac{\eta'-\eta}{|\eta'-\eta|} \, d\eta' = \sum_{k=0}^{N-1} \V{c}_k(\eta) q^{(S)}_k(\eta)=\V{c}^T(\eta) \V{q}^{(S)}(\eta),\\[6 pt]
\label{eq:qRPY}
\text{where} \qquad 
q_k^{(S)}(\eta) & = \int_{D'(\eta)} T_k(\eta') \frac{\eta'-\eta}{|\eta'-\eta|} \, d\eta'= -\int_{-1}^{\eta_\ell} T_k(\eta') \, d\eta' +\int_{\eta_h}^L T_k(\eta') \, d\eta',
\end{align}
are integrals that can be precomputed to high accuracy for each $\eta$ on the Chebyshev collocation grid. The precomputed integrals in\ \eqref{eq:qRPY} can then be used in an adjoint method to accelerate the repeated evaluation of\ \eqref{eq:expandmono}. Specifically, the coefficients of $\V{g}_\text{tt}(\eta)$ for each $\eta$ can be written as $\V{c}(\eta) = \M{V}^{-1} \V{g}_\text{tt}(\eta,\eta')$, where $\M{V}$ is the matrix that maps coefficients of the basis functions to values on the collocation grid. Then $\V{c}^T(\eta) = \V{g}^T(\eta,\eta') \M{V}^{-T}$, and $ \tt{\V{U}}^\text{(int, S)}(\eta)=\V{g}^T (\eta,\eta')\M{V}^{-T}\V{q}^{(S)}(\eta)$. The matrix vector product $\M{V}^{-T}\V{q}^{(S)}(\eta)$ can then be precomputed for each $\eta$, and so at each time step all that is needed to compute the integral\ \eqref{eq:StNSint} for each $\eta$ are the values of $ \V{g}_\text{tt}(\eta,\eta')$ on the collocation grid for $\eta'$ being another collocation point. Evaluating $\tt{\V{U}}^\text{(int, S)}$ therefore has a cost of $\mathcal{O}\left(N^2\right)$. 
 
We note that our use of Chebyshev polynomials differs from the approach of \cite{tornquad, maxian2021integral}, which uses a monomial basis for the expansion\ \eqref{eq:gmono}. Although the integrals\ \eqref{eq:qRPY} with monomials can be computed analytically, using the adjoint method with monomial coefficients requires inverting the Vandermonde matrix $\V{V}$, which is too ill-conditioned with $N \gtrsim 40$ points in double precision. The Chebyshev expansion is an improvement in the sense that computing coefficients from values on the Chebyshev grid is a well-conditioned problem, and so the technique we use to evaluate finite part integrals does not limit the number of points on the fiber or require the use of panels (composite quadrature). 

\subsection{Doublet integral \label{sec:DbNS}}
The next part of the translational mobility\ \eqref{eq:transmob} is the integral of the doublet kernel. We begin by applying the same singularity subtraction technique we used in\ \eqref{eq:Stsubtr} for the Stokeslet integral
\begin{align}
\label{eq:DbSS}
\tt{\V{U}}^{(D)}(s)=&\int_{D(s)} \Dlet{\V{X}(s),\V{X}(s')} \V{f}\left(s'\right)\, ds'\\[2 pt] \nonumber
= \EPMI & \int_{D(s)} \left(\frac{\M{I}-3\Xs(s)\Xs(s)}{|s'-s|^3}\right)\V{f}(s) \, ds' \\ \nonumber + &\int_{D(s)} \left(\Dlet{\V{X}(s),\V{X}(s')} \V{f}\left(s'\right) -\EPMI \left(\frac{\M{I}-3\Xs(s)\Xs(s)}{|s'-s|^3}\right)\V{f}(s)\right) \, ds'\\[2 pt] \nonumber
 =& \tt{\V{U}}^{(\text{inner, D})}(s)+\tt{\V{U}}^\text{(int, D)}(s).
\end{align}
The term $\tt{\V{U}}^{(\text{inner, D})}(s)$ is given for all $s$ by
\begin{gather}
\label{eq:UinnerDb}
8 \pi \mu \tt{\V{U}}^{(\text{inner,D})}(s) =  
\frac{2}{3}\left(\M{I}-3\Xs(s)\Xs(s)\right)\V{f}(s)  
\begin{cases}
\dfrac{1}{4} -\dfrac{\eps^2}{2s^2}-\dfrac{\eps^2}{2(L-s)^2}& 2\eps < s < L-2\eps \\ 
\dfrac{1}{8}-\dfrac{\eps^2}{2(L-s)^2}& s \leq 2\eps \\ 
\dfrac{1}{8}- \dfrac{\eps^2}{2s^2}& s \geq L-2\eps
\end{cases} 
\end{gather}

The nearly-singular doublet integral is handled in exactly the same way as the Stokeslet one, but the singular behavior is different, namely,
\begin{gather}
\nonumber
 \tt{\V{U}}^\text{(int, D)}(s) =\int_{D(s)} \frac{(s'-s)}{|s'-s|^3} \, \V{g}_D(s',s) \, ds',\\[4 pt] 
\label{eq:dint}
\text{where} \qquad \V{g}_D(s',s)=\left(\Dlet{\V{X}(s),\V{X}(s')} \V{f}\left(s'\right) -\EPMI \left(\frac{\M{I}-3\Xs(s)\Xs(s)}{|s'-s|^3}\right)\V{f}(s)\right) \frac{|s'-s|^3}{(s'-s)}
\end{gather}
is a smooth function with the finite limit
\begin{equation}
\lim_{s'\rightarrow s} \V{g}_D(s',s) =\EPMI \left( \left(\M{I}-3\Xs(s)\Xs(s)\right)\V{f}'(s)+ \frac{1}{2}\left(\Xs(s) \ds{\Xs}(s)+ \ds{\Xs}(s)\Xs(s)\right) \V{f}(s)\right).
\end{equation}
If we expand $\V{g}_D$ in terms of Chebyshev polynomials, the efficient evaluation of\ \eqref{eq:dint} requires precomputing integrals of the form
\begin{equation}
\label{eq:intMod}
q_k^\text{(D)}(\eta) =\int_{D'(\eta)} \frac{(\eta'-\eta)}{|\eta'-\eta|^3} T_k\left(\eta'\right) \, d\eta',
\end{equation}
for all $\eta$ on the Chebyshev grid. Notice that the integrals\ \eqref{eq:intMod} are not defined if $\eta'=\eta$ is included in the integration domain, but $D'(\eta)$ does \emph{not} contain $\eta$. 

\subsection{Rot-trans mobility \label{sec:RTNS}}
The long-range integral for the rot-trans mobility\ \eqref{eq:rottransmob} is given by an integral of the rotlet kernel over the domain $D(s)$. Using singularity subtraction, we compute this integral as 
\begin{align}
\nonumber
\rt{\Omega}^\text{(R)}(s)=&\int_{D(s)} \frac{\left(\V{f}(s') \times \V{R}(s')\right) \cdot \Xs(s) }{R(s')^3} ds'\\[2 pt] \nonumber
= &\int_{D(s)}\left(\frac{\left(\ds{\Xs}(s) \times \V{f}(s)\right) \cdot \Xs(s)}{2|s-s'|} \right) \, ds'\\ \nonumber + & \int_{D(s)}  \left(\frac{\left(\V{f}(s') \times \V{R}(s')\right) \cdot \Xs(s) }{R(s')^3} - \frac{\left(\ds{\Xs}(s) \times \V{f}(s)\right) \cdot \Xs(s)}{2|s-s'|} \right) \, ds' \\[2 pt] \label{eq:Rtltsubtr}
:=&\rt{\Omega}^{(\text{inner, R})}(s)+ \rt{\Omega}^{(\text{int, R})}(s),
\end{align}
where $\V{R}(s') = \V{X}(s)-\V{X}(s')$ and $R=\norm{\V{R}}$. Since the leading order singularity in $\rt{\Omega}^{(\text{inner, R})}$ of $1/|s-s'|$ is the same as that of the Stokeslet, the expression for $\rt{\Omega}^{(\text{inner, R})}$ is also given by\ \eqref{eq:UinnerSt}, but with $\left(\M{I}+\Xs(s)\Xs(s)\right)\V{f}(s)$ replaced by $1/2\left(\ds{\Xs}(s) \times \V{f}(s)\right) \cdot \Xs(s)=1/2\left(\Xs(s) \times \ds{\Xs}(s)\right) \cdot \V{f}(s)$. 

The nearly singular integral $\rt{\Omega}^{(\text{int, R})}(s)$ in\ \eqref{eq:Rtltsubtr} can then be put in the form\ \eqref{eq:FPre} (with scalar $g_\text{rt}(s,s')$) by setting
\begin{equation}
\label{eq:grt}
g_\text{rt}(s,s') = \left(\frac{|s-s'|}{R(s')^3}\left(\V{R}(s') \times \Xs(s)\right) \cdot \V{f}(s') -\frac{\left(\Xs(s) \times \ds{\Xs}(s)\right) \cdot \V{f}(s) }{2}\right)\frac{1}{s'-s},
\end{equation}
which has the finite limit
\begin{equation}
\label{eq:grtlimit}
 \lim_{s' \rightarrow s} g_\text{rt}(s,s') = -\frac{1}{2}\ds{\V{f}}(s)\cdot \left(\ds{\Xs}(s) \times \Xs(s) \right)-\frac{1}{6} \V{f}(s) \cdot \left(\ds^2{\Xs}(s) \times \Xs(s) \right).
\end{equation}
Now the scalar function $g_\text{rt}(s,s')$ can be expanded in a Chebyshev series as in\ \eqref{eq:expandmono}, and the precomputed integrals $q^{(S)}_k(\eta)$ in\ \eqref{eq:qRPY} used once again to evaluate the second integral in\ \eqref{eq:Rtltsubtr}. 

\subsection{Trans-rot mobility \label{sec:TRNS}}
For the trans-rot mobility\ \eqref{eq:transrotmob}, the singularity subtraction approach is to set 
\begin{align}
\nonumber
\tr{\V{U}}^{(\text{R})}(s)=&\int_{D(s)} \frac{\left(\Xs(s') \times \V{R}(s')\right)n^\parallel(s')}{R(s')^3} ds'\\[2 pt] \nonumber
= &\int_{D(s)}\left(\frac{\left(\Xs(s) \times \ds{\Xs}(s)\right) n^\parallel(s)}{2|s-s'|} \right) \, ds'\\ \nonumber + &\int_{D(s)}  \left(\frac{\left(\Xs(s') \times \V{R}(s')\right)n^\parallel(s') }{R(s')^3} - \frac{\left(\Xs(s) \times \ds{\Xs}(s)\right) n^\parallel(s)}{2|s-s'|}\right) \, ds' \\[2 pt] \label{eq:Rtltsubtr2}
:=&\tr{\V{U}}^{(\text{inner, R})}(s)+\tr{\V{U}}^{(\text{int, R})}(s).
\end{align}
where $\tr{\V{U}}^{(\text{inner, R})}$ is given by\ \eqref{eq:UinnerSt}, but with $\left(\M{I}+\Xs(s)\Xs(s)\right)\V{f}(s)$ replaced by $1/2\left(\Xs(s) \times \ds{\Xs}(s) \right) n^\parallel(s)$. The nearly singular integral $\tr{\V{U}}^{(\text{int, R})}(s)$ can then be put in the form\ \eqref{eq:FPre} by setting
\begin{align}
\label{eq:defg}
\V{g}_\text{tr}(s,s') &=  \left(\frac{|s-s'|}{R(s')^3}\left(\Xs(s') \times \V{R}(s')\right)n^\parallel(s') -\left(\frac{\Xs(s) \times \ds{\Xs}(s) }{2}\right)n^\parallel(s)\right)\frac{1}{s'-s},\\
 \text{where} \qquad &\lim_{s' \rightarrow s} \V{g}_\text{tr}(s,s') = -\frac{1}{2}\ds{n}^\parallel(s)\left(\ds{\Xs}(s) \times \Xs(s)\right) - \frac{1}{3}n^\parallel(s)\left(\ds^2{\Xs}(s) \times \Xs(s) \right).
\end{align}

\subsection{Rot-rot mobility}
Although we use asymptotic evaluation for the rot-rot mobility as discussed in Section\ \ref{sec:RPYasymp}, it is straightforward to derive a nearly singular quadrature scheme to evaluate the RPY integral for $\Mcrr$. Substituting the tensor $\Mbrr$ from\ \eqref{eq:MbrrRPY} into the integral\ \eqref{eq:PsiMobs} and using the approximation $\norm{\V{X}(s)-\V{X}(s')} \approx s-s'$ when $s-s'=\mathcal{O}(\eps)$, the mobility is given by
\begin{gather}
\label{eq:rfromr}
8\pi \mu \rr{\V{\Psi}}(s) = -\frac{1}{2}\int_{D(s)} \Dlet{\V{X}(s),\V{X}(s')} \V{n}(s')  \, ds'\\[2 pt] 
\nonumber
+\frac{1}{\eps^3}\int_{D^c(s)} \left(\left(1-\frac{27R(s')}{32\eps}+\frac{5 R(s')^3}{64 \eps^3}\right)\M{I}+\left(\frac{9}{32\eps R(s')}-\frac{3R(s')}{64\eps^3}\right)\left(\V{R}\V{R}\right)(s')\right)\V{n}(s') \, ds'.
\end{gather}
The first of these integrals can be evaluated to high accuracy using the doublet scheme of Section\ \ref{sec:DbNS}, while the second integral can be computed using direct Gauss-Legendre quadrature on $D^c(s)$ (split into two pieces). The parallel angular velocity $\rr{\Psi}^\parallel$ can then be obtained by taking the inner product of $\rr{\V{\Psi}}$ with $\Xs$. 

\subsection{Convergence of quadratures \label{sec:convquad}}
\begin{figure}
\centering
\includegraphics[width=0.49\textwidth]{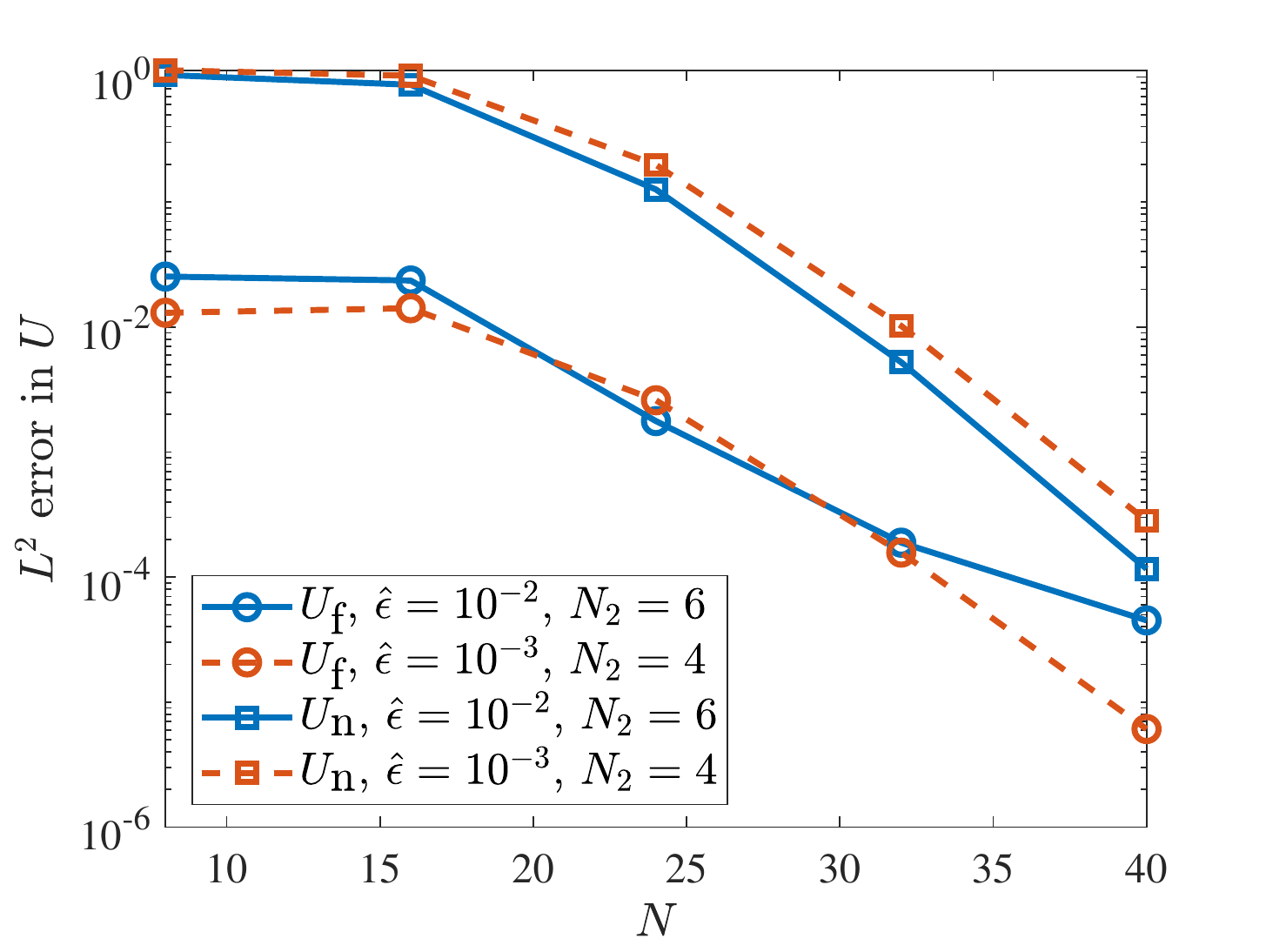}
\includegraphics[width=0.49\textwidth]{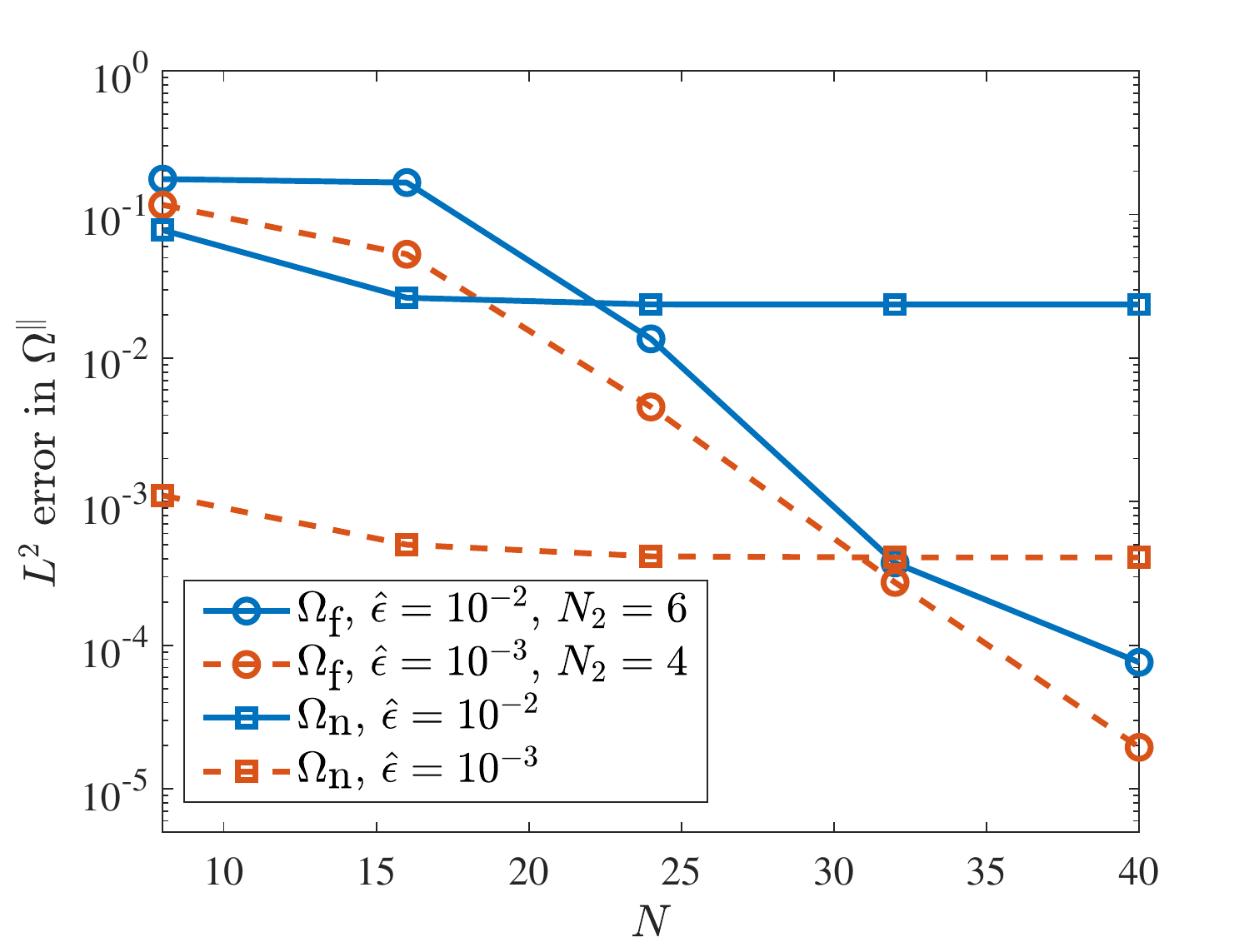}
\caption{\label{fig:QuadConv} Convergence of the nearly singular quadrature schemes for the (left) trans-trans mobility\ \eqref{eq:transmob} and trans-rot mobility\ \eqref{eq:transrotmob}, and (right) rot-trans mobility\ \eqref{eq:rottransmob} and rot-rot mobility\ \eqref{eq:rotrotmob} (here we are comparing the the full RPY integral of doublets as in\ \eqref{eq:PsiMobs}). We use the fiber\ \eqref{eq:Xsfibq} with $q=7$ and compute force and torque on the collocation grid of $N$ points using\ \eqref{eq:fkapd},\ \eqref{eq:nparBC}, and\ \eqref{eq:ftwmodC} with $\kappa=\gamma=1$. We then compute the velocities $\tt{\V{U}}$, $\tr{\V{U}}$, $\rt{\Omega}^\parallel$, and $\rr{\Omega}^\parallel$ on the collocation grid by applying the discrete mobility matrices given in Section\ \ref{sec:exactmob}. In all cases, we compute a reference solution by upsampled direct quadrature and report the \emph{relative} $L^2$ errors in the velocity.}
\end{figure}

We next test the accuracy of the (action of the) mobility matrices $\Mtt$, $\Mtr$, $\Mrt$, and $\Mrr$ by comparing against upsampled direct quadrature. We use the fiber\ \eqref{eq:Xsfibq} with $q=7$, set up a collocation grid of $N$ points, and compute the force $\V{f}=\V{f}^{(\kappa)}+\V{f}^{(\twmod)}$ and torque density $n^\parallel$ using\ \eqref{eq:fkapd},\ \eqref{eq:nparBC}, and\ \eqref{eq:ftwmodC} with $\kappa=\gamma=1$. We apply the forward mobility operators as discretized in Section\ \ref{sec:exactmob} to obtain the (translational and rotational) velocity at each of the collocation points, and compare the result to a reference solution which uses direct oversampled quadrature. The relative error is computed as the difference between the two, normalized by the $L^2$ norm of the reference solution. 

Figure\ \ref{fig:QuadConv} shows the results for the trans-trans mobility\ \eqref{eq:transmob} and trans-rot mobility\ \eqref{eq:transrotmob} (left), and rot-trans mobility\ \eqref{eq:rottransmob} and rot-rot mobility\ \eqref{eq:rotrotmob} (right) for two different values of $\epsRS$, and various values of $N$ (see Section\ \ref{sec:exactmob} for discretization details). In the case of the trans-trans, rot-trans, and trans-rot mobilities, we obtain spectral convergence and about 3 digits of accuracy when $N=32$. For $\epsRS=10^{-2}$, $N_2$ must be at least 6 to see this accuracy, but for $\epsRS=10^{-3}$ a smaller $N_2\approx 4$ is acceptable since the relative size of the domain $|s-s'| \leq 2\eps$ shrinks by a factor of 10. For $\epsRS=10^{-2}$, we previously found $N_2 \geq 8$ to give positive eigenvalues in the mobility, so accuracy is a weaker restriction than the eigenvalues of $\Mtt$ being strictly positive. 

In Fig.\ \ref{fig:QuadConv}, we also show the accuracy of the rot-rot mobility as defined in\ \eqref{eq:rotrotmob} relative to the actual RPY integral of doublets as defined in\ \eqref{eq:PsiMobs} and\ \eqref{eq:MbrrRPY}. This ``convergence'' plot is fundamentally different from the others, since in this case we are measuring the error of the asymptotic formula \emph{and} its approximation by a finite number of collocation points, rather than just the latter. Therefore, the relative error saturates around $N=16$ at a value that strongly depends on $\epsRS$. For $\epsRS=10^{-2}$, the relative error saturates around $2 \times 10^{-2}$, while for $\epsRS=10^{-3}$ the saturated relative error is about $4 \times 10^{-4}$. This reflects the fact that the full asymptotic formula for the rot-rot mobility \cite[Appendix~C.3]{TwistSBT} is $\mathcal{O}(\epsRS^2)$ in the fiber interior ($2\eps \leq s \leq L-2\eps$) and $\mathcal{O}(\epsRS)$ at the fiber endpoints. 

\end{appendices}

\bibliographystyle{plain}

\begin{thebibliography}{10}

\bibitem{barLud}
Ludvig af~Klinteberg and Alex~H Barnett.
\newblock Accurate quadrature of nearly singular line integrals in two and
  three dimensions by singularity swapping.
\newblock {\em BIT Numerical Mathematics}, pages 1--36, 2020.

\bibitem{af2017fast}
Ludvig af~Klinteberg, Davoud~Saffar Shamshirgar, and Anna-Karin Tornberg.
\newblock Fast ewald summation for free-space stokes potentials.
\newblock {\em Research in the Mathematical Sciences}, 4(1):1--32, 2017.

\bibitem{agarwal2017analysis}
Bhagwan~D Agarwal, Lawrence~J Broutman, and K~Chandrashekhara.
\newblock {\em Analysis and performance of fiber composites}.
\newblock John Wiley \& Sons, 2017.

\bibitem{alberts}
Bruce Alberts, Alexander Johnson, Julian Lewis, Martin Raff, Keith Roberts, and
  Peter Walter.
\newblock {\em Molecular biology of the cell}.
\newblock Garland Science, 2002.

\bibitem{andersson2021integral}
Helge~I Andersson, Elena Celledoni, Laurel Ohm, Brynjulf Owren, and Benjamin~K
  Tapley.
\newblock An integral model based on slender body theory, with applications to
  curved rigid fibers.
\newblock {\em Physics of Fluids}, 33(4):041904, 2021.

\bibitem{tref17}
Jared~L Aurentz and Lloyd~N Trefethen.
\newblock Block operators and spectral discretizations.
\newblock {\em SIAM Review}, 59(2):423--446, 2017.

\bibitem{balboa2017hydrodynamics}
Florencio Balboa~Usabiaga, Bakytzhan Kallemov, Blaise Delmotte, Amneet Bhalla,
  Boyce Griffith, and Aleksandar Donev.
\newblock Hydrodynamics of suspensions of passive and active rigid particles: a
  rigid multiblob approach.
\newblock {\em Communications in Applied Mathematics and Computational
  Science}, 11(2):217--296, 2017.

\bibitem{bergou2008discrete}
Mikl{\'o}s Bergou, Max Wardetzky, Stephen Robinson, Basile Audoly, and Eitan
  Grinspun.
\newblock Discrete elastic rods.
\newblock In {\em ACM SIGGRAPH 2008 papers}, pages 1--12. 2008.

\bibitem{blum1979biophysics}
Jacob~J Blum and Michael Hines.
\newblock Biophysics of flagellar motility.
\newblock {\em Quarterly reviews of biophysics}, 12(2):103--180, 1979.

\bibitem{brady1988stokesian}
John~F Brady and Georges Bossis.
\newblock Stokesian dynamics.
\newblock {\em Annual review of fluid mechanics}, 20(1):111--157, 1988.

\bibitem{brennen1977fluid}
Christopher Brennen and Howard Winet.
\newblock Fluid mechanics of propulsion by cilia and flagella.
\newblock {\em Annual Review of Fluid Mechanics}, 9(1):339--398, 1977.

\bibitem{ttbring08}
Thomas~T Bringley and Charles~S Peskin.
\newblock Validation of a simple method for representing spheres and slender
  bodies in an immersed boundary method for stokes flow on an unbounded domain.
\newblock {\em Journal of Computational Physics}, 227(11):5397--5425, 2008.

\bibitem{bruss2019twirling}
Isaac~R Bruss, Heena~K Mutha, Katherine Stoll, Brent Collins, Vinh Nguyen,
  David~JD Carter, Michael~P Brenner, and Kasey~J Russell.
\newblock Twirling, whirling, and tensioning: Plectoneme formation and
  suppression in flexible filaments.
\newblock {\em Physical Review Research}, 1(3):032020, 2019.

\bibitem{claessens2008helical}
Mireille Maria Anna~Elisabeth Claessens, C~Semmrich, L~Ramos, and AR~Bausch.
\newblock Helical twist controls the thickness of f-actin bundles.
\newblock {\em Proceedings of the National Academy of Sciences},
  105(26):8819--8822, 2008.

\bibitem{cortez2001method}
Ricardo Cortez.
\newblock The method of regularized stokeslets.
\newblock {\em SIAM Journal on Scientific Computing}, 23(4):1204--1225, 2001.

\bibitem{cortez2005method}
Ricardo Cortez, Lisa Fauci, and Alexei Medovikov.
\newblock The method of regularized stokeslets in three dimensions: analysis,
  validation, and application to helical swimming.
\newblock {\em Physics of Fluids}, 17(3):031504, 2005.

\bibitem{cortez2012slender}
Ricardo Cortez and Michael Nicholas.
\newblock Slender body theory for stokes flows with regularized forces.
\newblock {\em Communications in Applied Mathematics and Computational
  Science}, 7(1):33--62, 2012.

\bibitem{de1975low}
NJ~De~Mestre and WB~Russel.
\newblock Low-reynolds-number translation of a slender cylinder near a plane
  wall.
\newblock {\em Journal of Engineering Mathematics}, 9(2):81--91, 1975.

\bibitem{delong2013temporal}
Steven Delong, Boyce~E Griffith, Eric Vanden-Eijnden, and Aleksandar Donev.
\newblock Temporal integrators for fluctuating hydrodynamics.
\newblock {\em Physical Review E}, 87(3):033302, 2013.

\bibitem{dill1992kirchhoff}
Ellis~Harold Dill.
\newblock Kirchhoff's theory of rods.
\newblock {\em Archive for History of Exact Sciences}, pages 1--23, 1992.

\bibitem{dhale15}
Tobin~A Driscoll and Nicholas Hale.
\newblock Rectangular spectral collocation.
\newblock {\em IMA Journal of Numerical Analysis}, 36(1):108--132, 2015.

\bibitem{du2019dynamics}
Olivia Du~Roure, Anke Lindner, Ehssan~N Nazockdast, and Michael~J Shelley.
\newblock Dynamics of flexible fibers in viscous flows and fluids.
\newblock {\em Annu. Rev. Fluid. Mech.}, 51:539--572, 2019.

\bibitem{PSRPY}
Andrew~M Fiore, Florencio Balboa~Usabiaga, Aleksandar Donev, and James~W Swan.
\newblock Rapid sampling of stochastic displacements in brownian dynamics
  simulations.
\newblock {\em The Journal of chemical physics}, 146(12):124116, 2017.

\bibitem{fiore2018rapid}
Andrew~M Fiore and James~W Swan.
\newblock Rapid sampling of stochastic displacements in brownian dynamics
  simulations with stresslet constraints.
\newblock {\em The Journal of chemical physics}, 148(4):044114, 2018.

\bibitem{floyd2021stretching}
Carlos Floyd, Haoran Ni, Ravinda Gunaratne, Radek Erban, and Garegin~A Papoian.
\newblock On stretching, bending, shearing and twisting of actin filaments i:
  Variational models.
\newblock {\em arXiv preprint arXiv:2112.01480}, 2021.

\bibitem{gargslender}
Mohit Garg and Ajeet Kumar.
\newblock A slender body theory for the motion of special cosserat filaments in
  stokes flow.

\bibitem{gimbutas2015computational}
Zydrunas Gimbutas and Leslie Greengard.
\newblock Computational software: Simple fmm libraries for electrostatics, slow
  viscous flow, and frequency-domain wave propagation.
\newblock {\em Communications in Computational Physics}, 18(2):516--528, 2015.

\bibitem{gotz2001interactions}
Thomas G{\"o}tz.
\newblock {\em Interactions of fibers and flow: asymptotics, theory and
  numerics}.
\newblock dissertation. de, 2001.

\bibitem{huang2021numerical}
Weicheng Huang and M~Khalid Jawed.
\newblock Numerical simulation of bundling of helical elastic rods in a viscous
  fluid.
\newblock {\em Computers \& Fluids}, page 105038, 2021.

\bibitem{ishimoto2018elasto}
Kenta Ishimoto and Eamonn~A Gaffney.
\newblock An elastohydrodynamical simulation study of filament and spermatozoan
  swimming driven by internal couples.
\newblock {\em IMA Journal of Applied Mathematics}, 83(4):655--679, 2018.

\bibitem{inexRS}
Mehdi Jabbarzadeh and Henry~C Fu.
\newblock A numerical method for inextensible elastic filaments in viscous
  fluids.
\newblock {\em Journal of Computational Physics}, page 109643, 2020.

\bibitem{johnson}
Robert~E Johnson.
\newblock An improved slender-body theory for stokes flow.
\newblock {\em Journal of Fluid Mechanics}, 99(2):411--431, 1980.

\bibitem{kallemov2016immersed}
Bakytzhan Kallemov, Amneet Bhalla, Boyce Griffith, and Aleksandar Donev.
\newblock An immersed boundary method for rigid bodies.
\newblock {\em Communications in Applied Mathematics and Computational
  Science}, 11(1):79--141, 2016.

\bibitem{keaveny2011applying}
Eric~E Keaveny and Michael~J Shelley.
\newblock Applying a second-kind boundary integral equation for surface
  tractions in stokes flow.
\newblock {\em Journal of Computational Physics}, 230(5):2141--2159, 2011.

\bibitem{krub}
Joseph~B Keller and Sol~I Rubinow.
\newblock Slender-body theory for slow viscous flow.
\newblock {\em Journal of Fluid Mechanics}, 75(4):705--714, 1976.

\bibitem{koens2022tubular}
Lyndon Koens.
\newblock Tubular-body theory for viscous flows.
\newblock {\em Physical Review Fluids}, 7(3):034101, 2022.

\bibitem{varibook}
Cornelius Lanczos.
\newblock {\em The variational principles of mechanics}.
\newblock Courier Corporation, 2012.

\bibitem{lauga2009hydrodynamics}
Eric Lauga and Thomas~R Powers.
\newblock The hydrodynamics of swimming microorganisms.
\newblock {\em Reports on Progress in Physics}, 72(9):096601, 2009.

\bibitem{lee2014nonlinear}
Wanho Lee, Yongsam Kim, Sarah~D Olson, and Sookkyung Lim.
\newblock Nonlinear dynamics of a rotating elastic rod in a viscous fluid.
\newblock {\em Physical Review E}, 90(3):033012, 2014.

\bibitem{lim2008dynamics}
Sookkyung Lim, Anca Ferent, X~Sheldon Wang, and Charles~S Peskin.
\newblock Dynamics of a closed rod with twist and bend in fluid.
\newblock {\em SIAM Journal on Scientific Computing}, 31(1):273--302, 2008.

\bibitem{lim2004simulations}
Sookkyung Lim and Charles~S Peskin.
\newblock Simulations of the whirling instability by the immersed boundary
  method.
\newblock {\em SIAM Journal on Scientific Computing}, 25(6):2066--2083, 2004.

\bibitem{lindbo2011fast}
Dag Lindbo and Anna-Karin Tornberg.
\newblock Fast and spectrally accurate summation of 2-periodic stokes
  potentials.
\newblock {\em arXiv preprint arXiv:1111.1815}, 2011.

\bibitem{fcm03}
Sune Lomholt and Martin~R Maxey.
\newblock Force-coupling method for particulate two-phase flow: Stokes flow.
\newblock {\em J. Comput. Phys.}, 184(2):381--405, 2003.

\bibitem{ma2018structural}
Rui Ma and Julien Berro.
\newblock Structural organization and energy storage in crosslinked actin
  assemblies.
\newblock {\em PLoS computational biology}, 14(5):e1006150, 2018.

\bibitem{man2017bundling}
Yi~Man, William Page, Robert~J Poole, and Eric Lauga.
\newblock Bundling of elastic filaments induced by hydrodynamic interactions.
\newblock {\em Physical Review Fluids}, 2(12):123101, 2017.

\bibitem{TwistSBT}
Ondrej Maxian and Aleksandar Donev.
\newblock Slender body theories for rotating filaments.
\newblock {\em arXiv preprint arXiv:2203.12059}, 2022.

\bibitem{maxian2021bundling}
Ondrej Maxian, Aleksandar Donev, and Alex Mogilner.
\newblock Interplay between brownian motion and cross-linking kinetics controls
  bundling dynamics in actin networks.
\newblock {\em Biophysical Journal}, 2022.

\bibitem{maxian2021integral}
Ondrej Maxian, Alex Mogilner, and Aleksandar Donev.
\newblock Integral-based spectral method for inextensible slender fibers in
  stokes flow.
\newblock {\em Physical Review Fluids}, 6(1):014102, 2021.

\bibitem{maxian2021simulations}
Ondrej Maxian, Raul~P Pel{\'a}ez, Alex Mogilner, and Aleksandar Donev.
\newblock Simulations of dynamically cross-linked actin networks: morphology,
  rheology, and hydrodynamic interactions.
\newblock {\em PLOS Computational Biology}, 17(12):e1009240, 2021.

\bibitem{mori2020accuracy}
Yoichiro Mori and Laurel Ohm.
\newblock Accuracy of slender body theory in approximating force exerted by
  thin fiber on viscous fluid.
\newblock {\em Studies in Applied Mathematics}, 2021.

\bibitem{ehssan17}
Ehssan Nazockdast, Abtin Rahimian, Denis Zorin, and Michael Shelley.
\newblock A fast platform for simulating semi-flexible fiber suspensions
  applied to cell mechanics.
\newblock {\em J. Comput. Phys.}, 329:173--209, 2017.

\bibitem{nguyen2018impacts}
Frank~TM Nguyen and Michael~D Graham.
\newblock Impacts of multiflagellarity on stability and speed of bacterial
  locomotion.
\newblock {\em Physical Review E}, 98(4):042419, 2018.

\bibitem{oberbeck1876uber}
Anton Oberbeck.
\newblock Uber stationare flussigkeitsbewegungen mit berucksichtigung der inner
  reibung.
\newblock {\em J. reine angew. Math.}, 81:62--80, 1876.

\bibitem{olson2013modeling}
Sarah~D Olson, Sookkyung Lim, and Ricardo Cortez.
\newblock Modeling the dynamics of an elastic rod with intrinsic curvature and
  twist using a regularized stokes formulation.
\newblock {\em Journal of Computational Physics}, 238:169--187, 2013.

\bibitem{pareschi2005implicit}
Lorenzo Pareschi and Giovanni Russo.
\newblock Implicit--explicit runge--kutta schemes and applications to
  hyperbolic systems with relaxation.
\newblock {\em Journal of Scientific computing}, 25(1):129--155, 2005.

\bibitem{peskin1972flow}
Charles~S Peskin.
\newblock Flow patterns around heart valves: a numerical method.
\newblock {\em J. Comput. Phys.}, 10(2):252--271, 1972.

\bibitem{peskin2002acta}
Charles~S Peskin.
\newblock The immersed boundary method.
\newblock {\em Acta Numer.}, 11:479--517, 2002.

\bibitem{powers2010dynamics}
Thomas~R Powers.
\newblock Dynamics of filaments and membranes in a viscous fluid.
\newblock {\em Reviews of Modern Physics}, 82(2):1607, 2010.

\bibitem{reichert2006hydrodynamic}
Michael Reichert.
\newblock {\em Hydrodynamic interactions in colloidal and biological systems}.
\newblock PhD thesis, 2006.

\bibitem{RPYOG}
Jens Rotne and Stephen Prager.
\newblock Variational treatment of hydrodynamic interaction in polymers.
\newblock {\em J. Chem. Phys.}, 50(11):4831--4837, 1969.

\bibitem{rubin2013cosserat}
Mordecai~B Rubin.
\newblock {\em Cosserat theories: shells, rods and points}, volume~79.
\newblock Springer Science \& Business Media, 2013.

\bibitem{keavRPY}
Simon~F Schoeller, Adam~K Townsend, Timothy~A Westwood, and Eric~E Keaveny.
\newblock Methods for suspensions of passive and active filaments.
\newblock {\em Journal of Computational Physics}, 424:109846, 2021.

\bibitem{shelley1996nonlocal}
Michael~J Shelley and Tetsuji Ueda.
\newblock The nonlocal dynamics of stretching, buckling filaments.
\newblock {\em Advances in multi-fluid flows}, pages 415--425, 1996.

\bibitem{shelley2000stokesian}
Michael~J Shelley and Tetsuji Ueda.
\newblock The stokesian hydrodynamics of flexing, stretching filaments.
\newblock {\em Physica D: Nonlinear Phenomena}, 146(1-4):221--245, 2000.

\bibitem{swan2007simulation}
James~W Swan and John~F Brady.
\newblock Simulation of hydrodynamically interacting particles near a no-slip
  boundary.
\newblock {\em Physics of Fluids}, 19(11):113306, 2007.

\bibitem{tornquad}
Anna-Karin Tornberg.
\newblock Accurate evaluation of integrals in slender-body formulations for
  fibers in viscous flow.
\newblock {\em arXiv preprint arXiv:2012.12585}, 2020.

\bibitem{ts04}
Anna-Karin Tornberg and Michael~J Shelley.
\newblock Simulating the dynamics and interactions of flexible fibers in stokes
  flows.
\newblock {\em Journal of Computational Physics}, 196(1):8--40, 2004.

\bibitem{trefethen2000spectral}
Lloyd~N Trefethen.
\newblock {\em Spectral methods in MATLAB}, volume~10.
\newblock Siam, 2000.

\bibitem{vogel2012bacterial}
Reinhard Vogel.
\newblock The bacterial flagellum: Modeling the dynamics of the elastic
  filament and its transition between polymorphic helical forms.
\newblock 2012.

\bibitem{wada2006non}
Hirofumi Wada and Roland~R Netz.
\newblock Non-equilibrium hydrodynamics of a rotating filament.
\newblock {\em EPL (Europhysics Letters)}, 75(4):645, 2006.

\bibitem{wada2007model}
Hirofumi Wada and Roland~R Netz.
\newblock Model for self-propulsive helical filaments: kink-pair propagation.
\newblock {\em Physical review letters}, 99(10):108102, 2007.

\bibitem{wajnryb2013generalization}
Eligiusz Wajnryb, Krzysztof~A Mizerski, Pawel~J Zuk, and Piotr Szymczak.
\newblock Generalization of the rotne--prager--yamakawa mobility and shear
  disturbance tensors.
\newblock {\em Journal of Fluid Mechanics}, 731, 2013.

\bibitem{walker2020regularised}
Benjamin~J Walker and Eamonn~A Gaffney.
\newblock Regularised non-uniform segments and efficient no-slip
  elastohydrodynamics.
\newblock {\em Journal of Fluid Mechanics}, 915, 2021.

\bibitem{westwood2021coordinated}
Timothy~A Westwood and Eric~E Keaveny.
\newblock Coordinated motion of active filaments on spherical surfaces.
\newblock {\em arXiv preprint arXiv:2109.13578}, 2021.

\bibitem{wolgemuth2000twirling}
Charles~W Wolgemuth, Thomas~R Powers, and Raymond~E Goldstein.
\newblock Twirling and whirling: Viscous dynamics of rotating elastic
  filaments.
\newblock {\em Physical Review Letters}, 84(7):1623, 2000.

\bibitem{yan2018universal}
Wen Yan and Michael Shelley.
\newblock Universal image systems for non-periodic and periodic stokes flows
  above a no-slip wall.
\newblock {\em Journal of Computational Physics}, 375:263--270, 2018.

\bibitem{RPY_Shear_Wall}
PJ~Zuk, E~Wajnryb, KA~Mizerski, and P~Szymczak.
\newblock Rotne--prager--yamakawa approximation for different-sized particles
  in application to macromolecular bead models.
\newblock {\em Journal of Fluid Mechanics}, 741, 2014.

\end{thebibliography}

\end{document}